\documentclass[12pt]{article}
\usepackage{setspace}
\usepackage{mathtools}
\usepackage[utf8]{inputenc}
\usepackage{url}

\usepackage[square,numbers]{natbib}
\setcitestyle{authoryear,open={(},close={)}}
\usepackage[dvipsnames]{xcolor}        

\RequirePackage[colorlinks,citecolor=blue,urlcolor=blue]{hyperref}
\newcommand*{\doi}[1]{\href{http://dx.doi.org/#1}{#1}}
\usepackage{mathtools}
\usepackage{comment}

\usepackage[psamsfonts]{amssymb}
\usepackage{amsthm,amsmath,amsbsy,mathabx}
\usepackage{bbm, bm}
\usepackage[dvipsnames]{xcolor}
\usepackage{subcaption}
\usepackage[shortlabels]{enumitem}
\RequirePackage[colorlinks,citecolor=blue,urlcolor=blue]{hyperref}
\usepackage{graphicx,placeins}

\usepackage{inputenc}
\usepackage{graphicx,graphics,multirow}
\usepackage{subcaption}
\usepackage{supertabular}
\usepackage{amsmath,amscd,amsfonts,amssymb}
\usepackage{fancyhdr}
\usepackage{color}
\usepackage{verbatim}

\usepackage{pdflscape}
\usepackage{array,supertabular}
\usepackage{float}
\usepackage{rotating}
\usepackage{lscape}

\usepackage{booktabs,caption}
\usepackage[flushleft]{threeparttable}
\usepackage{amsthm}

\usepackage{dcolumn}
\newcolumntype{L}{D{.}{.}{2,5}}

\usepackage{chngpage} 

\newtheorem{theorem}{Theorem}[section]
\newtheorem{lemma}[theorem]{Lemma}

\newtheorem{proposition}[theorem]{Proposition}
\newtheorem{corollary}[theorem]{Corollary}

\newtheorem{assumption}{Assumption}
\newtheorem{remark}{Remark}[section]

\def\mds{\medskip}
\def\Rb{{\mathbb R}}

\def\Pc{{\mathcal P}}
\def\Fc{{\mathcal F}}
\def\Nc{{\mathcal N}}

\long\def\symbolfootnote[#1]#2{\begingroup%
\def\thefootnote{\fnsymbol{footnote}}\footnotetext[#1]{#2}\footnotemark[#1]\endgroup}

\def\ba#1\ea{\begin{align*}#1\end{align*}} 
\def\banum#1\eanum{\begin{align}#1\end{align}} 

\makeatletter
\@addtoreset{equation}{section}

\makeatother

\newcounter{Fig}[figure]

\newcounter{Tab}[table]

  {\refstepcounter{Tab}
   \addtocounter{table}{1}
   \samepage\vspace{0.2cm}
   \centerline{#1} \nobreak
   \begin{verse}{Table \thesection.\arabic{Tab}:}
  }%
  {\end{verse} \vspace{0.2cm}}

\relax

\newcommand{\dd}{\mbox{\boldmath $d$}}

\newcommand{\vv}{\mbox{\boldmath $v$}}

\def \Eb{{\mathbb E}}

\def \Lb{{\mathbb L}}

\def \Pb{{\mathbb P}}
\def \Rb{{\mathbb R}}

\def \Zb{{\mathbb Z}}

\def \Ac{{\mathcal A}}
\def \Rc{{\mathcal R}}

\def \Uc{{\mathcal U}}

\def \Fc{{\mathcal F}}

\def \Sc{{\mathcal S}}

\def \Nc{{\mathcal N}}
\def \Pc{{\mathcal P}}

\newcommand{\bqa}{\begin{eqnarray*}}
\newcommand{\eqa}{\end{eqnarray*}}
\newcommand{\bqan}{\begin{eqnarray}}
\newcommand{\eqan}{\end{eqnarray}}
\newcommand{\bqt}{\begin{quote}}
\newcommand{\eqt}{\end{quote}}
\newcommand{\bt}{\begin{tabbing}}
\newcommand{\et}{\end{tabbing}}
\newcommand{\bit}{\begin{itemize}}
\newcommand{\eit}{\end{itemize}}
\newcommand{\ben}{\begin{enumerate}}
\newcommand{\een}{\end{enumerate}}
\newcommand{\beq}{\begin{equation}}
\newcommand{\eeq}{\end{equation}}
\newcommand{\beqw}{\begin{equation*}}
\newcommand{\eeqw}{\end{equation*}}
\newcommand{\bdefi}{\begin{definition}}
\newcommand{\edefi}{\end{definition}}
\newcommand{\bpro}{\begin{proposition}}
\newcommand{\epro}{\end{proposition}}
\newcommand{\blem}{\begin{lemma}}
\newcommand{\elem}{\end{lemma}}
\newcommand{\bco}{\begin{corollary}}
\newcommand{\eco}{\end{corollary}}
\newcommand{\bdes}{\begin{description}}
\newcommand{\edes}{\end{description}}
\newcommand{\bre}{\begin{remark}}
\newcommand{\ere}{\end{remark}}

\newcommand{\eps}{\epsilon}

\usepackage{authblk}

\def\mds{\medskip}
\def\1{{\mathbf 1}}
\usepackage{xr}
\makeatletter
\newcommand*{\addFileDependency}[1]{
  \typeout{(#1)}
  \@addtofilelist{#1}
  \IfFileExists{#1}{}{\typeout{No file #1.}}
}
\makeatother

\makeatletter
\def\@seccntformat#1{\@ifundefined{#1@cntformat}%
   {\csname the#1\endcsname\quad}
   {\csname #1@cntformat\endcsname}
}
\makeatother

\begin{document}

\title{{\Large \bf Sparse factor models of high dimension}}

\medskip


\author{Benjamin Poignard\footnote{Osaka University, Graduate School of Economics, 1-7 Machikaneyama, Toyonaka, Osaka 560-0043, Japan; jointly affiliated at Riken AIP. bpoignard@econ.osaka-u.ac.jp.}, Yoshikazu Terada\footnote{Osaka University, Graduate School of Engineering Science, 1-3 Machikaneyama, Toyonaka, Osaka 560-8531, Japan; jointly affiliated at Riken AIP. terada.yoshikazu.es@osaka-u.ac.jp}}

\maketitle

\begin{abstract}
We consider the estimation of a sparse factor model where the factor loading matrix is assumed sparse. The estimation problem is reformulated as a penalized M-estimation criterion, while the restrictions for identifying the factor loading matrix accommodate a wide range of sparsity patterns. We prove the sparsistency property of the penalized estimator when the number of parameters is diverging, that is the consistency of the estimator and the recovery of the true zeros entries. 
These theoretical results are illustrated by finite-sample simulation experiments, and the relevance of the proposed method is assessed by  applications to portfolio allocation and macroeconomic data prediction.

\medskip


\medskip


\medskip
\noindent

\end{abstract}

\section{Introduction}\label{sec:1}

A major difficulty in the multivariate modeling of high-dimensional random vectors is the trade-off between a sufficiently rich parameterized model to capture complex relationships and yet parsimonious enough to avoid over-fitting issues. In light of this trade-off between parsimony and flexibility, factor modeling is a pertinent solution: it aims to summarize the information from a large data set of dimension $p$ through a small number $m$ of variables called factors. 
The quantity of interest in factor analysis is the variance-covariance matrix $\Sigma^\ast$ of the vector of observations, which is specified as $\Sigma^\ast = \Lambda^\ast M_{FF} \Lambda^{\ast\top} +\Psi^\ast$, with the factor loading matrix $\Lambda^\ast$, the variance-covariance of the factors $M_{FF}$ and the variance-covariance of the idiosyncratic errors $\Psi^\ast$. 
Numerous restrictions can be imposed for the sake of identification. Most works devoted to factor models constrain $\Psi^\ast$ to be diagonal when the dimension $p$ - the number of variables composing the vector of observations - is fixed. Under this constraint, \cite{anderson1988} derived the large sample properties of the Gaussian-based maximum likelihood (ML) estimator of the factor model. \cite{bai2012} extended these asymptotic results when $p$ is diverging and proposed different restrictions on $\Lambda^\ast$ and $M_{FF}$.
To relax the diagonal constraint on $\Psi^\ast$, \cite{chamberlain1983} developed the notion of approximate factor models, which allows for cross-correlation among the idiosyncratic errors. Under the condition of bounded eigenvalues of $\Psi^\ast$ non-diagonal, \cite{bai2016} derived the large sample properties of the Gaussian-based ML factor model estimator when $p$ diverges. 
The sparse estimation of a non-diagonal $\Psi^\ast$ was considered by \cite{bailiao2016}. Their approach gave rise to the notion of conditionally sparse factor models, in the sense that $\Psi^\ast$ is a sparse matrix with bounded eigenvalues. Under orthogonal strong factors, they computed the sparse approximate factor estimator, where $\Lambda^\ast$ and $\Psi^\ast$ are jointly estimated while penalizing $\Psi^\ast$ only through the adaptive Least Absolute Shrinkage and Selection Operator (LASSO) and Smoothly Clipped Absolute Deviation (SCAD) methods, and derived some consistency results. \cite{poignard2020} considered a two-step method for estimating a sparse $\Psi^\ast$ and derived consistency results together with the conditions for support recovery. 

\mds

\begin{figure}[htbp]
  \begin{minipage}[b]{0.45\textwidth}
    \centering
    \includegraphics[height = 42mm]{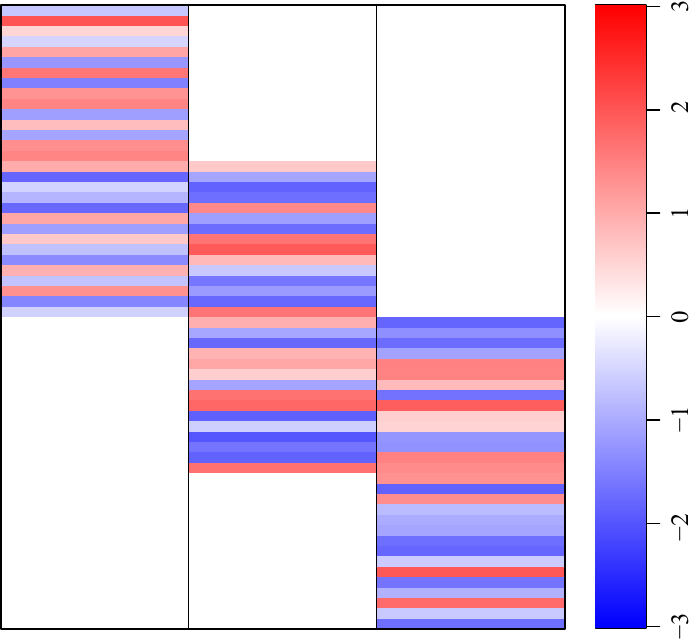}
    \subcaption{True loading matrix}
    \label{fig:(a)}
  \end{minipage}
  \begin{minipage}[b]{0.45\textwidth}
    \centering
    \includegraphics[height = 42mm]{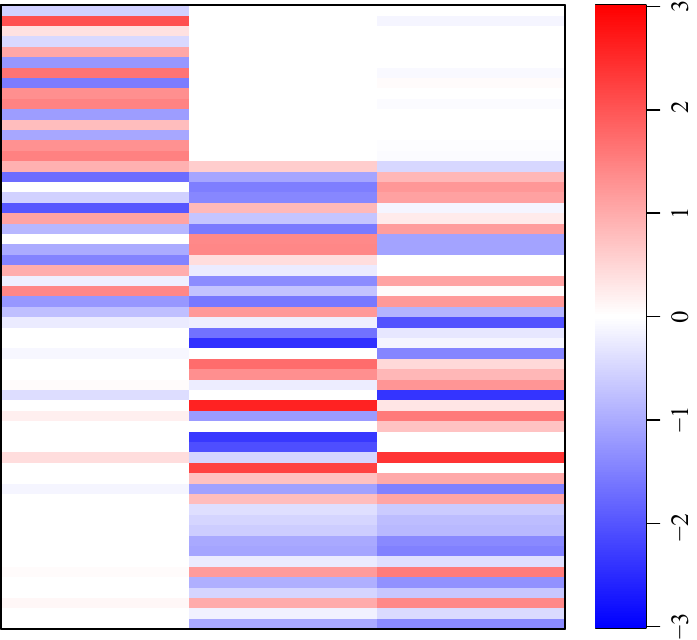}
    \subcaption{SOFAR estimator}
    \label{fig:(b)}
  \end{minipage}
  \caption{True loading matrix versus SOFAR estimator, with $p= 60$, the number of factors $m=3$ and sample size $n = 1000$.} 
\end{figure}

On the other hand, the studies on sparse factor models, where many
of the entries of $\Lambda^\ast$ are exactly equal to zero, have benefited from a limited attention so far. The key identification issue relating to $\Lambda^\ast$ is the rotational indeterminacy, where individual factors are only identified up to a rotation, or equivalently, different sparse column orderings in $\Lambda^\ast$ may provide an optimal solution, prohibiting any economic interpretation of the estimated factors. This identification issue has hampered the sparse modeling of $\Lambda^\ast$ and rendered its interpretation in terms of sparsity intricate. 
In a recent work, \cite{uematsu2023a} developed a Principal Component Analysis (PCA)-based sparse estimation procedure for $\Lambda^\ast$. Within the weak factor model setting, they specified a column-wise sparse assumption for $\Lambda^\ast$, where sparsity among the $k$-th column grows as $p^{\nu_k}$ for some unknown deterministic value $0<\nu_k\leq 1$. 
They proposed the Sparse Orthogonal Factor Regression (SOFAR) estimator, which corresponds to an adaptive LASSO-penalized least squares problem providing both the estimated factor variables and sparse factor loading matrix, but not the estimator of $\Psi^\ast$ due to the PCA nature of the framework, and showed its consistency and support recovery. 
Their analysis of the SOFAR is developed under the assumption of orthogonal factors and the restriction ``$\Lambda^{\ast\top}\Lambda^\ast$ diagonal''. The latter condition may greatly worsen the performances of the SOFAR method once $\Lambda^\ast$ does not satisfy this diagonal restriction. This constraint holds when $\Lambda^\ast$ satisfies a perfect simple structure, a case where each row of $\Lambda^\ast$ has at most one non-zero element, giving rise to a non-overlapping clustering of variables. 
To illustrate this point, Figure~\ref{fig:(a)} displays the true sparse loading matrix $\Lambda^\ast$ that satisfies a near perfect simple structure, thus violating the constraint ``$\Lambda^{\ast\top}\Lambda^\ast$ diagonal'' assumed by the SOFAR method. Based on the data generated from a centered Gaussian distribution with variance-covariance $\Sigma^\ast = \Lambda^\ast\Lambda^{\ast\top}+\Psi^\ast$, $\Psi^\ast$ diagonal with positive non-zero components, we applied the SOFAR procedure: Figure~\ref{fig:(b)} displays the loading matrix estimated by the SOFAR and highlights the substantial discrepancy between the SOFAR estimator and the true sparse loading matrix. The details relating to the data generating process are described in Section~\ref{simulations}. 
Following the work of \cite{bailiao2016}, \cite{daniele2019} considered the joint estimation of $(\Lambda^\ast,\Psi^\ast)$ based on a Gaussian quasi-ML with a LASSO penalization on $\Lambda^\ast$ and derived some consistency results. 
\cite{freyaldenhoven2023} proposed a $\ell_1$-rotation criterion and provided some theoretical properties. However, this is a rotation technique only; thus, there is no guarantee that this approach provides a sparse loading matrix. \cite{hirose2023} restricted their analysis to the perfect simple structure and proposed the prenet penalty function to suitably recover such a structure.

Another line of research relates to sparse Bayesian factor models. Recent theoretical studies essentially derived from the strong factor framework: \cite{pati2014} studied posterior contraction rates for variance-covariance estimation under a particular class of continuous shrinkage priors for $\Lambda^\ast$; \cite{ohn2024} derived an inferential framework to learn the factor dimensionality and the sparse structure. The applications of the Bayesian approach typically concern gene expression data, as in \cite{west2003}. 

In this paper, we aim to tackle the issues previously mentioned for the estimation of potentially any arbitrary sparse $\Lambda^\ast$. We consider the following problem: given $n$ observations of a $p$-dimensional random vector $X_t$, we jointly estimate $(\Lambda^\ast,\Psi^\ast)$ while penalizing $\Lambda^\ast$. We first establish the consistency of the factor decomposition-based estimator of $\Sigma^\ast$ for any arbitrary sparse $\Lambda^\ast$. Then, under a suitable identification condition, we derive the consistency of $\Lambda^\ast$ and $\Psi^\ast$, and prove the recovery of the true zero entries of $\Lambda^\ast$. Importantly, our identification restriction on $\Lambda^\ast$ allows us to consider the sparsest $\Lambda^\ast$ that does not necessarily satisfy the restrictions of \cite{bai2012} or \cite{uematsu2023a}: the latter restrictions may imply sparse loading matrices but may also restrict the diversity of their sparsity patterns. Since the seminal work of \cite{fan2001}, a significant literature on penalized M-estimators has been flourishing: see, e.g., \cite{fan2004} for an asymptotic analysis; \cite{loh2015} or \cite{poignard2022} for a non-asymptotic viewpoint. In particular, it has been applied to the estimation of variance-covariance and precision matrices of high dimension: e.g., \cite{rothman2008} derived consistency results for LASSO-penalized precision matrices; \cite{lam2009} established the ``sparsistency'' property of sparse variance-covariance, precision and correlation matrices based on the Gaussian ML, where sparsistency refers to the consistent estimation and the consistent recovery of all the true zero coefficients with probability one. Our proposed estimator for sparse factor models builds upon this framework. Our main contributions are as follows: we develop an estimation procedure for the estimation of sparse factor loading matrix; we prove the consistency and the correct recovery of the true zero entries when the dimension $p$ diverges; we provide an implementation procedure. 

\mds

The remainder of the paper is organized as follows. Section \ref{framework} outlines the sparse factor loading framework. Section \ref{asymptotic_theory} is devoted to the asymptotic properties of the penalized factor model estimator. Sections \ref{simulations} and \ref{real_data} illustrate these theoretical results through simulated and real data experiments. All the proofs of the main text, auxiliary results and implementation details are relegated to the Appendices. 

\mds

\textbf{\emph{Notations.}} Throughout this paper, we denote the cardinality of a set $E$ by $|E|$. For $\vv \in \Rb^{d}$, the $\ell_p$ norm is $\|\vv\|_p = \big(\sum^{\text{d}}_{k=1} |\vv_k|^p \big)^{1/p}$ for $p > 0$, and $\|\vv\|_{\infty} = \max_i|\vv_i|$. For a symmetric matrix $A$, $\lambda_{\min}(A)$ (resp. $\lambda_{\max}(A)$) is the minimum (resp. maximum) eigenvalue of $A$, and $\text{tr}(A)$ is the trace operator. For a matrix $B$, $\|B\|_s = \lambda^{1/2}_{\max}(B^\top B)$ and $\|B\|_F=\text{tr}^{1/2}(B^\top B)$ are the spectral and Frobenius norms, respectively, and $\|B\|_{\max} = \max_{ij}|B_{ij}|$ is the coordinate-wise maximum. For $B \in \Rb^{p \times m}$, the $\ell_1$ matrix norm is $\|B\|_1 = \max_{1 \leq j \leq m}\sum^p_{i=1}|B_{ij}|$.
We denote by $O \in \Rb^{d \times p}$ the $d \times p$ zero matrix. We denote by  $\mathcal{O}(m)$ the set of all $m\times m$ orthogonal matrices. 
$\text{vec}(A) \in \Rb^{dp}$ is the vectorization operator that stacks the columns of the matrix $A \in \Rb^{d \times p}$ on top of one another into a vector. 
For a square matrix $A \in \Rb^{d \times d}$, $\text{diag}(A)$ is the operator that stacks the diagonal elements of $A$ on top of one another into a vector.
The matrix $I_d$ denotes the $d$-dimensional identity matrix. 
For two matrices $A,B$, $A \otimes B$ is the Kronecker product. For two matrices $A,B$ of the same dimension, $A \odot B$ is the Hadamard product. For a function $f: \Rb^{d} \rightarrow \Rb$, we denote by $\nabla f$ the gradient or subgradient of $f$ and $\nabla^2 f$ the Hessian of $f$. We denote by $\partial f$ the component-by-component partial derivative of $f$.
$\Sc^c$ denotes the complement of the set $\Sc$. 

\section{The framework} \label{framework}

We consider a sequence of a $p_n$-dimensional random vector $X_t, t = 1,\ldots,n$ following the factor decomposition:
$$
\forall t = 1,\ldots,n, \;\; X_t = \Lambda^\ast_n F_t +\epsilon_t,
$$
where $F_t \in \Rb^m$ is the vector of factors, $\Lambda^\ast_n = (\lambda^\ast_1,\ldots,\lambda^\ast_{p_n})^\top \in \Rb^{p_n \times m}$ is the loading matrix, with $\lambda^\ast_j = (\lambda^\ast_{j1},\ldots,\lambda^\ast_{jm})^\top \in \Rb^m, j=1,\ldots,p_n$, and $\eps_t \in \Rb^{p_n}$ is the vector of idiosyncratic errors. Both $F_t,\eps_t$ are not observable; $X_t$ only is observable. Throughout the paper, the number of factors $m$ remains fixed; the sequence $(p_n)$ could tend to the infinity with $n$, i.e., $p_n\rightarrow \infty$ as $n\rightarrow \infty$. This motivates the indexing of the parameters by $n$. We make the following assumptions concerning the random variables $(F_t), (\eps_t)$.
\begin{assumption}\label{factor_assumption_1}
\begin{itemize}
    \item[(i)] $(F_t,\eps_t)_{t \geq 1}$ is stationary and ergodic with
    $\forall 1 \leq j \leq p_n, \; \Eb[\eps_{t,j}]=0$, and $\Eb[\eps_t\eps^\top_t]=\Psi^\ast_n \in \Rb^{p_n \times p_n}$ diagonal
    with positive diagonal elements $(\sigma^{\ast 2}_{1},\ldots,\sigma^{\ast 2}_{p_n})$, and
    $\forall 1 \leq j \leq m, \; \Eb[F_{t,j}] = 0$, $\Eb[F_tF^\top_t]=M_{FF} \in \Rb^{m \times m}$, and $\Eb[\eps_t F^\top_t]= O \in \Rb^{p_n \times m}$.
    \item[(ii)] $\exists r_1,r_2>0, b_1,b_2>0$ such that for any $s>0$: 
    $$\forall 1 \leq j \leq p_n, \; \Pb\big(|\eps_{t,j}| > s\big) \leq \exp(-(s/b_1)^{r_1}),\; \forall 1 \leq j \leq m,\; \Pb\big(|F_{t,j}| > s\big) \leq \exp(-(s/b_2)^{r_2}).$$
\end{itemize}
\end{assumption}

\begin{assumption}\label{factor_assumption_2}
$\exists r_3>0$ and $L>0$ satisfying, for all $n\in \mathbb{Z}^+$: $\alpha(n) \leq \exp(-Ln^{r_3})$, with $\alpha(\cdot)$ the mixing coefficient: $\alpha(n) = \sup_{A \in \Fc^0_{-\infty}, B\in \Fc^{\infty}_n } |\Pb(A)\Pb(B)-\Pb(A \cup B)|$, where $\Fc^0_{-\infty},\Fc^{\infty}_n$ the filtrations generated by $\{(F_t,\eps_t): - \infty \leq t \leq 0\}$ and $\{(F_t,\eps_t): n \leq i \leq \infty\}$. 
\end{assumption}

Assumptions \ref{factor_assumption_1}--\ref{factor_assumption_2} relate to the probabilistic properties of the factors and idiosyncratic errors. Both factors and idiosyncratic errors are allowed to be weakly dependent and to satisfy the strong mixing condition and exponential tail bounds. In particular, in Assumption \ref{factor_assumption_1}, we assume that the elements of $\eps_t$ are uncorrelated. This will facilitate the analysis under high-dimension. This restriction is also stated in Assumption B of \cite{bai2012}. 

\mds

From Assumptions \ref{factor_assumption_1}, $\Eb[X_t] = 0$ and the implied variance-covariance $\Sigma^\ast_n:=\text{Var}(X_t)$ becomes $\Sigma^\ast_n= \Lambda^\ast_n M_{FF} \Lambda^{\ast\top}_n + \Psi^\ast_n$. 
For the sake of interpreting sparsity which will be our core motivation, we assume the orthogonal factor model, i.e., $M_{FF}=I_m$.
This is a standard restriction to avoid some of the indeterminacy issues in factor models: see, e.g., \cite{anderson2003introduction}.
If we consider the oblique factor model, 
we may desire sparsity even for the factor covariance $M_{FF}$, thus complicating the problem formulation and the theoretical analysis. In light of the diverging nature of the parameters, the problem consists of a sequence of parametric factor models
$\Pc_n:=\{\Pb_{\theta_n}, \,\theta_n \in \Theta_n :=\Theta_{n\Lambda}\times \Theta_{n\Psi}\}$, with $\Theta_n\subseteq \Rb^{p_nm} \times ]0,\infty[^{p_n}$, $\theta_n = (\theta^\top_{n\Lambda},\theta^\top_{n\Psi})^\top$, with $\theta_{n\Lambda} = \text{vec}(\Lambda_n),\; \theta_{n\Psi}=\text{diag}(\Psi_n)=(\sigma^2_1,\ldots,\sigma^2_{p_n})^\top$. Therefore, the number of unknown parameters $p_n(m+1)$ may vary with the sample size. We assume that $\Pc_n$ contains the ``pseudo-true'' factor model parameter $\theta^\ast_n = (\theta^{\ast\top}_{n\Lambda},\theta^{\ast\top}_{n\Psi})^\top$. 

\mds

Because of the rotational indeterminacy, even under the orthogonal factor model, the model can be identified up to a rotation matrix, and at least $m^2$ restrictions are required for the sake of identification. Indeed, for any $R \in \mathcal{O}(m)$, $\Lambda^\ast_n\Lambda^{\ast\top}_n = \widetilde{\Lambda}^\ast_n\widetilde{\Lambda}^{\ast\top}_n$ with $\widetilde{\Lambda}^\ast_n:=\Lambda^\ast_n R$. In other words, we can rotate $\Lambda^\ast_n$ without altering the implied variance-covariance structure. Restrictions have thus been proposed to uniquely fix $\Lambda^\ast_n$: see Section 5 of \cite{anderson1956}; Section 4 of \cite{bai2012}. To prove the consistent estimation of the factor model parameters and the sparsistency, we will rely on \cite{anderson1956}'s condition which will be formulated in $\Theta_{n\Lambda}$.

\mds

To estimate $\theta^\ast_n$ under sparsity, we consider a global loss $\Lb_n:\Rb^{n\,p_n}\times\Theta_n \rightarrow \Rb$, where $\Lb_n(\theta_n)$ is evaluated under the parametric factor model $\Pc_n$, and we assume that $\Lb_n(\theta_n)$ is the empirical loss associated to a continuous function $(\ell_n)_{n\geq 1}$ with $\ell_n : \Rb^{p_n} \times \Theta_n \rightarrow \Rb$ and that we can write as $\Lb_n(\theta)=\sum_{t=1}^n\ell_n(X_t;\theta_n)$. We will consider two functions: the Gaussian loss and the least squares loss, respectively defined as $\ell_n(X_t;\theta) = \text{tr}\big(X_tX^\top_t \Sigma^{-1}_n\big) + \log(|\Sigma_n|)$ and $\ell_n(X_t;\theta_n) = \text{tr}\big((X_tX^\top_t - \Sigma_n)^2\big)$. The sparse factor model estimation problem is
\begin{equation}\label{stat_crit}
\widehat{\theta}_n \in \underset{(\theta_{n\Lambda},\theta_{n\Psi}) \in \Theta_{n\Lambda} \times \Theta_{n\Psi}}{\arg\,\min} \; \Lb_n(\theta_n) + n\, \sum_{k=1}^{p_nm}p(|\theta_{n\Lambda,k}|,\gamma_n),
\end{equation}
when such a minimizer exists. $p(x,\gamma_n), x \geq 0$ is a penalty function and $\gamma_n$ is the regularization parameter. 
We will consider the coordinate-separable non-convex SCAD and Minimax Concave Penalty (MCP) penalty functions. The SCAD of \cite{fan2001}, for $a_{\text{scad}}>2$, is defined as: for every $x \geq 0$,
\begin{equation*}
p(x,\gamma) = \gamma x \, \mathbf{1}\big(x \leq \gamma\big) + \frac{2 a_{\text{scad}}\gamma x-x^2-\gamma^2}{2(a_{\text{scad}}-1)}\, \mathbf{1}\big(\gamma < x \leq a_{\text{scad}} \gamma\big) + \frac{1}{2}(a_{\text{scad}}+1)\gamma^2 \, \mathbf{1}\big(x > a_{\text{scad}}\gamma\big),
\end{equation*}
where $a_{\text{scad}}>2$. The MCP due to~\cite{zhang2010} is defined for $b_{\text{mcp}}>0$ as: \textcolor{black}{for every $x\geq 0$}, 
\begin{equation*}
p(x,\gamma)  = \big(\gamma x-\frac{x^2}{2 b_{\text{mcp}}}\big) \, \mathbf{1}\big(x\leq b_{\text{mcp}}\gamma\big) + \gamma^2\frac{b_{\text{mcp}}}{2} \, \mathbf{1}
\big(x> b_{\text{mcp}}\gamma\big).
\end{equation*}
Throughout the paper, we will work under the following assumption, where we require the uniqueness of $\theta^\ast_n$ only.
\begin{assumption}\label{parameter_assumption}
$\forall n$, $\Theta_n$ is a borelian subset of $\Rb^{p_n(m+1)}$. The function $\theta_n \mapsto \Eb[\ell_n(X_t;\theta_n)]$ is uniquely minimized on $\Theta_n$ at $\theta^\ast_n$, and an open neighborhood of $\theta^\ast_n$ is contained in $\Theta_n$.
\end{assumption}

\section{Asymptotic properties}\label{asymptotic_theory}

To prove the large sample properties, we make the following assumptions. 
\begin{assumption}\label{assumption_parameters}
There exists $0<\kappa<\infty$ large enough such that the true parameter $\theta^{\ast}_{n}$ satisfies: $\text{\rm (i) } \forall 1 \leq j \leq p_n, \|\lambda^\ast_{j}\|_2\leq \kappa, \;\text{ and }\;
\text{\rm (ii) } \forall 1 \leq k \leq p_n, \kappa^{-2} \leq \sigma^{\ast2}_{k} \leq \kappa^2$.
\end{assumption}

\begin{assumption}\label{assumption_sparse_lambda}
For every $n$, $\Lambda^\ast_n \in \Rb^{p_n \times m}$, with $m$ fixed, is $s_n$-sparse, that is $s_n = \text{card}(\Sc_n)$ with
$\Sc_n:=\big\{k: \theta_{n\Lambda,k}^\ast\neq 0, k = 1,\ldots,m\,p_n\big\}$.
\end{assumption}

\begin{assumption}\label{assumption_regularity_loss_n}
The map $\theta_n\mapsto \ell_n(X_t;\theta_n)$ is twice differentiable on $\Theta_n$, for every $X_t \in \Rb^{p_n}$. Any pseudo-true value $\theta^\ast_{n}$ satisfies  
the first-order condition $\Eb[ \nabla_{\theta_n}\ell_n(X_t;\theta^\ast_{n})] = 0$.
\end{assumption}

We will denote by $\partial_1 p(x,\gamma_n)$ (resp. $\partial^2_{11}p(x,\gamma_n)$) the first order (resp. second order) derivative of $x\mapsto p(x,\gamma_n)$, for any $\gamma_n$.

\begin{assumption}\label{assumption_regularity_penalty_n}
Define 
$$A_n = \underset{1 \leq k \leq p_n m}{\max}\{|\partial_{1}p(|\theta^\ast_{n\Lambda,k}|,\gamma_n)|,\theta^\ast_{n\Lambda,k}\neq 0\}, B_n = \underset{1 \leq k \leq p_n m}{\max}\{|\partial^2_{11} p(|\theta^\ast_{n\Lambda,k}|,\gamma_n)|,\theta^\ast_{n\Lambda,k}\neq 0\}.$$
Then we assume $\sqrt{p_n}A_n \rightarrow 0$ and $B_n \rightarrow 0$ as $n \rightarrow \infty$. Moreover, $\exists\overline{K}_1, \overline{K}_2$ finite constants such that $|\partial^2_{11} p(\theta_1,\gamma_n) - \partial^2_{11} q(\theta_2,\gamma_n)| \leq \overline{K}_2 |\theta_1-\theta_2|$ for any real $\theta_1,\theta_2$ such that $\theta_1,\theta_2>\overline{K}_1\gamma_n$.
\end{assumption}
Assumption \ref{assumption_parameters} ensures that all the eigenvalues of the factor decomposition-based $\Sigma^\ast_n = \Lambda^\ast_n\Lambda^{\ast\top}_n+\Psi^\ast_n$ are bounded away from $0$ and that $\|\Sigma^{\ast-1}_n\|_s=O(1)$. 
The sparsity condition specified in Assumption \ref{assumption_sparse_lambda} implies $\|\Lambda^{\ast\top}_n\Lambda^\ast_n\|_s\leq \|\Lambda^\ast_n\|^2_F=s_n$, $\|\Lambda^\ast_n\|_1 \leq \sqrt{p_n}\|\Lambda^\ast_n\|_F=O(\sqrt{p_ns_n})$ and $\|\Lambda^\ast_n\|_1 \geq \|\Lambda^\ast_n\|_s=O(\sqrt{s_n})$. The strong factor model condition sets $s_n=p_n$: see, e.g., Assumption 3.5 in \cite{fan2011} or Assumption 3.3 in \cite{bailiao2016}. 
In contrast, the weak factor model viewpoint entails $s_n = p^{\nu}_n$ with $ \nu \leq 1$:  in Assumption A2 in \cite{bai2023}, $\nu \in ]0,1]$, whereas Assumption 1 in \cite{daniele2019} sets $\nu \in ]1/2,1]$; the latter work uses the QML framework and thus requires $\nu > 1/2$ for consistent estimation, whereas that work builds upon the PCA framework. 
Our approach for the consistent estimation of the factor model-based $\Sigma^\ast_n$ in Theorem \ref{Theorem_existence_consistent} and for the sparsistency property in Theorem \ref{sparsistency} under \cite{anderson1956}'s condition for identification below does not require any specific behavior for $s_n$ other than $s_n = O(p_n)$. However, to prove the consistency of the GLS estimator of the latent factor variables, we will work in the matrix space $\{\Lambda_n \in \Rb^{p_n \times m}\mid \underline{\mu}\leq \lambda_{\min}(p^{-\nu}_n\Lambda^\top_n\Lambda_n)\leq \lambda_{\max}(p^{-\nu}_n\Lambda^\top_n\Lambda_n) \leq \overline{\mu}\}$ with $\underline{\mu},\overline{\mu}>0$ and $\nu\in ]1/2,1]$, a space that includes the weak factor case. This will be helpful to obtain $\|(\Lambda^\top_n\Psi^{-1}_n \Lambda_n)^{-1}\|_s=O_p(p^{-\nu}_n)$.
Assumption \ref{assumption_regularity_loss_n} concerns the regularity of the loss function. Assumption \ref{assumption_regularity_penalty_n} relates to regularity of the penalty function and tuning parameter.



\mds

The following result establishes the existence of a consistent estimator of the factor model-based variance-covariance.

\begin{theorem}\label{Theorem_existence_consistent}
Suppose $3 r^{-1}_1 + 1.5r^{-1}_2+r^{-1}_3>1$. Let $\rho^{-1} = 3 r^{-1}_1 + 1.5r^{-1}_2+r^{-1}_3+1$. Assume $\log(p_n)^{6/\rho}=o(n)$ holds and Assumptions \ref{factor_assumption_1}-\ref{assumption_regularity_penalty_n} are satisfied. Then there exists a sequence of estimators $\widehat{\Sigma}_n=\widehat{\Lambda}_n\widehat{\Lambda}^\top_n+\widehat{\Psi}_n$ as defined in (\ref{stat_crit}) satisfying
\[
\|\widehat{\Sigma}_n-\Sigma^\ast_n\|_F = O_p\left(p_n\sqrt{\frac{s_n\log(p_n)}{n}}+\sqrt{p_n}A_n\right).
\]  
\end{theorem}

\mds


To ensure the uniqueness of the factor decomposition for the true covariance matrix, which will allow us to deduce the consistency of $\widehat{\theta}_n$ and the sparsistency property, we rely on the following identification condition.
\begin{assumption}\label{assumption_AR_condition}
\citet{anderson1956}'s condition: For $j = 1,\dots, p_n$, 
let $\Lambda_{n,-j}^\ast$ be the matrix obtained by removing the $j$-th row from $\Lambda_n^\ast$.
Assume that there exists a positive constant $c_L>0$ such that, for each $j=1,\dots,p_n$,  
$\Lambda_{n,-j}^\ast$ contains two submatrices whose determinants are bounded below by $c_L$.    
\end{assumption}
The factor loading matrix is not uniquely determined due to the rotational indeterminacy. 
Here, we assume that the sparsest true loading matrix $\Lambda_n^\ast$ satisfies the condition for identification introduced by \citet{anderson1956}, for which the upper square matrix of the sparsest loading matrix $\Lambda_n^\ast$ is lower triangular, i.e., denoting $\Lambda_{l,:} \in \Rb^{l \times m}$:
\begin{equation}
\Lambda_n= 
\begin{pmatrix}
\Lambda_{m,:}\\
\Lambda_{p_n-m,:}
\end{pmatrix}\; \text{ with } \; 
\Lambda_{m,:} = \begin{pmatrix} 
\lambda_{11} & 0 & \cdots & 0 \\
\lambda_{21} & \lambda_{22} & \cdots & 0 \\
\vdots & \vdots & \ddots & \vdots \\
\lambda_{m1} & \lambda_{m2} & \cdots & \lambda_{mm}
\end{pmatrix}
\; \text{ and } \; 
\lambda_{ii} > 0 \; (i=1,\dots,m).
\label{eq:IC5}
\end{equation}
The positivity of the diagonal elements is needed to avoid the sign indeterminacy.
Since the loading matrix in Figure~\ref{fig:(a)} does not satisfy this condition,
one might think that the condition is not suitable for sparse estimation.
With a certain variable ordering, even a perfect simple structure may not satisfy 
$\lambda_{jj}> 0$ for some $j$.
However, it is important to note that this identification condition can be slightly relaxed.
In fact, if there exists a $p_n\times p_n$ permutation matrix $P_n$ such that the permuted loading matrix 
$\tilde{\Lambda}_n^\ast = P_n\Lambda_n^\ast$ satisfies condition~(\ref{eq:IC5}), 
then we can eliminate the rotational indeterminacy.
Therefore, we employ the following modified identification condition throughout the paper. 

\begin{assumption}\label{assumption_space_lambda}
The parameter space of $\Lambda^\ast_n$ is 
\[
\Theta_{n\Lambda}=\{\Lambda_n \in \Rb^{p_n\times m} \mid \exists P_n \in \Pi(p_n);\; P_n\Lambda_n \text{ satisfies condition~(\ref{eq:IC5})} \},
\]
where $\Pi(p_n)$ is the set of all $p_n\times p_n$ permutation matrices.
\end{assumption}

We now can establish the consistency of $\widehat{\Lambda}_n$ and $\widehat{\Psi}_n$ under the conditions of Assumptions \ref{assumption_AR_condition} and \ref{assumption_space_lambda} and the recovery of the zero entries of $\Lambda^\ast_n$ - ``sparsity property'' -. We also prove the consistent estimation of the generalized least squares (GLS) estimator $\widehat{F}_t=(\widehat{\Lambda}^\top_n\widehat{\Psi}^{-1}_n \widehat{\Lambda}_n)^{-1} 
\widehat{\Lambda}^\top_n \widehat{\Psi}^{-1}_n X_t$ of the latent variable $F_t$, for any $t=1,\ldots,n$. See Section 14.7 of \cite{anderson2003introduction} for its derivation.

\begin{theorem}\label{sparsistency}
In addition to the conditions of Theorem \ref{Theorem_existence_consistent}, 
suppose Assumptions \ref{assumption_AR_condition} and \ref{assumption_space_lambda} are satisfied. 
Let $\widehat{\Lambda}_n$ be the estimator for $\Lambda_n^\ast$, which ignores the sign indeterminacy as usual.
The estimator  $\widehat{\theta}_n = (\widehat{\theta}^\top_{n\Lambda},\widehat{\theta}^\top_{n\Psi})^\top$ defined in (\ref{stat_crit}) satisfies:
\begin{equation*}
\|\widehat{\theta}_n-\theta^\ast_n\|_\infty = O_p\left(p_n\sqrt{\frac{s_n\log(p_n)}{n}}+\sqrt{p_n}A_n\right).
\end{equation*}
Assume that the penalty function satisfies $\gamma^{-1}_n\lim_{n \rightarrow \infty}\lim_{x\rightarrow 0}\partial_{x}p(x,\gamma_n) >0$. Moreover, assume $\gamma_n \rightarrow 0$, $p_n^2\sqrt{\log(p_n)}A_n = o(\gamma_n\sqrt{n})$ and $\sqrt{n/(p^3_ns_n\log(p_n))} \gamma_n \rightarrow \infty$. Then with probability tending to one, $\widehat{\theta}_{n\Lambda,k}=0$ for all $k \in \Sc^c_n$.\\
In addition to Assumption  \ref{assumption_space_lambda}, assume that the factor loading matrix belongs to the space $\{\Lambda_n\in\Rb^{p_n \times m} \mid \underline{\mu}\leq \lambda_{\min}(p^{-\nu}_n\Lambda^\top_n\Lambda_n)\leq \lambda_{\max}(p^{-\nu}_n\Lambda^\top_n\Lambda_n) \leq \overline{\mu}\}$ where $\underline{\mu},\overline{\mu}>0$ and $\nu \in ]1/2,1]$. Then the GLS estimator $\widehat{F}_t = (\widehat{\Lambda}^\top_n\widehat{\Psi}^{-1}_n \widehat{\Lambda}_n)^{-1} 
\widehat{\Lambda}^\top_n \widehat{\Psi}^{-1}_n X_t$ of $F_t$ satisfies, for $t=1,\ldots,n$, $\|\widehat{F}_t-F_t\| = o_p(1)$.
\end{theorem}
The consistent estimation of $\widehat{\theta}_n$ stated in Theorem \ref{sparsistency} relies on \cite{anderson1956}'s condition for identification. This restriction allows us to get $\|\widehat{\theta}_n-\theta^\ast_n\|_\infty \leq C_2 u_n$ from $\|\widehat{\Sigma}_n-\Sigma^\ast_n\|_2\leq C_1 u_n$ established in Theorem \ref{Theorem_existence_consistent}, with $C_1,C_2>0$ and $u_n = p_n\sqrt{s_n\log(p_n)/n}+\sqrt{p_n}A_n$: this is due to the ``strong identifiability'' property of the factor model stated in Lemma 1 and Theorem 1 in \cite{kano1983}, a property ensured by Assumptions \ref{assumption_AR_condition} and \ref{assumption_space_lambda}. Further details are provided in Section \ref{appendix_proofs} of the Appendix. The sparsistency property follows from the consistency of $\widehat{\theta}_n$ and suitable conditions on $p(\cdot,\gamma_n)$ and $\gamma_n$. 
The consistency of the GLS estimator $\widehat{F}_t$ is proved under weak factors.  


\section{Simulations}\label{simulations}

This section presents two data generating processes (DGPs). First, we generate the $p_n$-dimensional random vector $X_t$ based on the DGP:
\begin{equation}\label{dgp1}
X_t \sim \Nc_{\Rb^p}\big(0,\Sigma^\ast_n), \; \text{with} \; \Sigma^\ast_n = \Lambda^\ast_n \Lambda^{\ast\top}_n + \Psi^\ast_n, \; \text{for} \; t = 1,\ldots,n, 
\end{equation}
with $\Lambda^\ast_n$ assumed $s_n$-sparse, $s_n:=|\Sc_n|$. The sparsity patterns for $\Lambda^\ast_n$ are specified as follows:
\begin{itemize}
    \item[(i)] \textbf{Perfect simple structure with non-sparse blocks; $\Psi^\ast_n$ diagonal}\\
    $\Lambda^\ast_n$ satisfies the perfect simple structure for $p_n = 60, 120,180$ and $m = 3,4$, 
    and $s_n = p_n$ for any $p_n$. 
    When $p_n=60, m=3$ (resp. $m=4$), then $p_n \, m = 180$ (resp. $p_n\,m=240$).
    When $p_n=120, m=3$ (resp. $m=4$), then $p_n \, m=360$ (resp. $p_n\,m=480$).
    When $p_n=180, m=3$, then $p_n \,m=540$.
    \item[(ii)] \textbf{Perfect simple structure with overlaps, non-sparse blocks; $\Psi^\ast_n$ diagonal}\\
    $\Lambda^\ast_n$ satisfies a near perfect simple structure with overlaps  
    for $p_n = 60, 120,180$ and $m = 3,4$.
    We set the number of overlaps as $\lceil 0.5\times p_n/m \rceil$ in each case. 
    When $p_n=60, m=3$ (resp. $m=4$), then $p_n \, m= 180$ (resp. $p_n \, m=240$) and $s_n=80$ (resp. $s_n=84$),
    where the size of the overlaps is $10$ (resp. $8$). 
    When $p_n=120, m=3$ (resp. $m=4$), then $p_n \, m=360$ (resp. $p_n \, m=480$) and $s_n=160$ (resp. $s_n = 165$), where the size of the overlaps is $20$ (resp. $15$). 
    When $p_n=180, m=3$, then $p_n \, m=540$ and $s_n = 240$, where the size of the overlaps is $30$.
    \item[(iii)] \textbf{Perfect simple structure with overlaps, sparse blocks; $\Psi^\ast_n$ diagonal}\\
    $\Lambda^\ast_n$ satisfies the same pattern as in (ii), but we randomly replace some non-zero coefficients by zeros such that the ratio of zero coefficients in $\Lambda^\ast_n$ is $70\%$.
    When $p_n=60, m=3$ (resp. $m=4$), then $p_n \, m= 180$ (resp. $p_n \, m=240$) and $s_n=54$ (resp. $s_n=72$). 
    When $p_n=120, m=3$ (resp. $m=4$), then $p_n \, m=360$ (resp. $p_n \,m =480$) and $s_n=108$ (resp. $s_n = 144$). 
    When $p_n=180, m=3$, then $p_n \, m=540$ and $s_n = 162$.
    \item[(iv)] \textbf{General arbitrary sparse structure; $\Psi^\ast_n$ diagonal}\\
    The zero entries of $\Lambda^\ast_n$ are randomly set. When $p_n=60, m=3$ (resp. $m=4$), then $p\, m= 180$ (resp. $p\,m=240$) and $s_n=27$ (resp. $s_n=36$).\\
    When $p_n=120, m=3$ (resp. $m=4$), then $p_n \, m=360$ (resp. $p_n \, m =480$) and $s_n=54$ (resp. $s_n=72$). 
    When $p_n=180, m=3$, then $p_n \, m=540$ and $s_n=81$. 
    \item[(v)] \textbf{General arbitrary sparse structure; $\Psi^\ast_n$ sparse non-diagonal}\\
    In this final scenario, the sparsity pattern in $\Lambda^\ast_n$ is the same as in setting (iv). However, the idiosyncratic error variables are potentially correlated: the true variance-covariance $\Psi^\ast_n$ is a sparse non-diagonal matrix.  
\end{itemize}
Note that $s_n$ significantly depends on $m$ in (i)-(ii): in (ii), $44\%$ (resp. $35\%$) of the entries of $\Lambda^\ast_n$ are non-zero when $p_n=60,m=3$ (resp. $p_n=60,m=4$); the non-zero entries represent $15\%$ of all the coefficients of $\Lambda^\ast_n$ for any $p_n,m$,  in (iv)-(v).
For each sparsity pattern and $(p_n,m,n)$, we draw two hundred batches of $n$ independent samples from $X_t \sim \Nc_{\Rb^p}(0,\Sigma^\ast_n)$. For each batch, the diagonal elements of $\Psi^\ast_n$ are simulated in the uniform distribution $\Uc([0.5,1])$ for patterns (i)--(iv); in (v), the proportion of zero coefficients in the lower diagonal part of $\Psi^\ast_n$ is set as $75\%$ and the non-zero entries are generated as the sum of one or more normally distributed variables, so that $\lambda_{\min}(\Psi^\ast_n)>0$; the non-zero parameters of $\Lambda^\ast_n$ are simulated in $\Uc([-2,-0.5]\cup [0.5,2])$. For each dimension setting $(p_n,m)$, the sample size is $n=250, 500$. 

\mds

The second DGP generates the observations $X_t$ according to the time series dynamic:
\begin{equation}\label{dgp2}
X_t = \Lambda^\ast_n F_t + \eps_t, \; \text{for} \; t = 1,\ldots,n, 
\end{equation}
where $\eps_t \sim \Nc_{\Rb^{p_n}}(0,\Psi^\ast_n)$ and $F_{t,k} = \phi_k F_{t-1,k} + \zeta_{t,k}$ for $t = 1,\ldots,n$, $k=1,\ldots,m$, with $\phi_k \sim \Uc([0.5,0.85])$, $\zeta_{t,k} \sim \Nc_{\Rb}(0,(1-\phi^2_k))$. The parameters $\Lambda^\ast_n$ and $\Psi^\ast_n$ are generated as in setting (v), DGP (\ref{dgp1}). We also draw two hundred batches of $n$ observations. 

\mds

The sparsity pattern in $\Lambda^\ast_n$ remains identical for all batches and a given DGP, but the non-zero coefficients may have different locations in $\Lambda^\ast_n$ for the sparsity patterns (iii)--(v); under the sparsity patterns (i) and (ii), the location of the non-zero coefficients remains unchanged for all the batches; for all patterns, the cardinality $s_n$ remains unchanged. In setting (v), the proportion of zero entries remains unchanged for all the batches but the non-zero coefficients of $\Psi^\ast_n$ may have different locations. 

\mds

For each batch, we apply our sparsity-based method (\ref{stat_crit}) for the Gaussian and least squares losses, which will be denoted by SGF (Sparse Gaussian Factor) and SLSF (Sparse Least Squares Factor), respectively, hereafter. We specify $a_{\text{scad}}=3.7$, a value identified as optimal in \cite{fan2001} by cross-validated experiments. The MCP parameter is set as $b_{\text{mcp}}=3.5$, following Section 5 in \cite{loh2015}. The penalized problem is solved by a gradient descent algorithm. For DGP (\ref{dgp1}), we employ a $5$-fold cross-validation procedure to tune $\gamma_n$; when the data arises from (\ref{dgp2}) and due to its time-series nature, we employ an ``out-of-sample''-based cross-validation procedure: further details on the implementation and the cross-validation procedure can be found in Section \ref{implementation} of the Appendix. To emphasize the efficiency of our approach in terms of support recovery and accuracy, we employ the SOFAR estimator, in which case we use the R package ``rrpack'' of \cite{uematsu2023a}.

\mds

We report the variable selection performance through the percentage of zero coefficients in $\Lambda_n^\ast$ correctly estimated, denoted by $\text{C1}$, and the percentage of non-zero coefficients in $\Lambda_n^\ast$ correctly identified as such, denoted by $\text{C2}$. The mean squared error (MSE), defined as $\|\widehat\theta_{n\Lambda}-\theta^\ast_{n\Lambda}\|^2_2$, is reported as an estimation accuracy measure. Due to the rotational indeterminacy, we compel the SOFAR estimator to satisfy the specific identification condition to remove this indeterminacy. Thus, the target loading matrix of SOFAR could be different from the sparsest true loading $\Lambda^\ast_n$. More precisely, we compare the loading matrix estimated by SOFAR, denoted by $\widehat{\Lambda}_{n,\text{SOFAR}}$, with both the true $\Lambda^\ast_n$ and the rotated true $\widehat{\Lambda}_{n,\text{SOFAR}}^\ast$ defined as: $\widehat{\Lambda}_{n,\text{SOFAR}}^\ast := \Lambda^\ast_n \widehat{R}_{\text{SOFAR}}^\ast,\,
\widehat{R}_{\text{SOFAR}}^\ast 
:=
\mathop{\arg\min}_{R \in \mathcal{O}(m)}\|\widehat{\Lambda}_{n,\text{SOFAR}} - \Lambda^\ast_n R\|^2_F$.
In each table, in the SOFAR's column, the MSE values compared with $\widehat{\Lambda}_\text{n,SOFAR}^\ast$ are reported in the parenthesis. To be precise, for each target of SOFAR, we selected the best values of the regularization coefficient among various values in the sense of MSE. 

\mds

When the data is generated from DGP (\ref{dgp1}), setting~(i), the true factor loading matrix satisfies the orthogonal constraint, i.e., the restriction stated in \cite{uematsu2023a}. Thus, in light of Table~\ref{support_perfect_simple_no_over_no_sparse}, our proposed method and the SOFAR provide similar performances. 
In contrast to the perfect simple structure, $\Lambda^\ast_n$ generated from settings (ii)--(v) does not satisfy the orthogonal constraint, meaning that $\Lambda^\ast_n\Lambda^{\ast\top}_n$ is non-diagonal. 
Tables~\ref{support_perfect_simple_over_no_sparse}--\ref{support_arbitrary} suggest that the performances of our method are better than the SOFAR in terms of support recovery. From the viewpoint of MSE, it can be noted that our method outperforms the SOFAR, although the MSE results based on $\widehat{\Lambda}_{n,\text{SOFAR}}^\ast$ are smaller. However, if we consider the MSE with $\widehat\Lambda_{n,\text{SOFAR}}$, the performance of the SOFAR dramatically worsens as the sparsest true loading differs from the orthogonal structure. As detailed in Section~\ref{sec:1}, this is because the target of the SOFAR is not necessarily the same as the sparsest true loading matrix. Interestingly, when the data is generated based on setting (v), the recovery/MSE performances of our method displayed in Table \ref{support_arbitrary_misspec} are comparable to those reported in Table \ref{support_arbitrary}: although $\Psi^\ast_n$ is estimated under the diagonal constraint, the procedure does not suffer from misspecification. In the time-series case (\ref{dgp2}), our approach outperforms the SOFAR method according to Table \ref{support_arbitrary_misspec_ts}. Although the MSEs are larger compared to the i.i.d. case, the hierarchy in terms of model performances remains unchanged.

\mds

Another point worth mentioning is that the Gaussian loss-based estimation provides better performances in terms of MSE compared to the least squares loss, for any sparsity pattern in $\Lambda^\ast_n$. This is in line with the numerical findings of \cite{poignard2023_copula} related to the sparse estimation of the covariance matrix of the Gaussian copula by SCAD and MCP when the loss is the Gaussian likelihood and the least squares function.

\mds

\begin{figure}[t]
  \begin{minipage}[b]{0.31\hsize}
    \centering
    \includegraphics[width=0.9\textwidth]{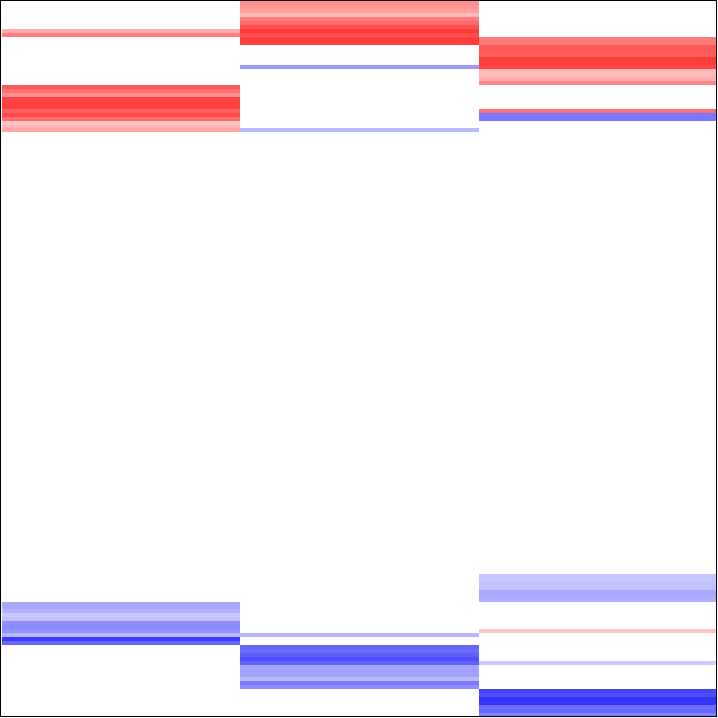}
    \subcaption{True loading matrix}\label{fig:2-true}
  \end{minipage}
  \begin{minipage}[b]{0.31\hsize}
    \centering
    \includegraphics[width=0.9\textwidth]{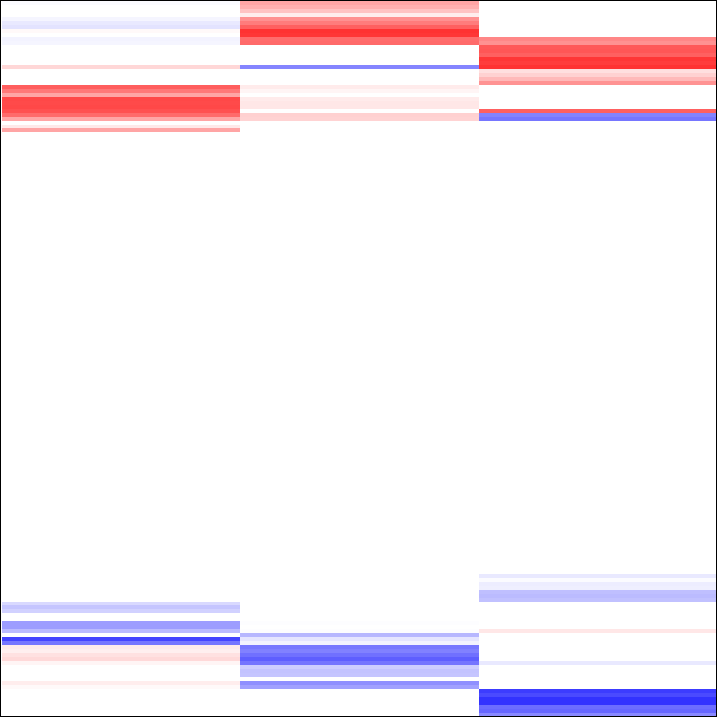}
    \subcaption{SOFAR estimator}\label{fig:2-sofar}
  \end{minipage}
  \begin{minipage}[b]{0.31\hsize}
    \centering
    \includegraphics[width=0.9\textwidth]{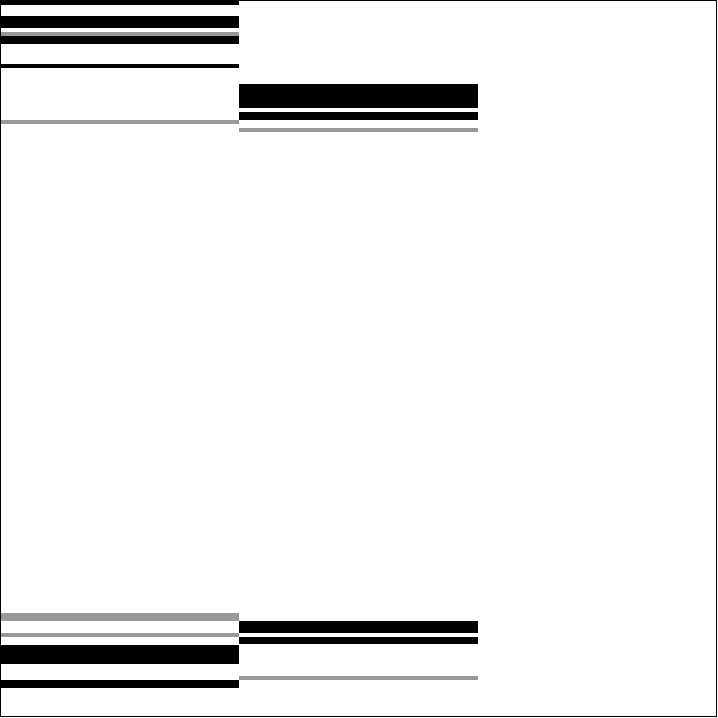}
    \subcaption{Misspecified elements}\label{fig:2-ic}
  \end{minipage}
  \begin{minipage}[b]{0.31\hsize}
    \centering
    \includegraphics[width=0.9\textwidth]{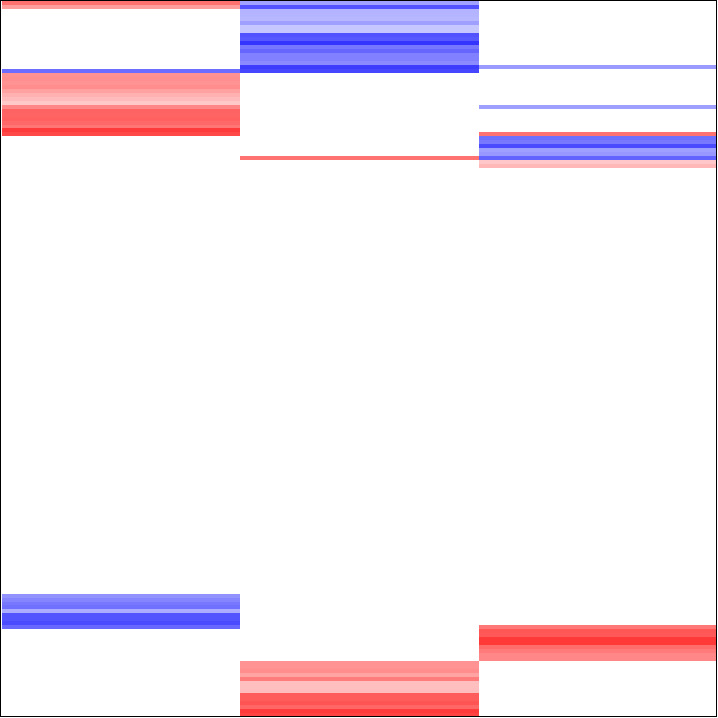}
    \subcaption{True loading matrix}\label{fig:7-true}
  \end{minipage}
  \hspace*{0.46cm}\begin{minipage}[b]{0.31\hsize}
    \centering
    \includegraphics[width=0.9\textwidth]{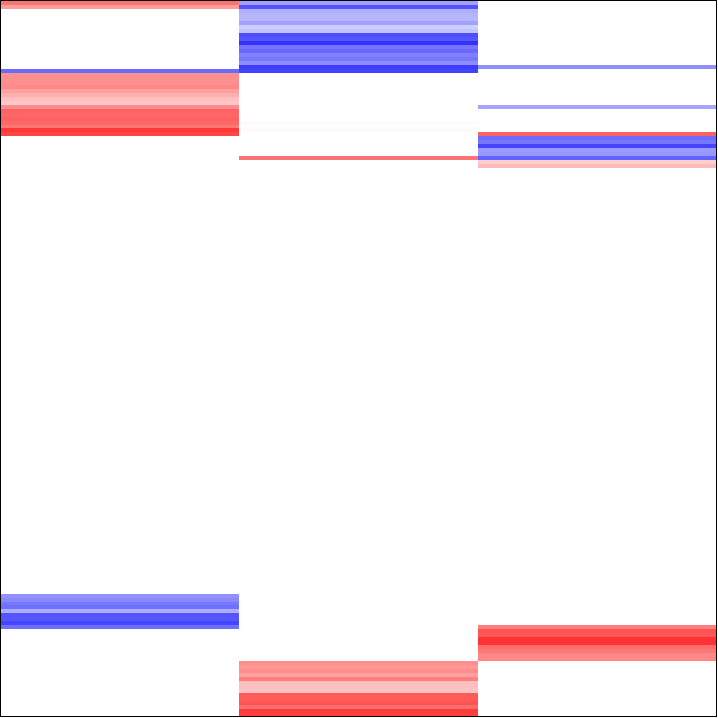}
    \subcaption{SOFAR estimator}\label{fig:7-sofar}
  \end{minipage}
  \hspace*{0.46cm}\begin{minipage}[b]{0.31\hsize}
    \centering
    \includegraphics[width=0.9\textwidth]{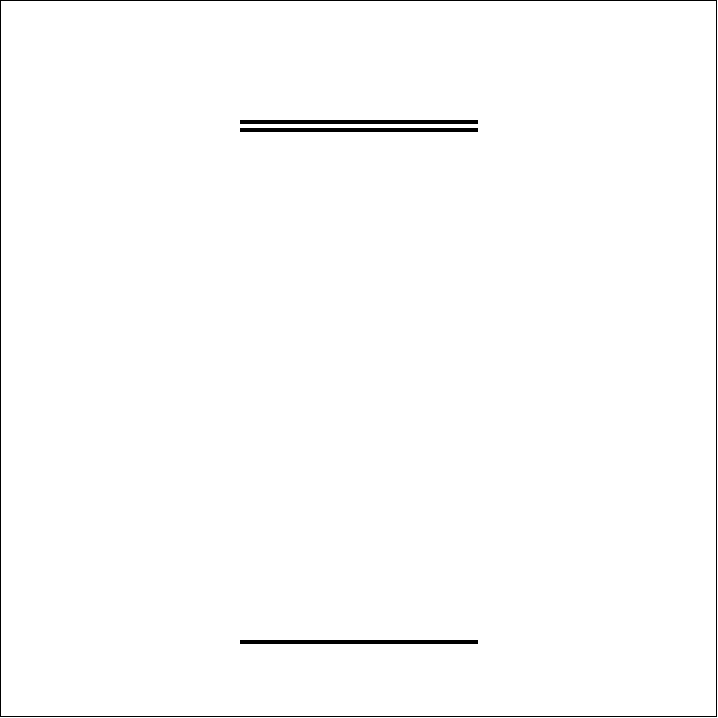}
    \subcaption{Misspecified elements}\label{fig:7-ic}
  \end{minipage}
  \caption{Two examples corresponding to Setting~(iv). (a): rearranged true loading matrix; (b): SOFAR estimator for (a); (c) miss-specified elements (black: incorrectly identified non-zero, grey: incorrectly identified zero). (d), (e), (f): alternative example.}\label{fig:setting-4}
\end{figure}

The results displayed in Setting~(iv) suggest that the SOFAR performs well in terms of recovery, where the C1 and C2 metrics as similar to our approach. 
The main reason is that the generated true $\Lambda^\ast_n$ are approximately column orthogonal for this Setting~(iv). That is, the off-diagonal elements of $\Lambda^{\ast\top}_n\Lambda^\ast_n$ are much smaller than its diagonal elements. To illustrate this point, Figure~\ref{fig:setting-4} displays two examples relating to Setting~(iv), with $n=1000$ and $p_n=180$. 
As displayed in (a, b, c), 
there are several misspecified elements in the SOFAR estimator when the true loading matrix slightly diverges from the case of being a column orthogonal matrix. 
On the other hand, as shown in (d), the SOFAR estimator recovers the true zero/non-zero entries almost perfectly when the true $\Lambda^\ast_n$ is nearly column orthogonal.

\begin{landscape}
\begin{table}[H]\centering\caption{Model selection and precision accuracy based on 200 replications, DGP (\ref{dgp1}), perfect simple structure without overlaps and non-sparse blocks (sparsity pattern (i)) with respect to $(n,p_n,m)$.\label{support_perfect_simple_no_over_no_sparse}}
\scalebox{0.9}{\begin{tabular}{cc|ccc|ccc}\hline\hline
 &    & \multicolumn{3}{c|}{$n=250$}       & \multicolumn{3}{c}{$n=500$}\\ 
\multicolumn{2}{c|}{$(p_n,m)$ }  & SGF & SLSF & SOFAR & SGF & SLSF & SOFAR \\ \hline            
                
 & C1 & $79.82 \; - \; 84.86$ & $ 75.42 \; - \; 82.02$ & $91.82$ 
      & $89.32 \; - \; 92.65$  & $81.96 \; - \;  87.16 $ & $94.14$  \\
\scalebox{0.95}{$(60,3)$}
& C2 & $ 100 \; - \; 100 $ & $ 100 \; - \; 100 $ & $98.39$ 
      & $100 \; - \; 100$ & $100 \; - \; 100$ & $99.78$ \\
 & MSE & $ 0.4588 \; - \; 0.4545 $ & $0.5766 \; - \; 0.5797$ &  $1.6931\;(0.6800)$
      & $0.1989 \; - \; 0.1978$ & $0.2459 \; - \; 0.2510$ & $0.5738\;(0.3060)$  \\
\hline
& C1 &  $82.48 \; - \; 87.94$ & $82.47 \; - \; 88.25$ & $94.49$ 
      & $91.30 \; - \; 94.18$ & $87.30 \; - \; 91.66$ & $96.65$  \\
\scalebox{0.95}{$(60,4)$}
& C2  & $100 \; - \;100 $ & $100 \; - \;100 $ &  $98.03$
      & $100 \; - \; 100$ & $100 \; - \; 100$ &  $99.63$\\
& MSE & $0.4751 \; - \; 0.4729$ & $0.6013 \; - \; 0.7343$ & $2.0320\;(0.8331)$ 
      & $0.2048 \; - \; 0.2036$ & $0.2512 \; - \; 0.2525$ & $0.7455\;(0.3679)$  \\
\hline   
                            
& C1 & $ 78.63 \; - \; 84.27$ & $78.16 \; - \; 84.54$ & $91.48$ 
      & $89.50 \; - \; 92.80$ & $82.06 \; - \; 87.08$ & $94.21$ \\
\scalebox{0.95}{$(120,3)$}
& C2 & $100 \; - \; 100$ & $100 \; - \; 100$ & $98.01$
      & $100 \; - \; 100$ & $100 \; - \; 100$ & $99.61$ \\
& MSE & $0.8781 \; - \; 0.8683$ & $1.1310 \; - \; 1.1321$ & $3.9879\;(1.4755)$ 
      & $0.4081 \; - \; 0.4050$ & $0.5108 \; - \; 0.5131$ & $1.4386\;(0.6683)$  \\
\hline        
& C1 & $80.12 \; - \; 85.67$ & $83.17 \; - \; 89.17$ & $93.94$ 
      & $91.13 \; - \; 94.18$ & $88.43 \; - \; 92.41$ & $95.85$  \\
\scalebox{0.95}{$(120,4)$}
& C2 & $100 \; - \; 100$ & $100 \; - \; 100$ &  $97.41$
      & $100 \; - \; 100$ & $100 \; - \; 100$ & $99.43$ \\
& MSE & $0.9378 \; - \; 0.9321$ & $1.1462 \; - \; 1.1482$ & $5.0538\;(1.8141)$ 
      & $0.4432 \; - \; 0.4395$ & $0.5273 \; - \; 0.5283$ & $1.9252\;(0.8381)$  \\
\hline
& C1 & $77.58 \; - \; 83.29$ & $75.44 \; - \; 82.03$ & $91.35$
      & $89.48 \; - \; 92.87$ & $81.46 \; - \; 86.39$ & $93.29$  \\
\scalebox{0.95}{$(180,3)$}
& C2 & $100 \; - \; 100$ & $100 \; - \; 100$ & $97.91$
      & $100 \; - \; 100$ & $100 \; - \; 100$ & $99.56$ \\
& MSE & $1.3781 \; - \; 1.3684$ & $1.7576 \; - \; 1.7638$ & $6.2790\;(2.3915)$ 
      & $0.6183 \; - \; 0.6144$ & $0.7745 \; - \; 0.7809$ & $2.4009\;(1.0297)$  \\
\hline\hline
\end{tabular}}
\begin{minipage}{21.5cm}
\vspace{3pt}
\footnotesize Note: the column ``SGF'' refers to the estimator deduced from the Gaussian loss with SCAD and MCP penalization, respectively. The column ``SLSF'' refers to the estimator deduced from the least squares loss with SCAD and MCP penalization, respectively.\\
The MSE values for SOFAR in parenthesis are based on $\widehat{\Lambda}_\text{n,SOFAR}^\ast$, with $\widehat{\Lambda}_\text{n,SOFAR}^\ast := \Lambda^\ast_n \widehat{R}_\text{SOFAR}^\ast,\;
\widehat{R}_\text{SOFAR}^\ast 
:=
\underset{R \in \mathcal{O}(m)}{\arg\min}\|\widehat{\Lambda}_{n,\text{SOFAR}} - \Lambda^\ast_n R\|^2_F$.
\end{minipage}
\end{table}
\end{landscape}

\begin{landscape}
\begin{table}[H]\centering\caption{Model selection and precision accuracy based on 200 replications, DGP (\ref{dgp1}), perfect simple structure with overlaps and non-sparse blocks (sparsity pattern (ii)) with respect to $(n,p_n,m)$.\label{support_perfect_simple_over_no_sparse}}
\scalebox{0.9}{\begin{tabular}{cc|ccc|ccc}\hline\hline
 &    & \multicolumn{3}{c|}{$n=250$}       & \multicolumn{3}{c}{$n=500$}\\ 
\multicolumn{2}{c|}{$(p_n,m)$ }  & SGF & SLSF & SOFAR & SGF & SLSF & SOFAR \\ \hline
 & C1 & $77.69 \; - \; 83.45$ & $70.36 \; - \; 77.34$ & $78.40$ 
      & $87.77 \; - \; 91.16$ & $76.97 \; - \; 82.49$ & $77.99$  \\
\scalebox{0.95}{$(60,3)$}
& C2 & $100 \; - \; 100$ & $100 \; - \; 100$ & $90.29$ 
      & $100 \; - \; 100$ & $100 \; - \; 100$ & $90.16$ \\
 & MSE & $0.5901 \; - \; 0.5915$ & $0.8379 \; - \; 0.8507$ & $16.4037\;(1.0354)$
      & $0.2794 \; - \; 0.2781$ & $0.4027 \; - \; 0.4047$ & $15.8031\;(0.5070)$  \\
\hline
& C1 &  $80.07 \; - \; 85.86$ & $75.67 \; - \; 82.70$ & $82.43$
      & $89.89 \; - \; 93.23$ & $81.90 \; - \; 87.04$ & $81.10$  \\
\scalebox{0.95}{$(60,4)$}
& C2 &  $100 \; - \; 100$ & $100 \; - \; 100$ &  $86.42$
      & $100 \; - \; 100$ & $100 \; - \; 100$ &  $87.80$\\
& MSE & $0.6667 \; - \; 0.6638$ & $1.0263 \; - \; 1.0450$ & $23.8262\;(1.3628)$ 
      & $0.3003 \; - \; 0.2989$ & $0.4574 \; - \; 0.4660$ & $24.9719\;(0.7080)$ \\
\hline            
& C1 & $76.29 \; - \; 82.36$ & $68.92 \; - \; 76.22$ & $80.27$ 
      & $87.64 \; - \; 91.57$ & $78.05 \; - \; 83.72$ & $80.27$  \\
\scalebox{0.95}{$(120,3)$}
& C2 & $100 \; - \; 100$ & $100 \; - \; 100$ & $90.65$ 
      & $100 \; - \; 100$ & $100 \; - \; 100$ & $92.77$ \\
& MSE & $1.0821 \; - \; 1.0805$ & $1.6303 \; - \; 1.6473$ & $30.6147\;(2.0938)$ 
      & $0.5558 \; - \; 0.5543$ & $0.7927 \; - \; 0.7980$ & $23.3241\;(1.0257)$  \\
\hline    
& C1 & $78.69 \; - \; 84.36$ & $75.80 \; - \; 82.48$ & $85.12$ 
      & $89.42 \; - \; 92.76$ & $83.62 \; - \; 88.69$ & $83.44$  \\
\scalebox{0.95}{$(120,4)$}
& C2 & $100 \; - \; 100$ & $100 \; - \; 100$ &  $87.18$
      & $100 \; - \; 100$ & $100 \; - \; 100$ & $90.55$ \\
& MSE & $1.2606 \; - \; 1.2517$ & $1.9371 \; - \; 1.9716$ & $39.1837\;(2.8719)$
      & $0.5805 \; - \; 0.5776$ & $0.8480 \; - \; 0.8564$ & $34.6621\;(1.3524)$  \\
\hline
& C1 & $76.17 \; - \; 81.87$ & $69.27 \; - \; 76.42$ & $80.91$ 
      & $87.44 \; - \; 91.19$ & $76.69 \; - \; 82.94$ & $82.03$  \\
\scalebox{0.95}{$(180,3)$}
& C2 & $100 \; - \; 100$ & $100 \; - \; 100$ & $92.01$ 
      & $100 \; - \; 100$ & $100 \; - \; 100$ & $91.81$ \\
& MSE & $1.7192 \; - \; 1.7126$ & $2.5399 \; - \; 2.5583$ & $41.5858\;(2.9853)$ 
      & $0.8787 \; - \; 0.8771$ & $1.2331 \; - \; 1.2419$ & $37.5806\;(1.4952)$  \\
\hline\hline
\end{tabular}}
\begin{minipage}{21.5cm}
\vspace{3pt}
\footnotesize Note: the column ``SGF'' refers to the estimator deduced from the Gaussian loss with SCAD and MCP penalization, respectively. The column ``SLSF'' refers to the estimator deduced from the least squares loss with SCAD and MCP penalization, respectively.\\
The MSE values for SOFAR in parenthesis are based on $\widehat{\Lambda}_\text{n,SOFAR}^\ast$, with $\widehat{\Lambda}_\text{n,SOFAR}^\ast := \Lambda^\ast_n \widehat{R}_\text{SOFAR}^\ast,\;
\widehat{R}_\text{SOFAR}^\ast 
:=
\underset{R \in \mathcal{O}(m)}{\arg\min}\|\widehat{\Lambda}_{n,\text{SOFAR}} - \Lambda^\ast_n R\|^2_F$.
\end{minipage}
\end{table}
\end{landscape}

\begin{landscape}
\begin{table}[H]\centering\caption{Model selection and precision accuracy based on 200 replications, DGP (\ref{dgp1}), perfect simple structure with overlaps and sparse blocks (sparsity pattern (iii)) with respect to $(n,p_n,m)$.\label{support_perfect_simple_over_sparse}}
\scalebox{0.9}{\begin{tabular}{cc|ccc|ccc}\hline\hline
 &    & \multicolumn{3}{c|}{$n=250$}       & \multicolumn{3}{c}{$n=500$}\\
\multicolumn{2}{c|}{$(p_n,m)$ }  & SGF & SLSF & SOFAR & SGF & SLSF & SOFAR \\ \hline                
 & C1 & $80.23 \; - \; 85.65$ & $82.79 \; - \; 88.06$ & $87.98$ 
      & $89.96 \; - \; 93.13$ & $88.00 \; - \; 91.79$ & $86.83$  \\
\scalebox{0.95}{$(60,3)$}
& C2 & $100 \; - \; 100$ & $100 \; - \; 100$ & $92.33$ 
      & $100 \; - \; 100$ & $100 \; - \; 100$ & $94.20$ \\
 & MSE & $0.4210 \; - \; 0.4177$ & $0.5458 \; - \; 0.5557$ & $7.8310\;(0.7397)$ 
      & $0.2069 \; - \; 0.2056$ & $0.2734 \; - \; 0.2756$ & $6.4450\;(0.3939)$  \\
\hline
& C1 & $80.83 \; - \; 86.61$ & $81.79 \; - \; 87.47$ & $85.93$ 
      & $90.53 \; - \; 93.68$ & $87.01 \; - \; 91.14$ & $85.81$  \\
\scalebox{0.95}{$(60,4)$}
& C2 & $100 \; - \; 100$ & $100 \; - \; 100$ &  $88.01$
      & $100 \; - \; 100$ & $100 \; - \; 100$ & $90.01$ \\
& MSE & $0.5867 \; - \; 0.5828 $ & $0.8388 \; - \; 0.8504$ & $18.39\;(1.2519)$ 
      & $0.2749 \; - \; 0.2731$ & $ 0.3852 \; - \; 0.3916$ & $15.7564\;(0.6239)$  \\
\hline  
& C1 & $79.48 \; - \; 85.03$ & $83.17 \; - \; 88.46$ & $89.00$ 
      & $90.01 \; - \; 93.27$ & $88.41 \; - \; 92.23$ & $89.16$  \\
\scalebox{0.95}{$(120,3)$}
& C2 & $100 \; - \; 100$ & $100 \; - \; 100$ & $93.90$ 
      & $100 \; - \; 100$ & $100 \; - \; 100$ & $94.75$ \\
& MSE & $0.8292 \; - \; 0.8196$ & $1.0983 \; - \; 1.1078$ & $11.8534\;(1.4481)$
      & $0.3654 \; - \; 0.3620$ & $0.4825 \; - \; 0.4852$ & $9.9180\;(0.7379)$  \\
\hline           
& C1 & $79.28 \; - \; 84.77$ & $80.85 \; - \; 86.85$ & $87.95$ 
      & $90.27 \; - \; 93.51$ & $86.52 \; - \; 90.59$ & $87.10$  \\
\scalebox{0.95}{$(120,4)$}
& C2 & $100 \; - \; 100$ & $99.99 \; - \; 99.98$ &  $89.67$
      & $100 \; - \; 100$ & $100 \; - \; 100$ & $91.07$ \\
& MSE & $1.1254 \; - \; 1.1185$ & $1.6724 \; - \; 2.0718$ & $28.2858\;(2.4809)$ 
      & $0.5237 \; - \; 0.5198$ & $0.7475 \; - \; 0.7565$ & $27.0152\;(1.2358)$  \\
\hline        
& C1 & $78.59 \; - \; 84.07$ & $82.81 \; - \; 88.17$ &  $90.16$
      & $89.80 \; - \; 93.19$ & $89.25 \; - \; 92.77$ & $89.78$  \\
\scalebox{0.95}{$(180,3)$}
& C2 & $100 \; - \; 100$ & $100 \; - \; 100$ & $94.15$ 
      & $100 \; - \; 100$ & $100 \; - \; 100$ & $94.88$ \\
& MSE & $1.2245 \; - \; 1.2107$ & $1.6198 \; - \; 1.6406$ & $16.3631\;(2.1675)$ 
      & $0.5678 \; - \; 0.5626$ & $0.7420 \; - \; 0.7530$ & $12.7403\;(1.0939)$ \\
\hline\hline
\end{tabular}}
\begin{minipage}{21.5cm}
\vspace{3pt}
\footnotesize Note: the column ``SGF'' refers to the estimator deduced from the Gaussian loss with SCAD and MCP penalization, respectively. The column ``SLSF'' refers to the estimator deduced from the least squares loss with SCAD and MCP penalization, respectively.\\
The MSE values for SOFAR in parenthesis are based on $\widehat{\Lambda}_\text{n,SOFAR}^\ast$, with $\widehat{\Lambda}_\text{n,SOFAR}^\ast := \Lambda^\ast_n \widehat{R}_\text{SOFAR}^\ast,\;
\widehat{R}_\text{SOFAR}^\ast 
:=
\underset{R \in \mathcal{O}(m)}{\arg\min}\|\widehat{\Lambda}_{n,\text{SOFAR}} - \Lambda^\ast_n R\|^2_F$.
\end{minipage}
\end{table}
\end{landscape}

\begin{landscape}
\begin{table}[H]\centering\caption{Model selection and precision accuracy based on 200 replications, DGP (\ref{dgp1}), arbitrary sparse structure and $\Psi^\ast_n$ diagonal (sparsity pattern (iv)) with respect to $(n,p_n,m)$. \label{support_arbitrary}}
\scalebox{0.9}{\begin{tabular}{cc|ccc|ccc}\hline\hline
 &    & \multicolumn{3}{c|}{$n=250$}       & \multicolumn{3}{c}{$n=500$}\\ 
\multicolumn{2}{c|}{$(p_n,m)$ }  & SGF & SLSF & SOFAR & SGF & SLSF & SOFAR \\ \hline
& C1 & $85.64 \; - \; 90.41$ & $94.71 \; - \; 96.49$ &  $93.75$
      & $92.55 \; - \; 95.02$ & $96.84 \; - \; 97.77$ & $93.39$  \\
\scalebox{0.95}{$(60,3)$}
& C2 & $100 \; - \; 100$ & $100 \; - \; 99.96$ &  $92.98$
      & $99.98 \; - \; 99.98$ & $99.98 \; - \; 99.98$ & $95.43$ \\
 & MSE & $0.2291 \; - \; 0.2248$ & $0.2498 \; - \; 0.2576$ & $4.0267\;(0.4398)$
      & $0.1176 \; - \; 0.1162$ & $0.1325 \; - \; 0.1347$ & $2.7778\;(0.2382)$  \\
\hline    
& C1 & $86.82 \; - \; 91.19$ & $94.45 \; - \; 96.43$ & $93.02$ 
      & $93.04 \; - \; 95.61$ & $96.57 \; - \; 97.55$ & $92.53$  \\
\scalebox{0.95}{$(60,4)$}
& C2 & $100 \; - \; 100$ & $100 \; - \; 99.99$ & $90.15$ 
      & $100 \; - \; 100$ & $100 \; - \; 100$ & $90.64$ \\
& MSE & $0.3290 \; - \; 0.3232$ & $0.3989 \; - \; 0.4080$ & $7.2800\;(0.7313)$ 
      & $0.1505 \; - \; 0.1489$ & $0.1837 \; - \; 0.1899$ & $7.7528\;(0.3875)$ \\
\hline 
& C1 & $83.91 \; - \; 89.20$ & $95.03 \; - \; 96.91$ &  $94.19$
      & $92.95 \; - \; 95.49$ & $96.83 \; - \; 97.89$ &  $93.88$ \\
\scalebox{0.95}{$(120,3)$}
& C2 & $100 \; - \; 100$ & $100 \; - \; 100$ &  $94.09$
      & $100 \; - \; 100$ & $100 \; - \; 100$ & $95.44$ \\
& MSE & $0.4946 \; - \; 0.4821$ & $0.5103 \; - \; 0.5207$ & $5.6643\;(0.8235)$ 
      & $0.1991 \; - \; 0.1960$ & $0.2318 \; - \; 0.2366$ & $5.3639\;(0.4273)$  \\
\hline
& C1 & $83.66 \; - \; 88.90$ & $94.65 \; - \; 96.81$ & $93.59$ 
      & $92.89 \; - \; 95.51$ & $96.88 \; - \; 97.89$ & $92.74$  \\
\scalebox{0.95}{$(120,4)$}
& C2 & $100 \; - \; 100$ & $100 \; - \; 100$ & $91.56$ 
      & $100 \; - \; 100$ & $100 \; - \; 100$ & $93.20$ \\
& MSE & $0.6346 \; - \; 0.6179$ & $0.7427 \; - \; 0.7609$ & $11.9326\;(1.4141)$ 
      & $0.2787 \; - \; 0.2746$ & $0.3481 \; - \; 0.3548$ & $10.6349\;(0.6905)$  \\
\hline       
& C1 & $82.21 \; - \; 87.69$ & $95.17 \; - \; 97.05$ & $94.60$ 
      & $92.42 \; - \; 95.16$ & $97.29 \; - \; 98.18$ & $94.27$  \\
\scalebox{0.95}{$(180,3)$}
& C2 & $100 \; - \; 100$ & $100 \; - \; 100$ &  $94.37$
      & $100 \; - \; 100$ & $100 \; - \; 100$ & $95.52$\\
& MSE & $0.7393 \; - \; 0.7294$ & $0.7443 \; - \; 0.7579$ & $7.8451\;(1.2388)$ 
      & $0.3018 \; - \; 0.2970$ & $0.3416 \; - \; 0.3478$ & $6.3912\;(0.6021)$  \\
\hline\hline
\end{tabular}}
\begin{minipage}{21.5cm}
\vspace{3pt}
\footnotesize Note: the column ``SGF'' refers to the estimator deduced from the Gaussian loss with SCAD and MCP penalization, respectively. The column ``SLSF'' refers to the estimator deduced from the least squares loss with SCAD and MCP penalization, respectively.\\
The MSE values for SOFAR in parenthesis are based on $\widehat{\Lambda}_\text{n,SOFAR}^\ast$, with $\widehat{\Lambda}_\text{n,SOFAR}^\ast := \Lambda^\ast_n \widehat{R}_\text{SOFAR}^\ast,\;
\widehat{R}_\text{SOFAR}^\ast 
:=
\underset{R \in \mathcal{O}(m)}{\arg\min}\|\widehat{\Lambda}_{n,\text{SOFAR}} - \Lambda^\ast_n R\|^2_F$.
\end{minipage}
\end{table}
\end{landscape}

\begin{landscape}
\begin{table}[H]\centering\caption{Model selection and precision accuracy based on 200 replications, DGP (\ref{dgp1}), arbitrary sparse structure and $\Psi^\ast_n$ non-diagonal  (sparsity pattern (v)) with respect to $(n,p_n,m)$. \label{support_arbitrary_misspec}}
\scalebox{0.9}{\begin{tabular}{cc|ccc|ccc}\hline\hline
 &    & \multicolumn{3}{c|}{$n=250$}       & \multicolumn{3}{c}{$n=500$}\\ 
\multicolumn{2}{c|}{$(p_n,m)$ }  & SGF & SLSF & SOFAR & SGF & SLSF & SOFAR \\ \hline
 & C1 & $84.66 \; - \; 89.67$ & $94.46 \; - \; 96.35$ & $93.93$ 
      & $89.61 \; - \; 93.31$ & $96.14 \; - \; 97.35$ & $93.54$  \\
\scalebox{0.95}{$(60,3)$}
& C2 & $100 \; - \; 100$ & $99.98 \; - \; 99.98$ & $92.91$ 
      & $100 \; - \; 100$ & $100 \; - \; 100$ & $94.44$  \\
 & MSE & $0.2592 \; - \; 0.2544$ & $0.2644 \; - \; 0.2713$ & $3.7779\;(0.4345)$ 
      & $0.1165 \; - \; 0.1138$ & $0.1284 \; - \; 0.1297$ & $3.0521\;(0.2431)$  \\
\hline      
& C1 & $84.36 \; - \; 89.64$ & $94.11 \; - \; 96.16$ & $92.58$ 
      & $90.56 \; - \; 93.48$ & $96.63 \; - \; 97.67$ & $92.46$  \\
\scalebox{0.95}{$(60,4)$}
& C2 &  $100 \; - \; 100$ & $99.94 \; - \; 99.93$ & $91.39$ 
      & $100 \; - \; 100$ & $100 \; - \; 100$ & $92.35$  \\
& MSE & $0.3534 \; - \; 0.3485$ & $0.4075 \; - \; 0.4206$ &  $7.4735\;(0.7177)$
      & $0.1552 \; - \; 0.1527$ & $0.1790 \; - \; 0.1846$ & $6.5640\;(0.3788)$  \\
\hline            
& C1 & $83.25 \; - \; 88.47$ & $94.93 \; - \; 96.72$ & $94.18$ 
      & $92.19 \; - \; 94.83$ & $96.76 \; - \; 97.74$ & $93.62$   \\
\scalebox{0.95}{$(120,3)$}
& C2 & $100 \; - \; 100$ & $100 \; - \; 100$ & $93.27$ 
      & $100 \; - \; 100$ & $100 \; - \; 100$ & $94.78$  \\
& MSE & $0.4862 \; - \; 0.4782$ & $0.4760 \; - \; 0.4920$ & $6.9414\;(0.8543)$ 
      & $0.1981 \; - \; 0.1947$ & $0.2363 \; - \; 0.2413$ & $5.3025\;(0.4199)$  \\
\hline      
                              
& C1 & $83.69 \; - \; 88.99$ & $95.01 \; - \; 96.89$ & $93.85$ 
      & $92.34 \; - \; 95.24$ & $96.79 \; - \; 97.86$ & $93.42$  \\
\scalebox{0.95}{$(120,4)$}
& C2 & $100 \; - \; 100$ & $100 \; - \; 100$ &  $89.72$
      & $100 \; - \; 100$ & $100 \; - \; 100$ & $92.81$  \\
& MSE & $0.6373 \; - \; 0.6264$ & $0.7082 \; - \; 0.7272$ & $13.5563\;(1.3894)$ 
      & $0.2805 \; - \; 0.2758$ & $0.3499 \; - \; 0.3582$ & $10.6029\;(0.7078)$  \\
\hline    
& C1 & $82.01 \; - \; 87.43$ & $95.43 \; - \; 97.16$ & $94.75$ 
      & $92.31 \; - \; 94.95$ & $96.90 \; - \; 97.94$ & $94.01$  \\
\scalebox{0.95}{$(180,3)$}
& C2 & $100 \; - \; 100$ & $100 \; - \; 100$ & $93.69$ 
      & $100 \; - \; 100$ & $100 \; - \; 100$ & $95.68$  \\
& MSE & $0.7208 \; - \; 0.7178$ & $0.7267 \; - \; 0.7481$ & $8.2446\;(1.1955)$ 
      & $0.3034 \; - \; 0.2977$ & $0.3569 \; - \; 0.3624$ & $6.3573\;(0.6182)$  \\
\hline\hline
\end{tabular}}
\begin{minipage}{21.5cm}
\vspace{3pt}
\footnotesize Note: the column ``SGF'' refers to the estimator deduced from the Gaussian loss with SCAD and MCP penalization, respectively. The column ``SLSF'' refers to the estimator deduced from the least squares loss with SCAD and MCP penalization, respectively.\\
The MSE values for SOFAR in parenthesis are based on $\widehat{\Lambda}_\text{n,SOFAR}^\ast$, with $\widehat{\Lambda}_\text{n,SOFAR}^\ast := \Lambda^\ast_n \widehat{R}_\text{SOFAR}^\ast,\;
\widehat{R}_\text{SOFAR}^\ast 
:=
\underset{R \in \mathcal{O}(m)}{\arg\min}\|\widehat{\Lambda}_{n,\text{SOFAR}} - \Lambda^\ast_n R\|^2_F$.
\end{minipage}
\end{table}
\end{landscape}

\begin{landscape}
\begin{table}[H]\centering\caption{Model selection and precision accuracy based on 200 replications, DGP (\ref{dgp2}), arbitrary sparse structure and $\Psi^\ast_n$ non-diagonal  (sparsity pattern (v)) with respect to $(n,p_n,m)$. \label{support_arbitrary_misspec_ts}}
\scalebox{0.9}{\begin{tabular}{cc|ccc|ccc}\hline\hline
 &    & \multicolumn{3}{c|}{$n=250$}       & \multicolumn{3}{c}{$n=500$}\\ 
\multicolumn{2}{c|}{$(p_n,m)$ }  & SGF & SLSF & SOFAR & SGF & SLSF & SOFAR \\ \hline         
     & C1 & $84.90 \; - \; 89.61$ & $91.40 \; - \; 93.87$ & $93.36$ 
      & $87.41 \; - \; 90.87$ & $88.96 \; - \; 90.73$ & $93.26$  \\
\scalebox{0.95}{$(60,3)$}
& C2 & $100 \; - \; 100$ & $98.98 \; - \; 98.93$ &  $92.37$
      & $100 \; - \; 100$ & $99.78 \; - \; 99.76$ & $94.83$  \\
 & MSE & $0.5828 \; - \; 0.5859$ & $1.1727 \; - \;1.0879$ &  $4.3293\;(0.7238)$
      & $0.2445 \; - \; 0.2564$ & $0.4521 \; - \; 0.4369$ & $3.6107\;(0.4049)$  \\
\hline          
& C1 & $84.80 \; - \; 89.59$ & $89.75 \; - \; 92.11$ & $92.83$ 
      & $90.25 \; - \; 93.18$ & $90.20 \; - \; 92.35$ & $91.88$  \\
\scalebox{0.95}{$(60,4)$}
& C2 & $100 \; - \; 100$ & $99.14 \; - \; 98.82$ & $89.43$ 
      & $100 \; - \; 100$ & $99.54 \; - \; 99.37$ & $92.82$  \\
& MSE & $0.8089 \; - \; 0.8162$ & $1.3863 \; - \; 1.4845$ & $8.7268\;(1.2287)$ 
      & $0.3473 \; - \; 0.3494$ & $0.8278 \; - \; 0.9286$ & $7.0553\;(0.6235)$  \\
\hline         
    
& C1 & $84.07 \; - \; 89.52$ & $90.48 \; - \; 93.10$ &  $94.36$
      & $92.17 \; - \; 94.74$ & $91.20 \; - \; 93.21$ & $93.62$  \\
\scalebox{0.95}{$(120,3)$}
& C2 & $100 \; - \; 100$ & $99.88 \; - \; 99.69$ & $92.77$ 
      & $100 \; - \; 100$ & $99.83 \; - \; 99.84$ & $95.23$  \\
& MSE & $1.2823 \; - \; 1.3134$ & $1.6126 \; - \; 1.7297$ & $7.9089\;(1.4714)$ 
      & $0.5463 \; - \; 0.5665$ & $1.0274 \; - \; 1.0343$ & $5.8175\;(0.7627)$  \\
\hline      
& C1 & $83.58 \; - \; 89.12$ & $90.85 \; - \; 93.45$ & $93.61$ 
      & $92.84 \; - \; 95.53$ & $93.07 \; - \; 94.81$ & $92.78$  \\
\scalebox{0.95}{$(120,4)$}
& C2 & $100 \; - \; 100$ & $99.79 \; - \; 99.55$ & $90.65$ 
      & $100 \; - \; 100$ & $99.88 \; - \; 99.81$ & $91.97$  \\
& MSE & $1.5176 \; - \; 1.5428$ & $2.2719 \; - \; 2.3588$ & $13.4791\;(2.4224)$ 
      & $0.5598 \; - \; 0.5742$ & $1.0539 \; - \; 1.1916$ & $12.8153\;(1.2390)$  \\
\hline   
     
& C1 & $82.39 \; - \; 87.79$ & $91.58 \; - \; 93.82$ & $94.62$ 
      & $91.85 \; - \; 94.84$ & $92.46 \; - \; 94.31$ & $94.00$  \\
\scalebox{0.95}{$(180,3)$}
& C2 & $100 \; - \; 100$ & $99.79 \; - \; 99.51$ &  $93.61$
      & $100 \; - \; 100$ & $99.91 \; - \; 99.86$ & $94.61$  \\
& MSE & $1.8582 \; - \; 1.9028$ & $2.3531 \; - \; 2.4386$ & $9.8763\;(2.2220)$ 
      & $0.8276 \; - \; 0.8543$ & $1.1838 \; - \; 1.1708$ & $8.0994\;(1.1578)$  \\
\hline\hline
\end{tabular}}
\begin{minipage}{21.5cm}
\vspace{3pt}
\footnotesize Note: the column ``SGF'' refers to the estimator deduced from the Gaussian loss with SCAD and MCP penalization, respectively. The column ``SLSF'' refers to the estimator deduced from the least squares loss with SCAD and MCP penalization, respectively.\\
The MSE values for SOFAR in parenthesis are based on $\widehat{\Lambda}_\text{n,SOFAR}^\ast$, with $\widehat{\Lambda}_\text{n,SOFAR}^\ast := \Lambda^\ast_n \widehat{R}_\text{SOFAR}^\ast,\;
\widehat{R}_\text{SOFAR}^\ast 
:=
\underset{R \in \mathcal{O}(m)}{\arg\min}\|\widehat{\Lambda}_{n,\text{SOFAR}} - \Lambda^\ast_n R\|^2_F$.
\end{minipage}
\end{table}
\end{landscape}

\section{Real data applications}\label{real_data}

\subsection{Portfolio allocation}\label{port_alloc}

In this section, we assess the variance-covariance forecast performances of multiple models using financial data, where we rank the forecast accuracy by the model confidence set (MCS) of \cite{hansen2011}, which provides a testing framework for the null hypothesis of forecast equivalence across subsets of models. We first start with a description of the data. Hereafter, we do not index by $n$ the model parameters and $p$ to ease the notations. 

\subsubsection{Data}

We consider hereafter the stochastic process $(r_t)_{t\in\Zb}$ in $\Rb^p$ of the log-stock returns, where $r_{t,j} = 100\times\log(P_{t,j}/P_{t-1,j}), 1 \leq j \leq p$ with $P_{t,j}$ the stock price of the $j$-th index at time $t$. We study three financial portfolios: a low-dimensional portfolio of daily log-returns composed with the MSCI stock index based on the sample December 1998-March 2018, which yields a sample size $n = 5006$, and for $p=23$ countries (Australia, Austria, Belgium, Canada, Denmark, Finland, France, Germany, Greece, Hong Kong, Ireland, Italy, Japan, Netherlands, New Zealand, Norway, Portugal, Singapore, Spain, Sweden, Switzerland, the United-Kingdom, the United-States); 
a mid-dimensional portfolio of daily log-returns based on the S\&P 100, consisting of firms that have been continuously included in the index over the period February 2010 -- January 2020, excluding AbbVie Inc., Dow Inc., General Motors, Kraft Heinz, Kinder Morgan and PayPal Holdings, which leaves $p=94$ with a sample size $n=2500$; 
a high-dimensional portfolio of daily log-returns based on the S\&P 500, consisting of the firms that have been continuously included in the index over the period September 2014 -- January 2022, excluding Etsy Inc., SolarEdge Technologies Inc., PayPal, Hewlett Packard Enterprise, Under Armour (class C), Fortive, Lamb Weston, Ingersoll Rand Inc., Ceridian HCM, Linde PLC, Moderna, Fox Corporation (class A and B), Dow Inc., Corteva Inc., Amcor, Otis Worldwide, Carrier Global, Match Group, Viatris Inc., which leaves $p=480$ assets with a sample size $n=1901$.  Hereafter, the portfolios are denoted by MSCI, S\&P 100 and S\&P 500, respectively\footnote{The MSCI, S\&P 100 and S\&P 500 data can be found on https://www.msci.com/, https://finance.yahoo.com and https://macrobond.com, respectively.}.



\subsubsection{Global Minimum Variance Portfolio (GMVP)}

The economic performances are assessed through the GMVP investment strategy. The latter problem at time $t$, in the absence of short-sales constraints, is defined as
\begin{equation}\label{portprob}
\min_{w_t} \; w^{\top}_t \; \widehat{\Sigma}_t \; w_t, \;\;\text{s.t.} \;\;\iota^{\top} w_t = 1,
\end{equation}
where $\widehat{\Sigma}_t$ is the predicted $p \times p$ variance-covariance matrix and $\iota$ is a $p \times 1$ vector of $1$'s. The explicit solution is given by $\omega_t = \widehat{\Sigma}^{-1}_t \iota/\iota^{\top}\widehat{\Sigma}^{-1}_t\iota$: as a function depending only on $\widehat{\Sigma}_t$, the GMVP performance essentially depends on the precise measurement of the variance-covariance matrix. The following out-of-sample performance metrics (annualized) will be reported: \textbf{AVG} as the average of the out-of-sample portfolio returns, multiplied by $252$; \textbf{SD} as the standard deviation of the out-of-sample portfolio returns, multiplied by $\sqrt{252}$; \textbf{IR} as the information ratio computed as $\textbf{AVG}/\textbf{SD}$. The key performance measure is the out-of-sample \textbf{SD}. The GMVP problem essentially aims to minimize the variance rather than to maximize the expected return. Hence, as emphasized in \cite{engle2019}, Section 6.2, high \textbf{AVG} and \textbf{IR} are desirable but should be considered of secondary importance compared with the quality of the measure of a variance-covariance matrix estimator. 

\subsubsection{Competing variance-covariance matrix estimators}

We will rank the following models according to the MCS procedure, where, excluding the DCC, the competing models provide fixed estimators of $\Sigma_t$, i.e., $\widehat{\Sigma}_t=\widehat{\Sigma}$:
\begin{itemize}
\item \textbf{Sample-VC}: the sample estimator of the variance-covariance. 
\item \textbf{covM1}: the predicted $\widehat{\Sigma}$ shrunk towards a one-factor market model, where the factor is defined as the cross-sectional average of all the components of $r_t$. The idiosyncratic volatility of the residuals allows the target to preserve the diagonal of the sample covariance matrix. This estimator was developed by \cite{ledoit2003}.
\item \textbf{GIS}: the geometric-inverse shrinkage estimator $\widehat{\Sigma}$, which is a nonlinear shrinkage estimator based on the symmetrized Kullback-Leibler loss: this static estimator can be considered as geometrically averaging the linear-inverse shrinkage (LIS) estimator with the quadratic-inverse shrinkage (QIS) estimator of \cite{ledoit2022}.
\item \textbf{SAF}: the factor model $\widehat\Sigma = \widehat\Lambda\widehat\Lambda^\top + \widehat\Psi$, where $\Psi$ is non-diagonal but sparse with bounded eigenvalues. This is the sparse approximate factor (SAF) model of \cite{bai2016} estimated by adaptive LASSO Gaussian QML under $\Lambda^\top \Psi^{-1}\Lambda$ diagonal. More details on the SAF and implementation are provided in Subsection \ref{safm} of the Appendix. 
\item \textbf{SGF}: the factor model $\widehat\Sigma = \widehat\Lambda\widehat\Lambda^\top+\widehat\Psi$, estimated by criterion (\ref{stat_crit}) with the Gaussian loss for a general arbitrary sparse $\Lambda$, $\Psi$ diagonal, with SCAD (SGF$^{\text{scad}}$) and MCP (SGF$^{\text{mcp}}$). 
\item \textbf{SLSF}: the factor model $\widehat\Sigma = \widehat\Lambda\widehat\Lambda^\top+\widehat\Psi$, estimated by criterion (\ref{stat_crit}) with the least squares loss for a general arbitrary sparse $\Lambda$, $\Psi$ diagonal with SCAD (SLSF$^{\text{scad}}$) and MCP (SLSF$^{\text{mcp}}$). 
\item \textbf{DCC}: the scalar DCC of \cite{engle2002} with GARCH(1,1) dynamics for the marginals. This is a dynamic estimator of $\Sigma_t$: once the model is estimated in-sample, the estimated parameters are plugged in the out-of-sample period to generate the out-of-sample sequence of $\Sigma_t$, providing a sequence of GMVP weights. The estimation is carried out by a standard two-step Gaussian QML. More details on the DCC and estimation, in particular the composite likelihood method, are provided in Subsection~\ref{dcc_model} of the Appendix.
\end{itemize}
It is worthwhile to include variance-covariance models that are not derived from a factor model, thereby motivating the inclusion of the scalar DCC in the analysis. 

As in Section \ref{simulations}, we set $a_{\text{scad}}=3.7$, $b_{\text{mcp}}=3.5$. The selection of an optimal $\gamma_n$ is based on the out-of-sample cross-validation procedure detailed in Subsection \ref{CV} of the Appendix. As emphasized in \cite{deNard2021}, there is no consensual selection procedure of the number of factors so that it is of interest to assess the sensitivity of the factor-model-based allocation performances with respect to $m$: we set $m \in \{1,2,3,4,5\}$. It is worth noting that the SOFAR procedure estimates a sparse loading matrix only and does not allow to obtain an estimator of the variance-covariance matrix: the SOFAR procedure is thus not under consideration. 
The economic performances are ranked out-of-sample by the MCS test, which takes the distance based on the difference of the squared returns of portfolios $i$ and $j$, defined as
$u_{ij,t} = \big(\mu_{i,t}-\overline{\mu}_i\big)^2 - \big(\mu_{j,t}-\overline{\mu}_j\big)^2$, with $\mu_{k,t}=w^{\top}_{k,t} y_t$ the portfolio return at time $t$, where $w_{k,t}$ represents the GMVP weight deduced from variance-covariance model $k$, and $\overline{\mu}_k$ is the average portfolio return over the period. The MCS test is evaluated at the $10\%$ level based on the range statistic and with block bootstrap with $10,000$ replications: see \cite{hansen2003} for further technical procedures. We report below
the out-of-sample periods used for the test accuracy and the number of factors $\widehat{m}$ selected in-sample by the eigenvalue-based procedure of \cite{onatski2010} on an indicative basis: 

\begin{itemize}
\item MSCI portfolio: in sample is 01/01/1999 -- 07/01/2010 with $\widehat{m}=1$ and out-of-sample is 07/02/2010 -- 03/12/2018.
\item S\&P 100 portfolio: in sample is 02/19/2010 -- 04/22/2014 with $\widehat{m}=3$ and out-of-sample is 04/23/2014 -- 01/23/2020.
\item S\&P 500 portfolio: in sample is 09/25/2014 -- 12/12/2018 with $\widehat{m}=1$ and out-of-sample is 12/13/2018 -- 01/27/2022.
\end{itemize}
The GMVP performance results are displayed in Table \ref{Metrics_performance}, and they can be summarized as follows - unless stated otherwise, the results are with respect to \textbf{SD} -: the sparse factor loading-based strategy enters the MCS in the low, mid and high-dimensional cases; the DCC provides the best performance for the low-dimensional portfolio only; for the mid-dimensional portfolio, the least squares-based SLSF estimator with $m=4,5$ only enters the MCS; in the high-dimensional portfolio, although the SAF seems more appropriate to measure the variance-covariance suggesting the existence of cross-correlation in $\eps_t$, the sparse factor loading SGF for $m=2,3$ still enters the MCS at $10\%$-level. 
   
\begin{table}[H]
\caption{Annualized GMVP performance metrics for various estimators.}\label{Metrics_performance}
\scalebox{0.82}{\begin{tabular}{c||cccc||cccc||cccc}\hline\hline
& \multicolumn{4}{c||}{MSCI} & \multicolumn{4}{c||}{S\&P 100} &  \multicolumn{4}{c}{S\&P 500} \\
&    \multicolumn{4}{c||}{\scalebox{0.9}{07/02/2010 - 03/12/2018}} &  \multicolumn{4}{c||}{\scalebox{0.9}{04/23/2014 -- 01/23/2020}} &  \multicolumn{4}{c}{\scalebox{0.9}{12/13/2018 - 01/27/2022}} \\
  &  \textbf{AVG} & \textbf{SD} & \textbf{IR} & \textbf{MCS}
  & \textbf{AVG} & \textbf{SD} & \textbf{IR} & \textbf{MCS} & \textbf{AVG} & \textbf{SD} & \textbf{IR} & \textbf{MCS}\\ 
\hline

    
DCC & 7.482  &  \textbf{9.551}  &  0.783   & \textbf{1.000}  &   17.293  & 11.941  &  1.448    & 0  &   15.278  & 20.146  &  0.758   &  0.024 \\
Sample-VC &  9.014 &  10.003 &   0.901  & 0  &  13.187 &  11.412 &   1.156    & 0  &  2.808 &  22.090  &  0.127   &  0.007 \\

 $\text{SGF}^{\text{scad}}_1$ &  7.629  & 11.172  &  0.683   & 0  &    9.426 &  12.617  &  0.747    & 0  &   4.288 &  21.472 &   0.200   &  0.002 \\
$\text{SGF}^{\text{mcp}}_1$ & 7.629 &  11.172  &  0.683  & 0  &   9.469 &  12.767 &   0.742    &  0 &  4.219  & 21.532  &  0.196   &  0.002 \\
$\text{SGF}^{\text{scad}}_2$ & 8.980  &  9.625 &   0.933   & \textbf{0.832}  &  11.405  & 13.309  &  0.857   & 0  &  7.904 &  17.833  &  0.443   & \textbf{0.231} \\
$\text{SGF}^{\text{mcp}}_2$ & 8.978 &   9.625  &  0.933  & \textbf{0.832}  &  11.555  & 13.265 &   0.871    & 0  &   7.914  & 17.883 &   0.443   & 0.024 \\
$\text{SGF}^{\text{scad}}_3$ & 7.891 &   9.955  &  0.793  & 0  &  15.702 &  12.403  &  1.266   &  0 &  7.124 &  17.602  &  0.405  & \textbf{0.798} \\
$\text{SGF}^{\text{mcp}}_3$ & 8.397 &   9.780  &  0.859   &  0 &  16.035  & 12.456 &   1.287    & 0  &   7.221  & 17.670  &  0.409   & \textbf{0.725} \\
$\text{SGF}^{\text{scad}}_4$ & 7.866 &   9.985  &  0.788   & 0  &  14.070  & 12.419  &  1.133    &  0 &  8.427 &  19.041  &  0.443   &  0.024 \\
$\text{SGF}^{\text{mcp}}_4$ &  8.511 &   9.809  &  0.868   & 0  &   14.774  & 12.338  &  1.197    & 0  &  7.682  & 18.745   & 0.410  & 0.024 \\
$\text{SGF}^{\text{scad}}_5$ & 9.281 &   9.977  &  0.930  &  0 &  11.373 &  12.001 &   0.948    &  0 &   8.748  & 19.438  &  0.450  & 0.024 \\
$\text{SGF}^{\text{mcp}}_5$ &  9.192  &  9.989  &  0.920  & 0  &  12.029 & 11.907  &  1.010    & 0  &   8.162 &  19.162 &   0.426   &  0.024 \\
       
 $\text{SLSF}^{\text{scad}}_1$ & 7.961 &  10.889 &   0.731   & 0  &   9.069 &  12.583  &  0.721    &  0 &    4.411 &  21.430  &  0.206  & 0.007 \\
$\text{SLSF}^{\text{mcp}}_1$ &  7.961 &  10.889  &  0.731   & 0  &   9.091 &  12.774 &   0.712   & 0  &  4.315 &  21.513 &   0.201  & 0.002 \\
$\text{SLSF}^{\text{scad}}_2$ & 8.552  &  9.730  &  0.879 &  0.002 &  9.017  &  11.823 &   0.763    &  0 &   3.677 &  21.346  &  0.172  & 0.007 \\
$\text{SLSF}^{\text{mcp}}_2$ &  9.021 &   9.655  &  0.934   & 0.018  &  9.025  & 12.041  &  0.750    & 0  &  3.645 &  21.387  &  0.170  & 0.007 \\
$\text{SLSF}^{\text{scad}}_3$ & 8.829  &  9.784  &  0.902   & 0.002  &  8.841 &  11.331  &  0.780   & 0  &  5.384 &  19.566  &  0.275  &  0.007\\
$\text{SLSF}^{\text{mcp}}_3$ & 8.858  &  9.723  &  0.911 & 0.017  &   8.882 & 11.454 &   0.776    & 0  &  5.881  & 19.300  &  0.305  & 0.007 \\
$\text{SLSF}^{\text{scad}}_4$  &  8.715 &   9.963 &   0.875   & 0  & 8.508  & 10.941  &  0.778    & \textbf{0.951}  &  4.275 &  19.195  &  0.223 &  0.007\\
$\text{SLSF}^{\text{mcp}}_4$ &  8.826 &   9.948  &  0.887   & 0  &  8.656 &  11.002 &   0.787    &  0.040 &   4.425 &  19.407 &   0.228   &  0.007\\
$\text{SLSF}^{\text{scad}}_5$  &  9.025 &   9.820  &  0.919  & 0  &  8.000 &  11.629  &  0.688    &  0 &   5.162  & 19.741  &  0.262  & 0.007 \\
$\text{SLSF}^{\text{mcp}}_5$ &  9.028 &  10.007  &  0.902   &  0 &    8.190  & \textbf{10.934}  &  0.749    & \textbf{1.000}  &   5.001  & 19.882   & 0.252  & 0.007 \\
    
$\text{SAF}_1$ & 8.063 &  10.446 &   0.772   & 0  &   11.031 &  12.137  &  0.909   & 0  &  3.660 &  22.041  &  0.166  &  0.002 \\
$\text{SAF}_2$ &  8.856 &   9.769  &  0.907  & 0  &   10.927  & 11.786  &  0.927    &  0 &  7.690  & 18.500  &  0.416  & 0.024 \\
$\text{SAF}_3$ & 8.837  &  9.902 &   0.893   & 0  &  14.595 &  11.815  &  1.235   & 0  &  8.014  & 25.121  &  0.319  &  0.007\\
$\text{SAF}_4$ & 8.679  &  9.939  &  0.873   & 0  &  14.412  & 11.759 &   1.226   &  0 &  12.330 &  19.223  &  0.641  & 0.007 \\
$\text{SAF}_5$ &  8.898  & 10.099  &  0.881  & 0  &  14.100  & 11.606 &   1.215    & 0  &   6.359 &  \textbf{17.456}  &  0.364  & \textbf{1.000} \\
     
$\text{GIS}$ & 9.004  &  9.995  &  0.901  & 0  &   13.073 &  11.285 &   1.158    &  0.040 &  5.820  & 19.849  &  0.293 &  0.007 \\
$\text{covM1}_1$ & 8.975  &  9.988  &  0.899   &  0 &   12.741 &  11.370  &  1.121    &  0.001 &  4.125  & 20.103  &  0.205  &  0.007 \\
\hline\hline
\end{tabular}}
\begin{minipage}{16.5cm}
\footnotesize Note: The lowest \textbf{SD} figure is in bold face. The index $k$ in $\text{SGF}^{\text{pen}}_k$, $\text{SLSF}^{\text{pen}}_k$ and $\text{SAF}_k$ represents the number of factors. The out-of-sample periods are indicated above \textbf{AVG}, \textbf{SD}, \textbf{IR} and \textbf{MCS}. The out-of-sample dates are reported below the portfolio name.
\end{minipage}
\end{table}

\subsection{Diffusion index data}\label{diffusion_index}

We now propose to investigate how the estimation of the factors impact the prediction accuracy. To do so, we follow the experiment of \cite{bailiao2016}, Section 5.3., which concerns the industrial production based on macroeconomic time series of the United States. It consists of a macroeconomic panel of $p=131$ series from 1959 to 2007, representing a sample of size $\mathcal{T}=528$. In the same spirit as in \cite{stock2002a} and \cite{bailiao2016}, we denote by $Y_t$ the scalar time series variable to be predicted and by $X_t$ the $p$-dimensional vector of candidate predictors, and assume that $(X_t,Y_{t})$ satisfies a factor model decomposition with $m$ common latent factors $F_t$:
\begin{equation*}
Y^h_{t+h} = \alpha_h + \beta_h F_t + \gamma_h W_{t,l} + u_{t+h}, \; X_t = \Lambda F_t + \eps_t,
\end{equation*}
where $h$ is the forecast horizon, $Y^h_{t+h} = h^{-1}\sum^h_{i=1}Y_{t+i}$ is the h-step-ahead variable to be
predicted and defined as the industrial production, $W_{t,l} = (Y_t,\ldots,Y_{t-l})^\top\in \Rb^l$ contains the lagged values of $Y_t$, and $u_{t+h}$ is the forecast error. We construct the forecasts of $Y^h_{t+h}$ according to a fixed rolling window of size $n = 422$ (that is $0.8\,\mathcal{T}$), where for each window we estimate the factor model of $X_t$ and estimate the latent factors $F_t$. Finally, we compute the forecast of $Y^h_{n+h}$ as $\widehat{Y}^h_{n+h} = \widehat{\alpha}_h + \widehat{\beta}_h \widehat{F}_n + \widehat{\gamma}_h W_{n,l}$,
where $\widehat{\alpha}_h, \widehat{\beta}_h,\widehat{\gamma}_h$ are obtained by regressing $Y^h_{t+h}$ on $\widehat{F}_t, W_{t,l}$ with intercept. The squared prediction error is defined as $(Y^h_{n+h}-\widehat{Y}^h_{n+h})^2$. 
All the data are standardized. \\
The following methods are employed to obtain $\widehat{F}_t$: PCA, SOFAR, SAF, SGF, SLSF. 
The estimator deduced from PCA follows from \cite{bai2003}, that is the estimator of the factors $\widehat{\mathbf{F}} \in \Rb^{n \times m}$ with $t$-th row $\widehat{F}^\top_t$ is $\sqrt{n}$ times eigenvectors corresponding to the $m$ largest eigenvalues of $\mathbf{X}\mathbf{X}^\top \in \Rb^{n \times n}$, and $\widehat{\Lambda} = \mathbf{X}^\top \widehat{\mathbf{F}}/n$. As for SAF, SGF,  SLSF, $\widehat{F}_t$ is obtained by the GLS estimator. Regarding SOFAR, its optimal regularization parameter was selected as the one that achieves the best performance in the simulation experiments. 
However, in the case of real data, the selection of the optimal tuning parameter that ensures the best out-of-sample performance is unfeasible. 
Therefore, from the candidate values $(e^{-2}, e^{-2.25}, \ldots, e^{+5})$, the regularization parameters are determined using BIC, as recommended by the authors. The out-of-sample prediction performances are assessed through the mean squared out-of-sample
forecasting error (MSE) relative to the PCA method for $h=12,24$, $l=1,3$ and $m=5,6,7$. The choice for $m$ follows from the setting of \cite{bailiao2016}, who considered $m=7,8$: we preferred to consider a slightly lower number of factors for parsimony. The mean squared MSE is defined as $s^{-1}\sum^{s-1}_{r=0}(Y^h_{r+n+h}-\widehat{Y}^h_{r+n+h})^2$, with $s=\mathcal{T}-n-h$. Table \ref{MSE_out_of_sample_IPdata} reports the MSE metrics
relative to the PCA method. It is observed that our proposed method outperforms SOFAR, excluding $h=12,l=3,m=6$ only. We notice SAF has the best performances for $m=5$, but when $m=6,7$, the sparse factor model is beneficial for prediction. We also observe that in both period-ahead forecasts $h$, $m=5$ is much better than $m=6,8$, for SAF, SGF and SLSF.  

\begin{table}[H]\centering\caption{Relative MSE for out-of-sample predictions. \label{MSE_out_of_sample_IPdata}}
\scalebox{0.85}{\begin{tabular}{cc|cccccc}\hline\hline
\multicolumn{2}{c|}{$(h,l,m)$ }   & SOFAR & SAF & $\text{SGF}^{\text{scad}}$ & $\text{SGF}^{\text{mcp}}$ & $\text{SLSF}^{\text{scad}}$ & $\text{SLSF}^{\text{mcp}}$ \\ \hline      
     & $m=5$  & 0.9267  & 0.7057 & 0.7317 & 0.7881 & 0.8799 &  0.8537 \\
$(12,1,m)$
& $m=6$  & 0.8811  & 0.9564  & 0.9448 & 0.9417 & 0.9191 & 0.8758 \\
 & $m=7$  & 0.9224 & 1.0733 & 0.7985 & 0.8961 & 1.0090 & 1.0173 \\
\hline          

     & $m=5$ & 0.9304 & 0.7136 & 0.7509 & 0.7979 & 0.9024 & 0.8734 \\
$(12,3,m)$
& $m=6$ & 0.8690 & 0.9588 & 0.9432 & 0.9453 & 0.9214 & 0.8787 \\
 & $m=7$ & 0.9176 & 1.0826  & 0.8106  & 0.9116 & 1.0254 &  1.0309 \\  \hline

     & $m=5$ & 0.9257 & 0.6987 & 0.7629  & 0.7385 & 0.7699 & 0.7695 \\
$(24,1,m)$
& $m=6$ & 0.9544 & 0.8958 & 0.8130  & 0.8510 & 0.7857 & 0.7919 \\
 & $m=7$ & 0.9230 & 0.8671  & 0.8672 & 0.8144 & 0.9118 & 0.9179 \\\hline

      & $m=5$ & 0.9296 &  0.7051 & 0.7658 & 0.7484 & 0.7853 & 0.7841 \\
$(24,3,m)$
& $m=6$ & 0.9542 &  0.8973 & 0.8172 & 0.8700 & 0.7901 & 0.7977 \\
 & $m=7$ & 0.9252 & 0.8599 & 0.8630 & 0.8103 & 0.9191 & 0.9305 \\
\hline\hline
\end{tabular}}
\end{table}                       

\section{Discussion and conclusion}

We study the asymptotic properties of sparse factor models. We establish the sparsistency property of the corresponding penalized M-estimator, where the identification condition allows for a broad range of sparsity patterns. 

Throughout the paper, we worked under a diagonal variance-covariance of the idiosyncratic error variables. An interesting direction would consist in replacing this condition by a non-diagonal but sparse variance-covariance matrix within the sparse factor loading framework. This would raise issues regarding the conditions for identification to disentangle the sparsity patterns of $\Lambda_n$ and $\Psi_n$. We shall leave it for the future research.

\bigskip

\noindent{\large\bf Acknowledgments}

\bigskip

This work was supported by JSPS KAKENHI Grant (JP22K13377 to BP; JP20K19756 and JP20H00601 to YT). The authors warmly thank Jan Magnus for his suggestions.

\bigskip

\newpage
\bibliography{biblio}

\newpage

\appendix
\numberwithin{equation}{section}
\makeatletter 
\newcommand{\section@cntformat}{Appendix \thesection\ }
\makeatother

\section{Preliminary results}\label{technical_appendix}

In this section, we provide three deviation bounds on the empirical variance-covariance matrix of $\eps_t,F_t$ and the covariance $\eps_tF^\top_t$.

\mds

\begin{lemma}\label{deviation_ineq}
Suppose Assumptions \ref{factor_assumption_1}--\ref{factor_assumption_2} of the main text are satisfied. Let $\rho^{-1} = 3 r^{-1}_1 + 1.5r^{-1}_2+r^{-1}_3+1$ and assume $\log(p_n)^{6/\rho}=o(n)$. Then, the following bounds are satisfied:
\begin{itemize}
    \item[(i)] $\underset{1 \leq k \leq m, 1 \leq l \leq p_n}{\max}|\frac{1}{n}\overset{n}{\underset{t=1}{\sum}} F_{t,k}\eps_{t,l} | = O_p(\sqrt{\log(p_n)/n})$.
    \item[(ii)] $\underset{1 \leq k,l \leq m}{\max}|\frac{1}{n}\overset{n}{\underset{t=1}{\sum}} F_{t,k}F_{t,l} - \Eb[F_{t,k}F_{t,l}]| =  O_p(\sqrt{1/n})$.
    \item[(iii)]  $\underset{1 \leq k,l \leq p_n}{\max}|\frac{1}{n}\overset{n}{\underset{t=1}{\sum}} \eps_{t,k}\eps_{t,l} - \Psi_{n,kl}| =  O_p(\sqrt{\log(p_n)/n})$.
\end{itemize}
\end{lemma}
\begin{proof}
The proofs of these deviation inequalities follow from Lemma A.3 and Lemma B.1 of \cite{fan2011}. 
\end{proof}


\section{Proofs}\label{appendix_proofs}


In this section, we provide the proofs of Theorem \ref{Theorem_existence_consistent} and Theorem \ref{sparsistency}, where $\widehat{\theta}_n$ is defined in (\ref{stat_crit}) and $\Lb_n(\theta)=\sum^n_{t=1}\ell_n(X_t;\theta)$. We will write $\widehat{S}_n=n^{-1}\sum^n_{t=1}X_tX^\top_t$.

\subsection{Proof of Theorem \ref{Theorem_existence_consistent}}

We first provide the proof when $\ell_n(X_t;\theta_n)$ is the Gaussian loss. We then consider the least squares case.

\mds

\noindent\textbf{\emph{Gaussian loss.}} Let $\ell_n(X_t;\theta_n)= \text{tr}\big(X_tX^\top_t \Sigma^{-1}_n\big) + \log(|\Sigma_n|)$, $\Sigma_n = \Lambda_n\Lambda^\top_n+\Psi_n$. Note that $n^{-1}\sum^n_{t=1}\ell_n(X_t;\theta_n) = \text{tr}(\widehat{S}_n\Sigma^{-1}_n)+\log(|\Sigma_n|)$. Define the symmetric positive-definite matrix 
$$V_n = \Lambda^\ast_n R^\top_{u} + R_{u} \Lambda^{\ast\top}_n + R_u R^\top_u + D_{u}, \; R_{u}=(r_{u,1},\ldots,r_{u,p_n})^\top \in \Rb^{p_n\times m}, \; D_{u} \in \Rb^{p_n\times p_n},$$ 
with $R_u$ a fixed matrix satisfying $ \forall 1 \leq j \leq p_n, \|r_{u,j}\|_2\leq \kappa <\infty$ and $D_{u}$ a fixed diagonal matrix with positive and bounded diagonal terms satisfying $\forall 1 \leq k \leq p_n, \kappa^{-2} \leq D_{u,kk} \leq \kappa^2$. Under the sparsity assumption, $\|\Lambda^\ast_n\|_F = O(s^{1/2}_n)$.
Define the re-scaled positive-definite and symmetric matrix $U_n = V_n/s_n^{1/2}=:L_{n,u}+R_u R^\top_u/s_n^{1/2} + D_{u}/s_n^{1/2}$ and let $\Delta_{n,u} = \alpha_n U_n$ with $\alpha_n = p_n\sqrt{s_n\log(p_n)/n}+\sqrt{p_n}A_n$. We aim to prove that for any $\delta>0$, there exists $C_{\delta}>0$ large and finite such that
\begin{equation*}
\Pb(\alpha^{-1}_n\|\widehat{\Sigma}_n-\Sigma^\ast_n\|_F>C_\delta) < \delta.
\end{equation*}
Define $\Lb^{\text{pen}}_n(\theta_n) = \Lb_n(\theta_n) + n\, \sum_{k=1}^{p_nm}p(|\theta_{n\Lambda,k}|,\gamma_n)$, we have
\begin{equation}\label{obj_prob}
\Pb(\alpha^{-1}_n\|\widehat{\Sigma}_n-\Sigma^\ast_n\|_F>C_\delta)
\leq \Pb(\exists U_n \in \Ac_n: \Lb^{\text{pen}}_n(\theta^\ast_n+\text{vec}(\Delta_{n,u})) \ge  \Lb^{\text{pen}}_n(\theta^\ast_n)),
\end{equation}
with $\Ac_n:=\{U_n:\|\Delta_{n,u}\|^2_F=\alpha^2_n C^2_{\delta} \}$. If the right-hand side of (\ref{obj_prob}) is smaller than $\eps$, this implies that there is a local minimum in $\{\Sigma^\ast_n+\Delta_{n,u}: \|\Delta_{n,u}\|^2_F \leq \alpha^2_nC^2_{\delta}\}$ such that $\|\widehat{\Sigma}_n-\Sigma^\ast_n\|_F=O_p(\alpha_n)$ for $n$ large enough, because $\Sigma^\ast_n+\Delta_{n,u}$ is positive-definite. Indeed, we have
\begin{equation*}
\lambda_{\min}(\Sigma^\ast_n+\Delta_{n,u}) \geq \lambda_{\min}(\Sigma^\ast_n)+\lambda_{\min}(\Delta_{n,u}) \geq \lambda_{\min}(\Sigma^\ast_n)-\|\Delta_{n,u}\|_F>0,
\end{equation*}
because $\Sigma^\ast_n$ is positive-definite by construction and $\|\Delta_{n,u}\|_F=O(\alpha_n)=o(1)$. The following expansion holds:
\begin{eqnarray*}
\lefteqn{\frac{1}{n}\Lb^{\text{pen}}_n(\theta^\ast_n+\text{vec}(\Delta_{n,u}))-\frac{1}{n}\Lb^{\text{pen}}_n(\theta^\ast_n)}\\
& \geq & \text{tr}(\widehat{S}_n(\Sigma^\ast_n+\Delta_{n,u})^{-1})+\log|\Sigma^\ast_n+\Delta_{n,u}| - \text{tr}(\widehat{S}_n(\Sigma^\ast_n)^{-1})-\log|\Sigma^\ast_n| \\
& & + \sum_{k \in \Sc_n}\Big(p(|\theta^\ast_{n\Lambda}+\alpha_n\text{vec}(L_{n,u}+R_uR^\top_u)|,\gamma_n)-p(|\theta^\ast_{n\Lambda}|,\gamma_n)\Big)\\
& = & \nabla_{\text{vec}(\Sigma_n)^\top}\frac{1}{n}\Lb_n(\theta^\ast_n) \text{vec}(\Delta_{n,u}) +\text{vec}(\Delta_{n,u})^\top \int^1_0 g(v,\Sigma_{n,v})(1-v)\text{d}v \, \text{vec}(\Delta_{n,u}) \\
& + & \underset{k \in \Sc_n}{\sum} \Big(\alpha_n \partial_{1}p(|\theta^\ast_{n\Lambda}|,\gamma_n)\text{vec}(L_{n,u}+R_uR^\top_u)_k + \alpha^2_n\partial^2_{11}p(|\theta^\ast_{n\Lambda}|,\gamma_n)\text{vec}(L_{n,u}+R_uR^\top_u)^2_k\{1+o(1)\}\Big)\\
& =: & K_1+K_2+K_3,
\end{eqnarray*}
with $\Sigma_{n,v} = \Sigma^\ast_n+ v \Delta_{n,u}$ and
\begin{equation*}
g(v,\Sigma_{n,v}) = \Sigma^{-1}_{n,v} \otimes \Sigma^{-1}_{n,v} \widehat{S}_n \Sigma^{-1}_{n,v} + \Sigma^{-1}_{n,v}\widehat{S}_n\Sigma^{-1}_{n,v} \otimes \Sigma^{-1}_{n,v} - \Sigma^{-1}_{n,v} \otimes \Sigma^{-1}_{n,v},
\end{equation*}
which is obtained by the identification of the Hessian matrix using, e.g., Theorem 18.6 of \cite{magnus2019}. We want to prove
\begin{equation}\label{obj_bound}
\Pb(\exists U_n \in \Ac_n: nK_1+nK_2+nK_3 \ge  0 ) < \delta.
\end{equation}
Hereafter, to simplify the notations, we divide by $n$ in (\ref{obj_bound}).
Take $K_1$:
$$K_1 = \text{vec}(\Sigma^{\ast-1}_n(\Sigma^\ast_n-\widehat{S}_n)\Sigma^{\ast-1}_n)^\top \text{vec}(\Delta_{n,u}) = \text{tr}((\Sigma^\ast_n-\widehat{S}_n)\Sigma^{\ast-1}_n(\Delta_{n,u})\Sigma^{\ast-1}_n).$$
$K_1$ can be upper bounded as
\begin{eqnarray*}
|K_1| = |\underset{ij}{\sum}\big[\Sigma^\ast_n-\widehat{S}_n\big]_{ij}(\Sigma^{\ast-1}_n \Delta_{n,u}\Sigma^{\ast-1}_n)_{ij}|\leq p_n\underset{1 \leq ij \leq p_n}{\max}|\big(\widehat{S}_n-\Sigma^\ast_n\big)_{ij}|\|\Sigma^{\ast-1}_n\|^2_s \|\Delta_{n,u}\|_F.
\end{eqnarray*}
We have the decomposition:
\begin{equation*}
\widehat{S}_n - \Sigma^\ast_n = \Lambda^\ast_n \widehat{S}_{FF} \Lambda^{\ast\top}_n + \Lambda^\ast_n \widehat{S}_{n,F\eps} + \widehat{S}_{n,\eps F}\Lambda^{\ast\top}_n + \widehat{S}_{n,\eps\eps} - \Lambda^\ast_n\Lambda^{\ast\top}_n - \Psi^\ast_n,
\end{equation*}
with $\widehat{S}_{FF}=n^{-1}\sum^n_{t=1}F_tF^\top_t$, $\widehat{S}_{n,F\eps} = n^{-1}\sum^n_{t=1}F_t\eps^\top_t$, $\widehat{S}_{n,\eps F} = \widehat{S}^\top_{n,F\eps}$, $\widehat{S}_{n,\eps\eps} = n^{-1}\sum^n_{t=1}\eps_t\eps^\top_t$. Now denoting by $\mathbf{e}_k$ the $p_n$-dimensional zero column vector except for the $k$-th component being one, we obtain
{\footnotesize{\begin{eqnarray*}
\lefteqn{\|\widehat{S}_n-\Sigma^\ast_n\|_{\max}}\\
&\leq &2\underset{1 \leq k,l \leq p_n}{\max}\|\mathbf{e}^\top_k\Lambda^\ast_n \Big(\frac{1}{n}\overset{n}{\underset{t=1}{\sum}}F_t\eps^\top_t\Big)\mathbf{e}_l\|_2 + \underset{1 \leq k,l\leq p_n}{\max}\Big|\Big(\frac{1}{n}\overset{n}{\underset{t=1}{\sum}} \eps_t\eps^\top_t -\Psi^\ast_n\Big)_{kl}\Big| + \underset{1 \leq k,l\leq p_n}{\max}\|\mathbf{e}^\top_k\Lambda^\ast_n\Big( \frac{1}{n}\overset{n}{\underset{t=1}{\sum}}F_tF^\top_t - I_m\Big)\Lambda^{\ast\top}_n\mathbf{e}_l\|_2\\
& \leq & 2\,\underset{1 \leq k\leq p_n}{\max}\|\lambda^\ast_{k}\|_2 \, \sqrt{m}\,\underset{1 \leq k\leq m,1\leq l \leq p_n}{\max}
\Big|\frac{1}{n}\overset{n}{\underset{t=1}{\sum}}F_{t,k}\eps_{t,l}\Big| + \|\frac{1}{n}\overset{n}{\underset{t=1}{\sum}} \eps_t\eps^\top_t -\Psi^\ast_n\|_{\max} + m^2\,\|\Lambda^\ast_n\|^2_{\max}\|\frac{1}{n}\overset{n}{\underset{t=1}{\sum}}F_tF^\top_t - \Eb[F_iF^\top_i]\|_{\max}.
\end{eqnarray*}}}
Therefore, under Assumption \ref{assumption_parameters}, we deduce
\begin{eqnarray*}
\lefteqn{|K_1| = |\underset{ij}{\sum}\big[\Sigma^\ast_n-\widehat{S}_n\big]_{ij}(\Sigma^{\ast-1}_n\Delta_{n,u}\Sigma^{\ast-1}_n)_{ij}|}\\
&\leq & p_n\Bigg(2\,\underset{1 \leq k\leq p_n}{\max}\|\lambda^\ast_{k}\|_2 \, \sqrt{m}\,\underset{1 \leq k\leq m,1\leq l \leq p_n}{\max}
\Big|\frac{1}{n}\overset{n}{\underset{t=1}{\sum}}F_{t,k}\eps_{t,l}\Big| + \|\frac{1}{n}\overset{n}{\underset{t=1}{\sum}} \eps_t\eps^\top_t -\Psi^\ast_n\|_{\max} \\
& &\qquad + m^2\,\|\Lambda^\ast_n\|^2_{\max}\|\frac{1}{n}\overset{n}{\underset{t=1}{\sum}}F_tF^\top_t - \Eb[F_tF^\top_t]\|_{\max}\Bigg)\|\Sigma^{\ast-1}_n\|^2_s \|\Delta_{n,u}\|_F\\
& \leq & p_n\Big(2\,\underset{1 \leq k\leq p_n}{\max}\|\lambda^\ast_{k}\|_2 \, \sqrt{m} O_p(\sqrt{\frac{\log(p_n)}{n}}) + O_p(\sqrt{\frac{\log(p_n)}{n}}) + m^2\,\|\Lambda^\ast_n\|^2_{\max}O_p(\sqrt{\frac{\log(p_n)}{n}})\Big)\|\Sigma^{\ast-1}_n\|^2_s \|\Delta_{n,u}\|_F,
\end{eqnarray*}
applying Lemma \ref{deviation_ineq}. 
Therefore, we obtain $|K_1| \leq O_p(\alpha^2_n C_{\delta})$. As for $K_2$, note that
\begin{equation*}
\|v \Sigma^{\ast-1}_n\Delta_{n,u}\|_s \leq \|\Sigma^{\ast-1}_n\|_s\|\Delta_{n,u}\|_F \leq O(1) C_\delta \alpha_n=o(1).
\end{equation*}
By the expansion $(A + \Delta)^{-1} = \overset{\infty}{\underset{k=0}{\sum}} (-A^{-1}\Delta)^k A^{-1}$, we write $\Sigma_{n,v} = \Sigma^\ast_n(I_{p_n}+\Sigma^{\ast-1}_nv\Delta_{n,u})$ and thus obtain
\begin{equation*}
\Sigma^{-1}_{n,v} = (I_{p_n}+\Sigma^{\ast-1}_nv\Delta_{n,u})^{-1}\Sigma^{\ast-1}_n = \overset{\infty}{\underset{k=0}{\sum}} (-\Sigma^{\ast-1}_nv\Delta_{n,u})^k \Sigma^{\ast-1}_n \Sigma^{\ast-1}_n = 
\Sigma^{\ast-1}_n - \Sigma^{\ast-1}_n v \Delta_{n,u} \Sigma^{\ast-1}_n + \Rc,
\end{equation*}
with $\Rc = (\Sigma^{\ast-1}_nv\Delta_{n,u})^2(I_{p_n}+\Sigma^{\ast-1}_nv\Delta_{n,u})^{-1} \Sigma^{\ast-1}_n$. The remainder $\Rc$ can be bounded as
\begin{equation*}
\|\Rc\|_s \leq \|\Sigma^{\ast-1}_n\Delta_{n,u}\|^2_s\|(I_{p_n}+\Sigma^{\ast-1}_nv\Delta_{n,u})^{-1}\|_s\|\Sigma^{\ast-1}_n\|_s \leq \frac{\|\Sigma^{\ast-1}_n\Delta_{n,u}\|^2_s \|O(1)}{1-\|\Sigma^{\ast-1}_n\Delta_{n,u}\|_s} =o(1). 
\end{equation*}
Therefore, we deduce $\Sigma^{-1}_{n,v} = 
\Sigma^{\ast-1}_n (I_{p_n}- v \Delta_{n,u} \Sigma^{\ast-1}_n + o(1))$, and so $\|\Sigma^{-1}_{n,v}\|_s = \underline{\mu}+O_p(\alpha_n)$. Moreover, since $\|\widehat{S}_n-\Sigma^\ast_n\|_s=O_p(\alpha_n)$, we get
\begin{equation*}
\widehat{S}_n\Sigma^{-1}_{n,v} = (\widehat{S}_n-\Sigma^\ast_n)\Sigma^{-1}_{n,v} + \Sigma^\ast_n\Sigma^{-1}_v = o_p(1) + I_{p_n} + O_p(\alpha_n) = I_{p_n}+o_p(1).
\end{equation*}
So we obtain $g(v,\Sigma_{n,v}) = \Sigma^{\ast-1}_n\otimes \Sigma^{\ast-1}_n + O_p(\alpha_n)$. This implies
\begin{eqnarray*}
\lefteqn{K_2 = \text{vec}(\Delta_{n,u})^\top \int^1_0 \Big[\Sigma^{\ast-1}_n\otimes \Sigma^{\ast-1}_n + o_p(1)\Big](1-v)\text{d}v \, \text{vec}(\Delta_{n,u})}\\
& \geq & \lambda_{\min}(\Sigma^{\ast-1}_n\otimes \Sigma^{\ast-1}_n )\|\Delta_{n,u}\|^2_F/2 (1+o_p(1))= \underline{\mu}^{-2} C^2_{\delta}\alpha^2_n/2(1+o_p(1)).
\end{eqnarray*}
Finally, we consider the penalty part $K_3$. By Assumption \ref{assumption_regularity_penalty_n}, we have
\begin{equation*}
|K_3| \leq \sqrt{s_n} \alpha_n \underset{1\leq k \leq p_nm}{\max}\{\partial_{1}p(|\theta^\ast_{n\Lambda}|,\gamma_n)\}\|\Delta_{n,u}\|_F + 2 \alpha^2_n \underset{1\leq k \leq p_nm}{\max}\{\partial^2_{11}p(|\theta^\ast_{n\lambda}|,\gamma_n)\} \|\Delta_{n,u}\|^2_F.
\end{equation*}
Therefore, putting the pieces together, 
for large enough $C_\delta>0$,
we deduce
\begin{equation*}
\Lb^{\text{pen}}_n(\theta^\ast_n+\text{vec}(\Delta_{n,u}))-\Lb^{\text{pen}}_n(\theta^\ast_n) = n\,\text{vec}(\Delta_{n,u})^\top \big(\Sigma^{\ast-1}_n\otimes \Sigma^{\ast-1}_n\big)\text{vec}(\Delta_{n,u}) (1+o_p(1)),
\end{equation*}
which is larger than $n\underline{\mu}^{-2} C^2_{\delta}\alpha^2_n/2 > 0$. 
We deduce (\ref{obj_bound}) and finally $\|\widehat\Sigma_n - \Sigma^\ast_n\|_F=O_p(\alpha_n)$.

\mds
\noindent\textbf{\emph{Least squares loss.}} We now establish the result when $\ell_n(X_t;\theta_n) = \text{tr}\big((X_tX^\top_t - \Sigma_n)^2\big)$. We note that $n^{-1}\sum^n_{t=1} \text{tr}\big((X_tX^\top_t - \Sigma_n)^2\big)$,
is equivalent to $\|\widehat{S}_n-\Sigma\|^2_F$ up to some constant terms that do not depend on $\Lambda_n,\Psi_n$. Indeed:
\begin{eqnarray*}
\text{tr}\big((\widehat{S}_n - \Sigma)^2\big)=\text{tr}\big(\widehat{S}^\top_n \widehat{S}_n\big)-\text{tr}\big(\widehat{S}^\top_n\Sigma_n+\Sigma^\top_n\widehat{S}_n\big)+\text{tr}\big(\Sigma^\top_n\Sigma_n\big),
\end{eqnarray*}
which is, up to constant terms that do not depend on $\Sigma$, equivalent to:
\begin{eqnarray*}
\frac{1}{n}\overset{n}{\underset{t=1}{\sum}}\text{tr}\big((X_tX^\top_t - \Sigma_n)^2\big) =  \frac{1}{n}\overset{n}{\underset{t=1}{\sum}}\text{tr}\big(\big(X_tX^\top_t\big)^\top X_tX^\top_t-\big(X_tX^\top_t\big)^\top\Sigma_n-\Sigma^\top_n X_tX^\top_t+\Sigma^\top_n\Sigma_n\big).
\end{eqnarray*}
Hereafter, we will work with the form $\|\widehat{S}_n-\Sigma_n\|^2_F$ and so we write $\Lb_n(\theta_n) = \|\widehat{S}_n-\Sigma_n\|^2_F$. Following the same steps and using the same notations as in the Gaussian case, we have the expansion
\begin{eqnarray*}
\lefteqn{\frac{1}{n}\Lb^{\text{pen}}_n(\theta^\ast_n+\text{vec}(\Delta_{n,u}))-\frac{1}{n}\Lb^{\text{pen}}_n(\theta^\ast_n)}\\
& \geq & \|\widehat{S}_n-(\Sigma^\ast_n+\Delta_{n,u})\|^2_F-\|\widehat{S}_n-\Sigma^\ast_n\|^2_F \\
& & + \underset{k \in \Sc_n}{\sum}\Big(p(|\theta^\ast_{n\Lambda}+\alpha_n\text{vec}(\Delta_{n,u,\Lambda})|,\gamma_n)-p(|\theta^\ast_{n\Lambda}|,\gamma_n)\Big)\\
& = & 2\text{vec}(\Sigma^\ast_n-\widehat{S}_n)^\top \text{vec}(\Delta_{n,u}) + \text{vec}(\Delta_{n,u})^\top 2(I_{p_n} \otimes I_{p_n}) \text{vec}(\Delta_{n,u}) \\
& & +\underset{k \in \Sc_n}{\sum} \Big(\alpha_n \partial_{1}p(|\theta^\ast_{n\lambda}|,\gamma_n)\text{vec}(L_{n,u})_k + \alpha^2_n\partial^2_{11}p(|\theta^\ast_{n\lambda}|,\gamma_n)\text{vec}(L_{n,u})^2_k\{1+o(1)\}\Big)\\
& =: & K_1+K_2+K_3.
\end{eqnarray*}
$K_1$ can be bounded as
\begin{equation*}
|K_1| \leq 2 \biggl|\underset{ij}{\sum}\big[\Sigma^\ast_n-\widehat{S}_n\big]_{ij}\Delta_{n,u}\biggr| \leq \alpha_n \|\Delta_{n,u}\|_F.
\end{equation*}
We thus get $|K_1| \leq O_p(\alpha^2_n C_{\delta})$. Now, since $\nabla^2_{\text{vec}(\Sigma_n)\text{vec}(\Sigma_n)^\top}\Lb_n(\theta^\ast_n) = I_{p_n}\otimes I_{p_n}$, it implies that $K_2 \geq \|\Delta_{n,u}\|^2_F$. We thus deduce $\|\widehat{\Sigma}_n-\Sigma^\ast_n\|_F=O_p(\alpha_n)$.

\subsection{Proof of Theorem \ref{sparsistency}}


\noindent\emph{\textbf{Consistency of $\widehat{\theta}_n$.}} From Theorem \ref{Theorem_existence_consistent}, 
we have 
\[
\|\widehat{\Lambda}_n\widehat{\Lambda}^\top_n+\widehat{\Psi}_n-\Sigma^\ast_n\|_F=O_p(p_n\sqrt{s_n\log(p_n)/n}+\sqrt{p_n}A_n).
\]
Define $\alpha_n = \sqrt{p_n} \big(\sqrt{p_ns_n\log(p_n)/n}+A_n\big)$. Then, under Assumptions \ref{assumption_AR_condition}--\ref{assumption_space_lambda}, by Theorem 1 of \cite{kano1983}, we obtain:
\begin{equation*}
\|\widehat{\theta}_n-\theta^\ast_n\|_{\infty} = O_p(\alpha_n).
\end{equation*}
Therefore, for any $\delta>0$, there exists $C_{\delta}>0$ such that $\Pb(\|\widehat{\theta}_n-\theta^\ast_n\|_2 > C_{\delta} \sqrt{p_n}\alpha_n) <\delta$. The $\|\cdot\|_{\infty}$-consistency of $\widehat{\theta}_n$ follows from the ``strong identifiability'' property of the factor model defined in \cite{kano1983}: the latter property holds if and only if, for any $\eps>0$, $\exists \delta>0$ such that if $\|\Sigma_n-\Sigma^\ast_n\|_{\max} < \delta$, 
then $\|\Lambda_n - \Lambda^\ast_n\|_{\max}<\eps$ and $\|\Psi_n-\Psi^\ast_n\|_{\max}<\eps$, with $\Sigma_n = \Lambda_n\Lambda_n^\top + \Psi_n$.  

The proof of the sparsity property is performed in the same vein as in \cite{lam2009}. We first provide the proof for the estimator deduced from the Gaussian loss. The least squares case will follow. 

\mds

\noindent\emph{\textbf{Sparsisty property.}}\\
\noindent\emph{Gaussian loss.} We consider an estimator $\widehat{\theta}_n = (\widehat{\theta}^\top_{n\Lambda,1},\widehat{\theta}^\top_{n\Lambda,2},\widehat{\theta}^\top_{n\Psi})^\top$ of $\theta^\ast_n$ satisfying $\|\widehat{\theta}_n-\theta^\ast_n\|_2=O_p(\sqrt{p_n}\alpha_n)$. We aim to show
\begin{equation}\label{sparsistency_obj}
\Lb^{\text{pen}}_n((\widehat{\theta}^\top_{n\Lambda,1},0^\top,\widehat{\theta}^\top_{n\Psi})^\top) = \underset{\|\widehat{\theta}_{n\Lambda,2}\|_2\leq C\sqrt{p_n}\alpha_n}{\min}\, \Lb^{\text{pen}}_n((\widehat{\theta}^\top_{n\Lambda,1},\widehat{\theta}^\top_{n\Lambda,2},\widehat{\theta}^\top_{n\Psi})^\top),
\end{equation}
for any constant $C>0$ with probability tending to one. Let $u_n = C \sqrt{p_n}\alpha_n$. To prove (\ref{sparsistency_obj}), it is sufficient to show that for any $\theta_n$ such that $\|\theta_n-\theta^\ast_n\|_2 \leq u_n$, with probability tending to one, we have:
\begin{equation*}
\partial_{\theta_{n\Lambda,j}}\Lb^{\text{pen}}(\theta_n)>0 \; \text{when} \; 0 < \theta_{n\Lambda,j} < u_n, \;\; \partial_{\theta_{n\Lambda,j}}\Lb^{\text{pen}}(\theta_n)<0 \; \text{when} \; -u_n < \theta_{n\Lambda,j} < 0,
\end{equation*}
for $j \in \Sc^c_n$. Using the derivative formulas of Section \ref{derivatives} of the Appendix, for a minimizer $(\Lambda_n,\Psi_n)$, we have:
\begin{equation*}
\partial_{\theta_{n\Lambda,j}}\Lb^{\text{pen}}(\theta_n) = n\,2\,\text{vec}\big(\Sigma^{-1}_n(\Sigma_n-\widehat{S}_n)\Sigma^{-1}_n\Lambda_n\big)_j + n\,\partial_{1}p(|\theta_{n\Lambda,j}|,\gamma_n)\,\text{sgn}(\theta_{n\Lambda,j}),
\end{equation*}
where $\text{sgn}(x)$ is the sign of $x$. We aim to show that the sign of $\partial_{\theta_{n\Lambda,j}}\Lb^{\text{pen}}(\theta_n)$ depends on $\text{sgn}(\theta_{n\Lambda,j})$ only with probability tending to one: this means that the estimated parameter is at $0$, implying $\widehat{\theta}_{n\Lambda,j}=0$ for all $j\in \Sc^c_n$. First, by Lemma 1 of \cite{lam2009}, we have
\begin{eqnarray*}
\lefteqn{\text{vec}\big(\Sigma^{-1}_n(\Sigma_n-\widehat{S}_n)\Sigma^{-1}_n\Lambda_n\big)_j \leq \|\Sigma^{-1}_n(\Sigma_n-\widehat{S}_n)\Sigma^{-1}_n\Lambda_n\|_s}\\
& \leq & \|\Sigma^{-1}_n\|^2_s \|\Sigma_n-\widehat{S}_n\|_s \|\Lambda_n\|_s \leq O(1) \big(\|\Sigma_n-\Sigma^\ast_n\|_s+\|\Sigma^\ast_n-\widehat{S}_n\|_s\big) \big(\|\Lambda_n-\Lambda^\ast_n\|_s+\|\Lambda^\ast_n\|_s\big).
\end{eqnarray*}
Since $\|\Sigma^{-1}_n\|^2_s \leq \|\Psi^{-1}_n\|_s=O(1)$, we deduce
\begin{equation*}
\text{vec}\big(\Sigma^{-1}_n(\Sigma_n-\widehat{S}_n)\Sigma^{-1}_n\Lambda_n\big)_j
\leq O(1)\left\{O_p(\alpha_n) + O_p\left(p_n \sqrt{\frac{\log(p_n)}{n}}\right)\right\} \big(O_p(\sqrt{p_n}\alpha_n)+\sqrt{s_n}\big).
\end{equation*}
Therefore, 
\begin{eqnarray*}
\lefteqn{\partial_{\theta_{n\Lambda,j}}\Lb^{\text{pen}}(\theta_n) = n\,O_p(\sqrt{\frac{p^3_ns_n\log(p_n)}{n}}+p_n^2\sqrt{\frac{\log(p_n)}{n}}A_n) + n\,\partial_{\Lambda_{n,ij}}p(|\theta_{n,\Lambda_{ij}}|,\gamma_n)\,\text{sgn}(\theta_{n\Lambda,j})}\\
& = & n\,\gamma_n \Big\{O_p(\gamma^{-1}_n\sqrt{\frac{p^3_ns_n\log(p_n)}{n}}+\gamma^{-1}_np_n^2\sqrt{\frac{\log(p_n)}{n}}A_n)+\gamma^{-1}_n\partial_{1}p(|\theta_{n\Lambda,j}|,\gamma_n)\,\text{sgn}(\theta_{n\Lambda,j})\Big\}.
\end{eqnarray*}
Under $\gamma^{-1}_n\sqrt{\frac{p^3_ns_n\log(p_n)}{n}} \rightarrow 0$, 
$p_n^2\sqrt{\log(p_n)}A_n = o(\gamma_n\sqrt{n})$, $\underset{n \rightarrow \infty}{\lim \inf} \, \underset{x \rightarrow 0^+}{\lim \inf} \, \gamma^{-1}_n\partial_{1}p(x,\gamma_n) >0$, the sign of $\partial_{\Lambda_{n,ij}}\Lb^{\text{pen}}(\theta_n)$ depends on $\text{sgn}(\theta_{n,\Lambda_{ij}})$ with probability tending to one.  

\mds

\noindent\emph{Least squares loss.} Following the same steps as in the Gaussian loss, for any $\theta_n$ such that $\|\theta_n-\theta^\ast_n\|_2 \leq u_n$ and $u_n = C \sqrt{p_n}\alpha_n$, we show that the sign of $\partial_{\theta_{n\Lambda,j}}\Lb^{\text{pen}}(\theta_n)$ depends on $\text{sgn}(\theta_{n\Lambda,j})$ and $j \in \Sc^c_n$ only. For a minimizer $(\Lambda_n,\Psi_n)$, we have:
\begin{equation*}
\partial_{\theta_{n\Lambda,j}}\Lb^{\text{pen}}(\theta_n) = n\,4 \, \text{vec}\big((\Sigma_n-\widehat{S}_n)\Lambda_n\big)_j + n\,\partial_{1}p(|\theta_{n\Lambda,j}|,\gamma_n)\,\text{sgn}(\theta_{n\Lambda,j}),
\end{equation*}
for any $j \in \Sc^c_n$. The first order term can be bounded as $\text{vec}\big((\Sigma_n-\widehat{S}_n)\Lambda_n\big)_j \leq \|\Sigma_n-\widehat{S}_n\|_s \|\Lambda_n\|_s$. Therefore, we deduce
\begin{eqnarray*}
\partial_{\theta_{n\Lambda,j}}\Lb^{\text{pen}}(\theta_n) 
=  n\,\gamma_n \Bigg\{
O_p\bigg(\gamma^{-1}_n\sqrt{\frac{p^3_ns_n\log(p_n)}{n}}+\gamma^{-1}_np_nA_n\bigg)
+\gamma^{-1}_n\partial_{1}p(|\theta_{n\Lambda,j}|,\gamma_n)\,\text{sgn}(\theta_{n\Lambda,j})\Bigg\}.
\end{eqnarray*}
Hence, $\widehat{\theta}_{n\Lambda,k}=0$ for any $k \in \Sc^c_n$ with probability one.

\mds

\noindent\emph{\textbf{Consistency of the GLS estimator.}} We now move to the consistency of the GLS estimator $\widehat{F}_t$, for any fixed $t=1,\ldots,n$. The proof can be done in the same spirit as in \cite{bailiao2016}, Theorem 3.1. Based on the factor decomposition $X_t = \Lambda^\ast_n F_t + \eps_t$, for any $t$, we have
\begin{equation*}
\widehat{F}_t-F_t = -(\widehat{\Lambda}^\top_n\widehat{\Psi}^{-1}_n \widehat{\Lambda}_n)^{-1} 
\widehat{\Lambda}^\top_n \widehat{\Psi}^{-1}_n(\widehat{\Lambda}_n-\Lambda^\ast_n)F_t + (\widehat{\Lambda}^\top_n\widehat{\Psi}^{-1}_n \widehat{\Lambda}_n)^{-1} 
\widehat{\Lambda}^\top_n \widehat{\Psi}^{-1}_n \eps_t =:M_1+M_2.
\end{equation*}
We first consider $M_1$. Define $\widehat{\xi}^{(1)}_k$ as the $k$-th column of $\widehat{\Lambda}^\top_n\widehat{\Psi}^{-1}_n$. It can be bounded as
\begin{eqnarray*}
\lefteqn{|M_1| \leq \|(\widehat{\Lambda}^\top_n\widehat{\Psi}^{-1}_n \widehat{\Lambda}_n)^{-1}\|_s\sqrt{\sum^{p_n}_{k=1}\|\widehat{\xi}^{(1)}_k(\widehat{\lambda}_{k}-\lambda^\ast_k)\|^2_2}\|F_t\|}\\
& \leq & \lambda_{\max}(\widehat{\Psi}_n) \lambda_{\min}(\widehat{\Lambda}^\top_n\widehat{\Lambda}_n)^{-1}  O_p(1) \|\widehat{\Lambda}_n-\Lambda^\ast_n\|_F\|F_t\| \leq O_p(p^{-\nu}_n) O_p(1)  \|\widehat{\Lambda}_n-\Lambda^\ast_n\|_F \|F_t\| =o_p(1),
\end{eqnarray*}
since $\max_k\|\widehat{\xi}^{(1)}_k\|_2=O_p(1)$
and under $\lambda_{\min}(\widehat{\Lambda}^\top_n\widehat{\Lambda}_n)^{-1} = O_p(p^{-\nu}_n)$, with $1/2<\nu\leq 1$.
As for $M_2$, it can be rewritten as
\begin{equation*}
M_2 = (\widehat{\Lambda}^\top_n\widehat{\Psi}^{-1}_n \widehat{\Lambda}_n)^{-1} 
\Lambda^{\ast\top}_n \Psi^{\ast-1}_n \eps_t + (\widehat{\Lambda}^\top_n\widehat{\Psi}^{-1}_n \widehat{\Lambda}_n)^{-1} 
(\widehat{\Lambda}^\top_n \widehat{\Psi}^{-1}_n-\Lambda^{\ast\top}_n \Psi^{\ast-1}_n)\eps_t =: L_1+L_2.
\end{equation*}
We have
\begin{eqnarray*}
\lefteqn{|L_2|  \leq  O_p(p^{-\nu}_n)\|(\widehat{\Lambda}^\top_n \widehat{\Psi}^{-1}_n-\Lambda^{\ast\top}_n \Psi^{\ast-1}_n)\eps_t\|_F}\\
&  & \leq O_p(p^{-\nu}_n)\Big(\|(\widehat{\Lambda}^\top_n-\Lambda^\ast_n)\widehat{\Psi}^{-1}_n\eps_t\|_F+\|\Lambda^{\ast\top}_n(\widehat{\Psi}^{-1}_n-\Psi^{\ast-1}_n)\eps_t\|_F\Big)=:P_1+P_2.
\end{eqnarray*}
First, we treat $P_1$. Define by $\widehat{v}_{t,k}$ as the $k$-th entry of $\widehat{\Psi}^{-1}_n\eps_t$. Then
\begin{equation*}
P_1 \leq O_p(p^{-\nu}_n)\sqrt{\sum^{p_n}_{k=1}\|(\widehat{\lambda}_{k}-\lambda^\ast_k)\widehat{v}_{t,k}\|^2_2} =O_p(p^{-\nu}_n)O_p(\sqrt{p_n})\|\widehat{\Lambda}^\top_n-\Lambda^\ast_n\|_F = o_p(1).
\end{equation*}
To bound $P_2$, first note $P_2 = O_p(p^{-\nu}_n)\|\Lambda^{\ast\top}_n\widehat{\Psi}^{-1}_n(\Psi^{\ast}_n-\widehat{\Psi}_n)\Psi^{\ast-1}_n\eps_t\|_F$. Then define $\widehat{\xi}^{(2)}_k$ as the $k$-th column of $\Lambda^{\ast\top}_n\widehat{\Psi}^{-1}_n$ and $v_{t,k}$ the $k$-th entry of $\Psi^{\ast-1}_n\eps_t$. We obtain
\begin{equation*}
\|\Lambda^{\ast\top}_n(\widehat{\Psi}^{-1}_n-\Psi^{\ast-1}_n)\eps_t\|_F = \sqrt{\sum^{p_n}_{k=1}\|\widehat{\xi}^{(2)}_k(\Psi^\ast_{n,kk}-\widehat{\Psi}_{n,kk})v_{t,k}\|^2_2} \leq \sqrt{\sum^{p_n}_{k=1}\|\widehat{\xi}^{(2)}_k\|^2|v_{t,k}|^2|\widehat{\Psi}_{n,kk}-\Psi^\ast_{n,kk}|^2}.
\end{equation*}
Since $\max_k\|\widehat{\xi}^{(2)}_k\|_2=O_p(1)$ and $\max_k\|v_{t,k}\|_2 = O_p(\log(p_n))$, we deduce 
\begin{equation*}
\|\Lambda^{\ast\top}_n(\widehat{\Psi}^{-1}_n-\Psi^{\ast-1}_n)\eps_t\|_F \leq O_p(\log(p_n)\|\widehat{\Psi}_n-\Psi^\ast_n\|_F).
\end{equation*}
Therefore, we get 
\begin{equation*}
P_2 = O_p(p^{-\nu}_n) O_p(\log(p_n)\|\widehat{\Psi}_n-\Psi^\ast_n\|_F) = o_p(1).
\end{equation*}
Hence, denoting by $\xi_k$ the $k$-th column of $\Lambda^{\ast\top}_n\Psi^{\ast-1}_n$, we can deduce 
\begin{equation*}
\|\widehat{F}_t-F_t\|\leq 
O_p(p^{-\nu}_n) \sum^{p_n}_{k=1}\|\eps_{t,k}\xi_k\| + o_p(1) = O_p(p^{-\nu+1/2}_n) + o_p(1) = o_p(1).
\end{equation*}

\section{Derivative formulas}\label{derivatives}

In this section, to ease the notations, we do not index by $n$ the model parameters, the function $\ell(\cdot;\cdot)$ and the sample variance-covariance $\widehat{S}$. We provide here the formulas of the first and second derivatives of $\Lb_n(\cdot)$ with respect to the factor model parameters $\theta=(\theta^\top_{\Lambda},\theta^\top_{\Psi})$. Although the proofs of Theorem \ref{Theorem_existence_consistent} and Theorem \ref{sparsistency} do not require the second and cross-derivatives, it is still of interest to derive these quantities. We will present the derivative formulas when $\Lb_n(\theta)$ is under the form $\Lb_n(\theta)=n^{-1}\sum^n_{t=1}\ell(X_t;\theta)$. Hereafter, $\theta = (\theta^\top_{\Lambda},\theta^\top_{\Psi})^\top$, where $\theta_{\Lambda} = \text{vec}(\Lambda)$ and $\Sigma = \Lambda\Lambda^\top+\Psi$, $\Lambda \in \Rb^{p \times m}, \Psi \in \Rb^{p \times p}$. As for the $\Psi$ parameter, we will provide the first derivative with respect $\theta_\Psi = \text{vec}(\Psi) \in \Rb^{p^2}$, i.e., for any $\Rb^{p \times p}$ matrix, and with respect to $\theta_{\Psi} = \text{diag}(\Psi)=(\sigma^2,\ldots,\sigma^2_p)^\top \in \Rb^p$, i.e., when accounting for the diagonal restriction on $\Psi$. The second order derivative will be provided with respect to $\theta_{\Psi}=\text{vec}(\Psi) \in \Rb^{p^2}$ and $\theta_{\Psi}=(\sigma^2,\ldots,\sigma^2_p)^\top \in \Rb^p$. 

\subsection{Gaussian loss function}\label{Gaussian_loss}

The Gaussian loss $\Lb_n(\theta) = n^{-1}\sum^n_{t=1}\ell(X_t;\theta)$, with $\ell(X_t;\theta) = \text{tr}\big((X_t X^\top_t) \Sigma^{-1}\big)+\log(|\Sigma|)$, can be written as $\Lb_n(\theta) = \text{tr}\big(\widehat{S}\Sigma^{-1}\big)+\log(|\Sigma|)$, $\widehat{S}=n^{-1}\sum^n_{t=1}X_tX^\top_t$, which will be used hereafter. The first differential of $\Lb_n(\cdot)$ is:
\begin{equation*}
\dd \Lb_n(\theta)= \text{tr}\big(-\widehat{S}\Sigma^{-1}(\dd\Sigma)\Sigma^{-1} + \Sigma^{-1}(\dd\Sigma)\big) =: \text{tr}\big(V (\dd\Sigma)\big), \, V= \Sigma^{-1}- \Sigma^{-1}\widehat{S}\Sigma^{-1}.
\end{equation*}
The differential of $V$ is:
{\small{\begin{eqnarray*}
\dd V = -\Sigma^{-1}(\dd\Sigma)\Sigma^{-1} + \Sigma^{-1}(\dd\Sigma)\Sigma^{-1}\widehat{S}\Sigma^{-1} + \Sigma^{-1}\widehat{S}\Sigma^{-1}(\dd\Sigma)\Sigma^{-1} = \Sigma^{-1}(\dd\Sigma)\Sigma^{-1} - \Sigma^{-1}(\dd\Sigma)V - V(\dd\Sigma)\Sigma^{-1}.
\end{eqnarray*}}}
The second order differential of $\Lb_n(\cdot)$ becomes:
\begin{eqnarray*}
\dd^2 \Lb_n(\theta)=\text{tr}\big((\dd V)(\dd\Sigma)+V(\dd^2\Sigma)\big)
= \text{tr}\big(\Sigma^{-1}(\dd\Sigma)\Sigma^{-1}(\dd\Sigma)\big)-2\,\text{tr}\big(\Sigma^{-1}(\dd\Sigma)V(\dd\Sigma)\big)+\text{tr}\big(V(\dd^2\Sigma)\big).
\end{eqnarray*}
Here, $\dd\Sigma = (\dd\Lambda)\Lambda^\top+\Lambda(\dd\Lambda)^\top+\dd\Psi, \;\; \dd^2\Sigma=2(\dd\Lambda)(\dd\Lambda)^\top$, since $\Lambda$ and $\Psi$ are both linear in the parameters, so that $\dd^2\Lambda$ and $\dd^2\Psi$ are both zero. We deduce that the first differential of $\Lb_n(\cdot)$ is $\dd \Lb_n(\theta) = 2\, \text{tr}\big(\Lambda^\top V(\dd\Lambda)\big) + \text{tr}\big(V(\dd\Psi)\big)$, and the second order differential is
\begin{eqnarray*}
\lefteqn{\dd^2\Lb_n(\theta)=\text{tr}\big(\Sigma^{-1}(\dd\Sigma)\Sigma^{-1} (\dd\Sigma)\big) - 2\,\text{tr}\big(\Sigma^{-1}(\dd\Sigma)V (\dd\Sigma)\big) + \text{tr}\big(V (\dd^2\Sigma)\big)}\\
& = & \text{tr}\big(\Sigma^{-1}(\dd\Sigma)\Sigma^{-1} (\dd\Sigma)\big) - 2\,\text{tr}\big(\Sigma^{-1}(\dd\Sigma)V (\dd\Sigma)\big) + 2\,\text{tr}\big(V (\dd \Lambda)(\dd\Lambda)^\top\big).
\end{eqnarray*}
Let us explicit the expression $(\dd\Sigma)\Sigma^{-1}(\dd\Sigma)$:
\begin{eqnarray*}
\lefteqn{(\dd\Sigma)\Sigma^{-1}(\dd\Sigma)=\left[(\dd\Lambda)\Lambda^\top+\Lambda(\dd\Lambda)^\top+\dd\Psi\right]\Sigma^{-1}\left[(\dd\Lambda)\Lambda^\top+\Lambda(\dd\Lambda)^\top+\dd\Psi\right]}\\
& = & (\dd\Lambda)\Lambda^\top\Sigma^{-1}(\dd\Lambda)\Lambda^\top + (\dd\Lambda)\Lambda^\top\Sigma^{-1}\Lambda(\dd\Lambda)^\top+(\dd\Lambda)\Lambda^\top\Sigma^{-1}(\dd\Psi) + \Lambda (\dd\Lambda)^\top \Sigma^{-1}(\dd\Lambda)\Lambda^\top\\
& & +\Lambda(\dd\Lambda)^\top\Sigma^{-1}\Lambda(\dd\Lambda)^\top +\Lambda(\dd\Lambda)^\top\Sigma^{-1}(\dd\Psi) + (\dd\Psi)\Sigma^{-1}(\dd\Lambda)\Lambda^\top+(\dd\Psi)\Sigma^{-1}\Lambda(\dd\Lambda)^\top+(\dd\Psi)\Sigma^{-1}(\dd\Psi).
\end{eqnarray*}
By the properties of the trace operator:
\begin{eqnarray*}
\lefteqn{\text{tr}\big((\dd\Sigma)\Sigma^{-1}(\dd\Sigma)V\big)=2\,\text{tr}\big(\Lambda^\top V (\dd\Lambda)\Lambda^\top\Sigma^{-1}(\dd\Lambda)\big)+\text{tr}\big(\Lambda^\top\Sigma^{-1}\Lambda(\dd\Lambda)^\top V(\dd\Lambda)\big)}\\
& & + \text{tr}\big(\Lambda^\top V \Lambda (\dd\Lambda)^\top\Sigma^{-1}(\dd\Lambda)\big) + 2\, \text{tr}\big(V(\dd\Lambda)\Lambda^\top \Sigma^{-1}(\dd\Psi)\big)+2\,\text{tr}\big(\Sigma^{-1}(\dd\Lambda)\Lambda^\top V (\dd\Psi)\big)\\
&&+\text{tr}\big((\dd\Psi)\Sigma^{-1}(\dd\Psi)V\big),
\end{eqnarray*}
and 
\begin{eqnarray*}
\lefteqn{\text{tr}\big((\dd\Sigma)\Sigma^{-1}(\dd\Sigma)\Sigma^{-1}\big)=2\,\text{tr}\big(\Lambda^\top\Sigma^{-1}(\dd\Lambda)\Lambda^\top\Sigma^{-1}(\dd\Lambda)\big)}\\
& & +2\,\text{tr}\big(\Lambda^\top\Sigma^{-1}\Lambda(\dd\Lambda)^\top\Sigma^{-1}(\dd\Lambda)\big)+4\,\text{tr}\big(\Sigma^{-1}(\dd\Lambda)\Lambda^\top\Sigma^{-1}(\dd\Psi)\big)+\text{tr}\big((\dd\Psi)\Sigma^{-1}(\dd\Psi)\Sigma^{-1}\big).
\end{eqnarray*}
We thus deduce:
{\small{\begin{eqnarray*}
\lefteqn{\dd^2\Lb_n(\theta)=2\,\text{tr}\big(\Lambda^\top\Sigma^{-1}(\dd\Lambda)\Lambda^\top\Sigma^{-1}(\dd\Lambda)\big)-4\,\text{tr}\big(\Lambda^\top V(\dd\Lambda)\Lambda^\top\Sigma^{-1}(\dd\Lambda)\big)}\\
& & -2\,\text{tr}\big(\Lambda^\top\Sigma^{-1}\Lambda(\dd\Lambda)^\top V (\dd\Lambda)\big)+2\,\text{tr}\big((\dd\Lambda)^\top V (\dd\Lambda)\big)+2\,\text{tr}\big(\Lambda^\top\Sigma^{-1}\Lambda(\dd\Lambda)^\top\Sigma^{-1}(\dd\Lambda)\big)\\
&& -2\,\text{tr}\big(\Lambda^\top V \Lambda (\dd\Lambda)^\top \Sigma^{-1}(\dd\Lambda)\big)
+ 4\,\text{tr}\big(\Sigma^{-1}\Lambda(\dd\Lambda)^\top\Sigma^{-1}(\dd\Psi)\big)-4\,\text{tr}\big(\Sigma^{-1}\Lambda(\dd\Lambda)^\top V (\dd\Psi)\big)\\
&&-4\,\text{tr}\big(V \Lambda (\dd\Lambda)^\top\Sigma^{-1}(\dd\Psi)\big)+ \text{tr}\big(\Sigma^{-1}(\dd\Psi)\Sigma^{-1}(\dd\Psi)\big)-2\,\text{tr}\big(\Sigma^{-1}(\dd\Psi)V (\dd\Psi)\big)\\
& = & \begin{pmatrix}
\dd\text{vec}(\Lambda)\\
\dd\text{vec}(\Psi)
\end{pmatrix}^\top\begin{pmatrix}
H_{\Lambda\Lambda}  & H_{\Lambda\Psi} \\
H_{\Psi\Lambda} & H_{\Psi\Psi}
\end{pmatrix}\begin{pmatrix}
\dd\text{vec}(\Lambda)\\
\dd\text{vec}(\Psi)
\end{pmatrix}.
\end{eqnarray*}}}
Note that the following relationships hold, for any matrices $A,B$ with suitable dimensions, by the formulas for the identification of the Hessian in Theorem 18.6 of \cite{magnus2019}, we have:
\begin{equation*}
\text{tr}\big(A (\dd \Lambda)^\top B(\dd\Lambda)\big) = \big(\dd\text{vec}(\Lambda)\big)^\top\Big[\frac{1}{2}\big(A^\top\otimes B\big)+\big(A\otimes B^\top\big)\Big]\dd\text{vec}(\Lambda),
\end{equation*}
\begin{equation*}
\text{tr}\big(A (\dd\Lambda)B(\dd\Lambda)\big) =\big(\dd\text{vec}(\Lambda)\big)^\top\Big[\frac{1}{2}K_{mp}\big( \big(A^\top \otimes B\big) + \big(B^\top \otimes A\big)\big)\Big]\dd\text{vec}(\Lambda),
\end{equation*}
\begin{equation*}
\text{tr}\big(A(\dd\Lambda)^\top B(\dd\Psi)\big) = \big(\dd\text{vec}(\Psi)\big)^\top  \Big[\frac{1}{2}\big(A^\top\otimes B\big)+\big(A\otimes B^\top\big)\Big] \dd\text{vec}(\Lambda),
\end{equation*}
\begin{equation*}
\text{tr}\big(A(\dd\Psi)^\top B(\dd\Psi)\big)=\big(\dd\text{vec}(\Psi)\big)^\top\Big[\frac{1}{2}\big(A^\top\otimes B\big)+\big(A\otimes B^\top\big)\Big]\dd\text{vec}(\Psi).
\end{equation*}
Then we can write $\dd^2 \Lb_n(\theta)=\sum^4_{k=1} P_k$, with
\begin{eqnarray*}
\lefteqn{P_1 =-2\,\text{tr}\big(\Lambda^\top\Sigma^{-1}\Lambda(\dd\Lambda)^\top V (\dd\Lambda)\big)+2\,\text{tr}\big((\dd\Lambda)^\top V (\dd\Lambda)\big)}\\
&&+2\,\text{tr}\big(\Lambda^\top\Sigma^{-1}\Lambda(\dd\Lambda)^\top\Sigma^{-1}(\dd\Lambda)\big)-2\,\text{tr}\big(\Lambda^\top V \Lambda (\dd\Lambda)^\top \Sigma^{-1}(\dd\Lambda)\big)\\
& = & 2\big(\dd\text{vec}(\Lambda)\big)^\top \Big[ -\big(\Lambda^\top\Sigma^{-1}\otimes V\big)+\big(I_m\otimes V\big)+\big(\Lambda^\top\Sigma^{-1}\Lambda \otimes \Sigma^{-1}\big)-\big(\Lambda^\top V\Lambda\otimes\Sigma^{-1}\big)\Big]\dd\text{vec}(\Lambda),
\end{eqnarray*}
\begin{eqnarray*}
\lefteqn{P_2 = 2\,\text{tr}\big(\Lambda^\top\Sigma^{-1}(\dd\Lambda)\Lambda^\top\Sigma^{-1}(\dd\Lambda)\big)-4\,\text{tr}\big(\Lambda^\top V(\dd\Lambda)\Lambda^\top\Sigma^{-1}(\dd\Lambda)\big)}\\
& = & 2\big(\dd\text{vec}(\Lambda)\big)^\top K_{mp}\Big[\big(\Sigma^{-1}\Lambda\otimes \Lambda^\top\Sigma^{-1}\big)-\big(V\Lambda\otimes\Lambda^\top\Sigma^{-1}\big)-\big(\Sigma^{-1}\Lambda\otimes\Lambda^\top V\big)\Big]\dd\text{vec}(\Lambda),
\end{eqnarray*}
\begin{eqnarray*}
\lefteqn{P_3 =4\,\text{tr}\big(\Sigma^{-1}\Lambda(\dd\Lambda)^\top\Sigma^{-1}(\dd\Psi)\big)-4\,\text{tr}\big(\Sigma^{-1}\Lambda(\dd\Lambda)^\top V (\dd\Psi)\big)-4\,\text{tr}\big(V \Lambda (\dd\Lambda)^\top\Sigma^{-1}(\dd\Psi)\big)}\\
& = & 4\big(\dd\text{vec}(\Psi)\big)^\top \Big[\big(\Sigma^{-1}\Lambda\otimes\Sigma^{-1}\big)-\big(\Sigma^{-1}\Lambda\otimes V\big)-\big(V\Lambda\otimes\Sigma^{-1}\big)\Big]\dd\text{vec}(\Lambda),
\end{eqnarray*}
and 
{\small{\begin{equation*}
P_4 =\text{tr}\big(\Sigma^{-1}(\dd\Psi)\Sigma^{-1}(\dd\Psi)\big)-2\,\text{tr}\big(\Sigma^{-1}(\dd\Psi)V (\dd\Psi)\big)= 4\big(\dd\text{vec}(\Psi)\big)^\top \Big[\big(\Sigma^{-1}\otimes\Sigma^{-1}\big)-2\big(\Sigma^{-1}\otimes V\big)\Big]\dd\text{vec}(\Psi).
\end{equation*}}}
Therefore, $\nabla_{\Lambda}\Lb_n(\theta)= 2\Sigma^{-1}(\Sigma-\widehat{S})\Sigma^{-1}\Lambda, \, \nabla_{\Psi}\Lb_n(\theta)=  \Sigma^{-1}(\Sigma-\widehat{S})\Sigma^{-1}$,
so that in vector form:
\begin{equation*}
\nabla_{\theta_{\Lambda}}\Lb_n(\theta)= 2\,\text{vec}\big(\Sigma^{-1}(\Sigma-\widehat{S})\Sigma^{-1}\Lambda\big)\in\Rb^{pm}, \, \nabla_{\theta_{\Psi}}\Lb_n(\theta)=  \text{vec}\big(\Sigma^{-1}(\Sigma-\widehat{S})\Sigma^{-1}\big) \in \Rb^{p^2}.
\end{equation*}
As for the Hessian matrix, we have:
\begin{eqnarray*}
\lefteqn{\nabla^2_{\theta_{\Lambda}\theta^\top_{\Lambda}}\Lb_n(\theta)= 2\,K_{mp}\left(\Sigma^{-1}\Lambda \otimes \Lambda^\top \Sigma^{-1}\right)-2\,K_{mp}\left[\left(V\Lambda\otimes\Lambda^\top\Sigma^{-1}\right)+\left(\Sigma^{-1}\Lambda\otimes\Lambda^\top V\right)\right]}\\
& & +2\,\left[-\left(\Lambda^\top\Sigma^{-1}\Lambda\otimes V \right)+ \left(I_m \otimes V\right) +  \left(\Lambda^\top\Sigma^{-1}\Lambda\otimes \Sigma^{-1} \right) - \left(\Lambda^\top V \Lambda\otimes \Sigma^{-1} \right)\right],
\end{eqnarray*}
and $\nabla^2_{\theta_{\Psi}\theta^\top_{\Psi}}\Lb_n(\theta) =\left(\Sigma^{-1}\otimes \Sigma^{-1}\right)-2\,\left(\Sigma^{-1}\otimes V\right) \in \Rb^{p^2\times p^2}$. Finally, the cross-derivative is
\begin{equation*}
\nabla^2_{\theta_{\Psi}\theta^\top_{\Lambda}}\Lb_n(\theta) = 2\,\left[\left(\Sigma^{-1}\Lambda \otimes \Sigma^{-1} \right) -\left(\Sigma^{-1}\Lambda \otimes V\right) - \left(V \Lambda \otimes \Sigma^{-1}\right)\right], \; \big(\nabla^2_{\theta_{\Psi}\theta^\top_{\Lambda}} \Lb_n(\theta)\big)^\top = \nabla^2_{\theta_{\Lambda}\theta^\top_{\Psi}} \Lb_n(\theta).
\end{equation*}
Under the diagonality constraint w.r.t. $\Psi$, i.e., $\Psi$ is diagonal with $\theta_{\Psi}=(\sigma^2_1,\ldots,\sigma^2_p)^\top\in\Rb^p$, the gradient becomes $\nabla_{\theta_{\Psi}}\Lb_n(\theta)=\text{diag}\big(I_p\odot\big(\Sigma^{-1}(\Sigma-\widehat{S})\Sigma^{-1}\big)\big)$. To take the diagonality constraint w.r.t. $\Psi$ in the Hessian matrix, we introduce the selection matrix $D^\ast_p\in\Rb^{p^2 \times p}$, defined by $D^\ast_p=\sum^p_{k=1}\left(e_k \otimes e_k\right)e^\top_k$, with $e_k\in\Rb^p$ the $p$-dimensional vector with zero entries only, excluding the $k$-th element which is $1$. Therefore, the Hessian matrix becomes:
$$\nabla^2_{\theta_{\Psi}\theta^\top_{\Psi}}\Lb_n(\theta) =D^{\ast\top}_p\left[\left(\Sigma^{-1}\otimes \Sigma^{-1}\right)-2\,\left(\Sigma^{-1}\otimes V\right)\right]D^\ast_p \in \Rb^{p \times p},$$
and $\nabla^2_{\theta_{\Psi}\theta^\top_{\Lambda}}\Lb_n(\theta) = 2\,D^{\ast\top}_p\left[\left(\Sigma^{-1}\Lambda \otimes \Sigma^{-1} \right) -\left(\Sigma^{-1}\Lambda \otimes V\right) - \left(V \Lambda \otimes \Sigma^{-1}\right)\right]$.

\subsection{Least squares loss function}\label{LS_loss}

Following the preliminary remark in the proof of Theorem \ref{Theorem_existence_consistent}, least squares case, 
$\Lb_n(\theta) = n^{-1}\sum^n_{t=1}\text{tr}\big((X_tX^\top_t - \Sigma)^2\big)$
is equivalent to $\|\widehat{S}-\Sigma\|^2_F$, up to some constant terms that do not depend on $\Lambda,\Psi$. 
Hereafter, to obtain our derivatives, we work with $\Lb_n(\theta) = \|\widehat{S}-\Sigma\|^2_F$. The first differential of $\Lb_n(\cdot)$ is:
\begin{equation*}
\dd \Lb_n(\cdot) = \text{tr}\big(2\Sigma(\dd\Sigma)-2\widehat{S}(\dd\Sigma)\big)=:\text{tr}\big(V(\dd\Sigma)\big), \; V = 2(\Sigma-\widehat{S}).
\end{equation*}
Since $\dd V = 2\dd \Sigma$, the second order differential becomes:
\begin{equation*}
\dd^2\Lb_n(\theta) = 2\,\text{tr}\big((\dd\Sigma)(\dd\Sigma)\big)+\text{tr}\big(V(\dd^2\Sigma)\big), \, \dd\Sigma = (\dd\Lambda)\Lambda^\top+\Lambda(\dd\Lambda)^\top+\dd\Psi, \;\; \dd^2\Sigma=2(\dd\Lambda)(\dd\Lambda)^\top.
\end{equation*}
We have:
\begin{eqnarray*}
\lefteqn{(\dd\Sigma)(\dd\Sigma)=\big((\dd\Lambda)\Lambda^\top+\Lambda(\dd\Lambda)^\top+\dd\Psi\big)\big((\dd\Lambda)\Lambda^\top+\Lambda(\dd\Lambda)^\top+\dd\Psi\big)}\\
& = & (\dd\Lambda)\Lambda^\top(\dd\Lambda)\Lambda^\top+(\dd\Lambda)\Lambda^\top\Lambda(\dd\Lambda)+(\dd\Lambda)\Lambda^\top(\dd\Psi)+\Lambda(\dd\Lambda)^\top(\dd\Lambda)\Lambda^\top+\Lambda(\dd\Lambda)^\top\Lambda(\dd\Lambda)^\top\\
&&+\Lambda(\dd\Lambda)^\top(\dd\Psi)+ (\dd\Psi)(\dd\Lambda)\Lambda^\top+(\dd\Psi)\Lambda(\dd\Lambda)^\top+(\dd\Psi)(\dd\Psi).
\end{eqnarray*}
We then deduce:
\begin{eqnarray*}
\dd \Lb_n(\theta)=\text{tr}\big(V\big((\dd\Lambda)\Lambda^\top+\Lambda(\dd\Lambda)^\top+\dd\Psi\big)\big)=2\,\text{tr}\big(\Lambda^\top V (\dd\Lambda)\big)+\text{tr}\big(V (\dd\Psi)\big),
\end{eqnarray*}
and
{\small{\begin{eqnarray*}
\lefteqn{\dd^2 \Lb_n(\theta)=2\,\text{tr}\big(\big((\dd\Lambda)\Lambda^\top+\Lambda(\dd\Lambda)^\top+\dd\Psi\big)\big((\dd\Lambda)\Lambda^\top+\Lambda(\dd\Lambda)^\top+\dd\Psi\big)\big)+2\,\text{tr}\big(V(\dd\Lambda)(\dd\Lambda)^\top\big)}\\
& = & 4\,\text{tr}\big(\Lambda^\top(\dd\Lambda)\Lambda^\top(\dd\Lambda)\big)+4\,\text{tr}\big(\Lambda^\top\Lambda(\dd\Lambda)^\top(\dd\Lambda)\big)+2\,\text{tr}\big((\dd\Lambda)^\top V(\dd\Lambda)\big)+8\,\text{tr}\big(\Lambda(\dd\Lambda)^\top(\dd\Psi)\big)+2\,\text{tr}\big((\dd\Psi)(\dd\Psi)\big)\\
& = & \begin{pmatrix}
\dd\text{vec}(\Lambda)\\
\dd\text{vec}(\Psi)
\end{pmatrix}^\top\begin{pmatrix}
H_{\Lambda\Lambda}  & H_{\Lambda\Psi} \\
H_{\Psi\Lambda} & H_{\Psi\Psi}
\end{pmatrix}\begin{pmatrix}
\dd\text{vec}(\Lambda)\\
\dd\text{vec}(\Psi)
\end{pmatrix}
\end{eqnarray*}}}
Following the same steps as in the Gaussian loss function, the gradients are $\nabla_{\Lambda}\Lb_n(\theta)=4\,\big(\Sigma-\widehat{S}\big)\Lambda,\, \nabla_{\Psi}\Lb_n(\theta)=2\,\big(\Sigma-\widehat{S}\big)$, i.e., $\nabla_{\theta_{\Lambda}}\Lb_n(\theta)=4\,\text{vec}\big(\big(\Sigma-\widehat{S}\big)\Lambda\big),\, \nabla_{\theta_{\Psi}}\Lb_n(\theta)=2\,\text{vec}\big(\Sigma-\widehat{S}\big)$. The Hessian matrix is:
\begin{eqnarray*}
\nabla^2_{\theta_{\Lambda}\theta^\top_{\Lambda}} \Lb_n(\theta)  = 4\,K_{m,p}
\big(\Lambda\otimes\Lambda^\top\big)+4\,\big(\Lambda^\top\Lambda\otimes I_p\big)+2\,\big(I_m\otimes V\big),\,\nabla^2_{\theta_{\Psi}\theta^\top_{\Psi}} \Lb_n(\theta)  = 2\big(I_p \otimes I_p\big),
\end{eqnarray*}
and the cross-derivative is given by $\nabla^2_{\theta_{\Psi}\theta^\top_{\Lambda}} \Lb_n(\theta) = 4\,\big(\Lambda \otimes I_p\big)$, $\big(\nabla^2_{\theta_{\Psi}\theta^\top_{\Lambda}} \Lb_n(\theta)\big)^\top = \nabla^2_{\theta_{\Lambda}\theta^\top_{\Psi}} \Lb_n(\theta)$. Under the diagonality constraint w.r.t. $\Psi$, i.e., $\Psi$ is diagonal with $\theta_{\Psi}=(\sigma^2_1,\ldots,\sigma^2_p)^\top\in\Rb^p$, the gradient becomes $\nabla_{\theta_{\Psi}}\Lb_n(\theta)=2\,\text{diag}\big(I_p\odot \big(\Sigma-\widehat{S}\big)\big)$ and the Hessian becomes $\nabla^2_{\theta_{\Psi}\theta^\top_{\Psi}} \Lb_n(\theta)  = 2D^{\ast\top}_p\big(I_p \otimes I_p\big)D^\ast_p$ and $\nabla^2_{\theta_{\Psi}\theta^\top_{\Lambda}} \Lb_n(\theta) = 4\,D^{\ast\top}_p\big(\Lambda \otimes I_p\big)$.

\section{Implementations}\label{implementation}

In this section, we discuss the computational issues to estimate the penalized factor model. To ease the notations, the factor model parameters and the variables relating to the dimension are not indexed by $n$. The code was implemented on Matlab and run on Intel(R) Xeon(R) Gold 6242R CPU, 3.10GHz, with Installed RAM 128 GB. We make the code for the sparse factor model and corresponding simulated experiments publicly available in the following Github repository: \url{https://github.com/Benjamin-Poignard/sparse-factor-models} 

\subsection{Algorithm}

Recall that $\theta = (\theta^\top_{\Lambda},\theta^\top_{\Psi})^\top$, where $\theta_{\Lambda} = \text{vec}(\Lambda)$ and $\Sigma = \Lambda\Lambda^\top+\Psi$, $\Lambda \in \Rb^{p \times m}, \Psi \in \Rb^{p \times p}$. The penalized problem is solved according to the following steps:
\begin{itemize}
    \item[\textbf{Step 1}:] Initialization.
    \item[(i)] Solve $\widehat{\theta} := (\widehat{\theta}^\top_{\Lambda},\widehat{\theta}^\top_{\Psi})^\top = \arg\;\min_{\theta\in \Theta_{\text{ic}}}\,\Lb_n(\theta)$,
    with $\Theta_{\text{ic}}$ the parameter set $\Theta_{\text{ic}} = \big\{\theta = (\theta^\top_{\Lambda},\theta^\top_{\Psi}): \; \Lambda^\top \Psi^{-1}\Lambda/p \; \text{diagonal}, \; \Psi \; \text{diagonal}\big\}$. This estimator corresponds to the strict factor model estimator of \cite{bai2012} and is obtained by the EM algorithm following their Section 8.
    \item[(ii)] For a given $\gamma_n$, set a grid size $K$ . $\forall 1 \leq l \leq K$, solve
    \begin{equation*}
    Q^l = \underset{Q: \; Q^\top \, Q = I_m}{\arg\;\min}\;\big\{\sum_{k=1}^{p\,m}p(|\theta_{Q,k}|,\gamma_n)\big\}, \; \text{where} \; \theta_{Q} = \text{vec}(\widehat{\Lambda}\,Q),
    \end{equation*}
    with $\widehat{\Lambda}$ the estimator of step (i). The problem is solved using the algorithm developed by \cite{wen2013} for minimization problems subject to orthogonal constraints. The optimal value $\widetilde{Q}_{\gamma_n}$ is set as the argument $Q^l,1 \leq l \leq K$, providing the minimum value of the previous loss.
    \end{itemize}
    \begin{itemize}
    \item[\textbf{Step 2}:] Iteration. Set $\theta_{\Psi_{(0)}}=\text{diag}(\widehat{\Psi})$, $\widehat{\Psi}$ obtained in \textbf{Step 1}. For each given $\gamma_n$, repeat (i)-(ii) until convergence:
    \begin{itemize}
        \item[(i)] Solve
        $\Lambda_{(k+1)} = \arg\;\min_{\theta_\Lambda} \big\{\Lb_n((\theta^\top_{\Lambda},\theta^\top_{\Psi_{(k)}})^\top)+\sum_{k=1}^{p\,m}p(|\theta_{\Lambda,k}|,\gamma_n)\big\}$.
        The penalized problem is optimized by a gradient descent algorithm based on the updating formulas of Section 4.2 in~\cite{loh2015}. The initial value for $\Lambda$ is set as the rotated matrix $\widehat{\Lambda}\widetilde{Q}_{\gamma_n}$, with $\widehat{\Lambda}$ obtained in \textbf{Step 1}.
        \item[(ii)] Solve $\theta_{\Psi_{(k+1)}} = \arg\;\min_{\theta_\Psi} \Lb_n((\theta^\top_{\Lambda_{(k+1)}},\theta^\top_{\Psi})^\top)$ s.t. $\Psi$ diagonal. The problem is optimized by the EM algorithm.
    \end{itemize}    
\end{itemize}
\textbf{Step 1} is crucial to provide suitable starting values for $\Lambda$. In our experiments, we set $K=200$. We set $a_{\text{mcp}}=3.7$ for the SCAD and  $b_{\text{mcp}}=3.5$ for the MCP.

\subsection{Selection of \texorpdfstring{$\gamma_n$}{Lg}}\label{CV}

The tuning parameter $\gamma_n$ controls the model complexity and must be calibrated for each penalty. The selection strategy depends on the nature of the data: i.i.d. or dependent. In the case of i.i.d. data, we employ a $5$-fold cross-validation procedure, in the same spirit as in, e.g., Section 7.10 of \cite{Hastie2009}. To be specific, we divide the data into $5$ disjoint subgroups of roughly the same size - the so-called folds - and denote the index of observations in the $k$-th fold by $T_k$, $1 \leq k \leq 5$ and the size of the $k$-th fold by $n_k$. The $5$-fold cross-validation score is defined as $\text{CV}(\gamma_n) = \sum^5_{k=1}\Big\{n^{-1}_k\sum_{i \in T_k}\ell(X_t;\widehat{\theta}_{-k}(\gamma_n))\Big\}$,
with $n^{-1}_k\sum_{i \in T_k}\ell(X_t;\widehat{\theta}_{-k}(\gamma_n))$ the non-penalized loss evaluated over the $k$-th fold $T_k$ of size $n_k$, which serves as the test set, and $\widehat{\theta}_{-k}(\gamma_n)$ is the estimator of the factor model based on the sample $\big(\bigcup^5_{k=1}T_k\big)\setminus T_k$ - the training set -. The optimal tuning parameter $\gamma^\ast_n$ is then selected according to: $\gamma^\ast_n = \arg \,\min_{\gamma_n}\, \text{CV}(\gamma_n)$. Then, $\gamma^\ast_n$ is used to obtain the final estimate over the whole data $\big(\bigcup^5_{k=1}T_k\big)$. Here, the minimization of the cross-validation score is performed over a $\{c\sqrt{\log(\text{dim})/n}\}$, with $c$ a user-specified grid of values, set, e.g., as $c=0.001,0.005, 0.01,0.015,\ldots,4$, and $\text{dim}$ the problem dimension, i.e., $\text{dim}=pm$. The cross-validation score also be performed over $c$ simply. 

\mds

Several procedures have been proposed to take into account the dependency among observations, but there is no consensual approach. Based on thorough simulated and real data experiments, \cite{cerqueira2020} concluded that good estimation performances are obtained by the "blocked" cross-validation when the process is stationary; under non-stationary time series, the most accurate estimations are produced by "out-of-sample" methods, where the last part of the data is used for testing. In our experiments, we follow the out-of-sample approach for cross-validation. To be specific, we split the full sample into training and test sets. The training sample corresponds to the first $75\%$ of the entire sample and the test sample to the last $25\%$. For given values of $\gamma_n$, we estimate the sparse factor model in the training set and then compute the loss function in the test set using the estimator computed on the training set. This procedure is repeated for all the $\gamma_n$ candidates and we select the one minimizing the loss. Then, we estimate the model over the full sample using the optimal $\gamma_n$. Formally, the procedure is as follows: divide the data into two sets, the training set and test set; define the cross-validation score as $\text{CV}(\gamma_n) = n^{-1}_{\text{oos}}\sum_{i \in T_{\text{oos}}}\ell(X_t;\widehat{\theta}_{\text{in}}(\gamma_n))$, where $n_{\text{oos}}$ is the size of the test set, $T_{\text{oos}}$ is the index of observations in the test set, and $\widehat{\theta}_{\text{in}}(\gamma_n)$ is the factor model estimator based on the training set with regularization parameter $\gamma_n$. The optimal regularization parameter $\gamma^\ast_n$ is then selected according to: $\gamma^\ast_n = \arg \,\min_{\gamma_n}\, \text{CV}(\gamma_n)$. Then, $\gamma^\ast_n$ is used to obtain the final estimator over the whole data. Here, the minimization of the cross-validation score is performed over $\{c\sqrt{\log(\text{dim})/n}\}$, with $c$ a grid of values as in the i.i.d. case.

\section{Competing models}

\subsection{DCC model}\label{dcc_model}

Rather than a direct specification of the variance-covariance matrix process $(\Sigma_t)$ of the time series $(r_t)$, an alternative approach is to split the task into two parts: individual conditional volatility dynamics on one side, and conditional correlation dynamics on the other side. The most commonly used correlation process is the Dynamic Conditional Correlation (DCC) of \cite{engle2002}. In its BEKK form, the general DCC model is specified as: $r_t = \Sigma^{1/2}_t \eta_t$, with $\Sigma_t := \Eb[r_t r^\top_t | \Fc_{t-1}]$ positive definite, where:
\begin{equation*}
\Sigma_t = D_t R_t D_t, \; R_t =  Q^{\star-1/2}_t Q_t Q^{\star-1/2}_t, \;
Q_t = \Omega + \overset{q}{\underset{k=1}{\sum}} M_k Q_{t-k} M^\top_k + \overset{r}{\underset{l=1}{\sum}} W_l u_{t-l} u^\top_{t-l} W^\top_l,
\end{equation*}
where $D_t = \text{diag}\left(\sqrt{h_{t,1}},\sqrt{h_{t,2}},\ldots,\sqrt{h_{t,p}}\right)$ with $(h_{t,j})$ the univariate conditional variance dynamic of the $j$-th component, $u_t = \left(u_{t,1},\ldots,u_{t,p}\right)$ with $u_{t,j} = r_{t,j}/\sqrt{h_{t,j}}$, $Q_t = \left[q_{t,uv}\right]$, $Q^{\star}_t = \text{diag}\left(q_{t,11},q_{t,22},\ldots,q_{t,pp}\right)$. $\Fc_{t-1}$ is the sigma field generated by the past information of the $p$-dimensional process $(r_t)_{t \in \mathbb{Z}}$ until (but including) time $t-1$. The model is parameterized by some deterministic matrices $(M_k)_{k = 1,\ldots,q}$, $(W_l)_{l = 1,\ldots,r}$ and a positive definite $p \times p$ matrix $\Omega$. In the original DCC of \cite{engle2002}, the process $(Q_t)$ is defined by $Q_t = \Omega^\star +\sum^q_{k=1} B_k \odot Q_{t-k} +\sum^r_{l=1} A_l \odot u_{t-l} u^\top_{t-l}$, where the deterministic matrices $(B_k)_{k = 1,\ldots,q}$ and $(A_l)_{l = 1,\ldots,r}$ are positive semi-definite. Since the number of parameters of the latter models is of order $O(p^2)$, our applications are restricted to scalar $B_k$'s and $A_l$, say $b_k, a_l$, and we work with $q=r=1$. Moreover, we employ a correlation targeting strategy, that is we replace $\Omega^\star$ by $(1-a_1-b_1)\overline{Q}$ with $\overline{Q}$ the sample covariance of the standardized returns $u_t$. We specify a GARCH(1,1) model for each individual conditional variance dynamic $h_{tj}$. The estimation is performed by a standard two-step Gaussian QML for the MSCI and S\&P 100 portfolio. For the S\&P 500 portfolio, we employ the composite likelihood method proposed in~\cite{pakel2021} to fix the bias problem caused by the two-step quasi likelihood estimation and to make the likelihood decomposition plausible for large vectors. The composite likelihood method consists of averaging the second step likelihood (correlation part) over $2 \times 2$ correlation-based likelihoods. To be precise, using the $n$-sample of observations, the second step composite likelihood becomes:
\begin{equation*}
\Lb_{2n}(a_1,b_1) = \frac{1}{n} \overset{n}{\underset{t=1}{\sum}}\overset{L}{\underset{l=1}{\sum}} \Big[\log(|R^{(l)}_t|)+u^{(l)\top}_t\big(R^{(l)}_t\big)^{-1}u^{(l)}_t\Big],
\end{equation*}
where $R^{(l)}_t$ is the $2 \times 2$ correlation matrix with $l$ corresponding to a pre-specified pair of indices in $\{1,\ldots,p\}$, $u^{(l)}_t$ is a $2 \times 1$ sub-vector of the standardized residuals $u_t$, with $k$ selecting the pre-specified pair of indices in $\{1,\ldots,p\}$. There are several ways of selecting the pairs: every distinct pair of the assets, i.e, $u^{(1)}_t=(u_{t,1},u_{t,2})^\top,u^{(2)}_t=(u_{t,1},u_{t,3})^\top,\ldots,u^{p(p-1)/2}_t=(u_{t,p-1},u_{t,p})^\top$, so that $L = p(p-1)/2$ in $L_{2n}(a_1,b_1) $; or contiguous overlapping pairs, i.e., $u^{(1)}_t=(u_{t,1},u_{t,2})^\top,u^{(2)}_t=(u_{t,2},u_{t,3})^\top,\ldots,u^{(p-1)}_t=(u_{t,p-1},u_{t,p})^\top$ so that $L=p-1$. We use the contiguous overlapping pairs for computation gains, where the complexity is $O(p)$, in contrast to the $O(p^2)$ complexity if $L=p(p-2)/2$: this pair construction is the 2MSCLE method of \cite{pakel2021}. 

\subsection{Sparse approximate factor model}\label{safm}

Based on the same factor model decomposition $X_t = \Lambda F_i + \eps_i$, rather than assuming sparsity w.r.t. $\Lambda$, \cite{bailiao2016} assume a non-diagonal but sparse $\Psi$. Assuming bounded eigenvalues for $\Psi$ and under the constraint for identification $\Lambda^\top \Psi^{-1}\Lambda$ diagonal, \cite{bailiao2016} consider the joint estimation of $(\Lambda,\Psi)$ by Gaussian QML with adaptive-LASSO penalization w.r.t. $\Psi$: this is the sparse approximate factor model (SAFM) estimator. In all our applications, we employ the Matlab code of \cite{bailiao2016} which solves the SAFM by the EM algorithm: the code can be downloaded in \url{https://econweb.rutgers.edu/yl1114/papers/factor3/factor3.html}. The initial values $(\widehat{\Lambda}_{\text{init}},\widehat{\Psi}_{\text{init}})$ are set the Gaussian QML estimators of \cite{bai2012} with $\widehat{\Psi}_{\text{init}}$ diagonal. We implemented the adaptive LASSO where the selection of the tuning parameter $\gamma_n$ follows from same method as in Subsection \ref{CV} in the time series case. The macroeconomic data used in Section \ref{diffusion_index} can be retrieved from the aforementioned website. 

\end{document}